\documentclass[10pt,oneside,a4paper,english]{book}

\title{Many Models for Water Waves}
\author{Vincent Duchêne}
\date{ \today}

\newcommand{\commentedref}[1]{\vspace{-10pt}}

\usepackage[french,main=english]{babel}
\usepackage[utf8x]{inputenc} 
\usepackage[T1]{fontenc}   
\usepackage[ec]{aeguill}     
\usepackage{enumerate}

\usepackage{titling}  
\usepackage{epigraph}
\usepackage{booktabs}

\usepackage{amsmath, amsthm} 
\usepackage{amsfonts,amssymb}
\usepackage{mathtools}
\usepackage{bm} 
\usepackage{mathrsfs}
\usepackage{stmaryrd}
\usepackage[super]{nth}

\usepackage{float}          
\usepackage[justification=centering]{caption}
\usepackage{subcaption}

\usepackage[pdftex]{graphicx}        

\usepackage[usenames,dvipsnames]{xcolor}
\definecolor{dgreen}{rgb}{0,.5,0}

\usepackage{tikz}
\newcommand{\mytikzstyle}{%
	\tikzstyle{model}=[rectangle,draw,rounded corners=4pt,fill=gray!50,minimum width=10em,minimum height=3em]
	\tikzstyle{conjecture}=[->,dotted,thin,>=stealth,color=gray ]
	\tikzstyle{formal}=[->,dotted,very thick,>=stealth,color=gray ]
	\tikzstyle{proved}=[->,dotted,very thick,>=stealth ]
	\tikzstyle{obvious}=[->,thick,>= stealth ,rounded corners=5pt]
}

\newcommand{\movie}[3][width=\textwidth]{
{ 
		\vspace{0.5em}
		\begin{center}
			\small Overlapped initial and final time snapshots \normalsize 
			\includegraphics[#1]{./#2-0-#3.png}
		\end{center}
	}
}

\usepackage{chngcntr}

\counterwithout{footnote}{chapter} 
\numberwithin{equation}{section} 
\numberwithin{figure}{section}
\numberwithin{table}{section}

\usepackage[Sonny]{fncychap}
\usepackage{chngcntr}\counterwithout{section}{chapter}

\usepackage{fancyhdr}        

\newcommand{\logo}{{\sf \raisebox{\depth}{\rotatebox{180}{WW}}fWW}}
\renewcommand{\logo}{}

\fancypagestyle{plain}         
{
	\fancyhf{}
	\renewcommand{\headrulewidth}{0pt}
}

\fancypagestyle{empty}         
{
	\fancyhf{}
	\renewcommand{\headrulewidth}{0pt}
}

\fancypagestyle{intro}{
	\fancyfoot{}
	\fancyhead[C]{} %
	\fancyhead[R]{}  
	\fancyhead[L]{\scriptsize\thepage}
} 
\fancypagestyle{prologue}{%
	\fancyfoot{}
	\fancyhead[C]{\small\nouppercase{\scshape Prologue}} %
	\fancyhead[R]{\scriptsize\logo}  
	\fancyhead[L]{\scriptsize\thepage}
}

\fancypagestyle{mystyle}{%
	\fancyfoot{}
	\fancyhead[C]{\small\nouppercase{\scshape Many Models for Water Waves}} %
	\fancyhead[R]{\scriptsize\logo}  
	\fancyhead[L]{}
}

\fancypagestyle{appendix}{%
	\fancyfoot{}
	\fancyhead[C]{\small\nouppercase{\scshape Appendix}} %
	\fancyhead[R]{\scriptsize\logo}  
	\fancyhead[L]{}
}

\fancypagestyle{biblio}{%
	\fancyfoot{}
	\fancyhead[C]{\small\nouppercase{\scshape Bibliography}} %
	\fancyhead[R]{\scriptsize\logo}  
	\fancyhead[L]{}
}

\usepackage[final=true, pdfpagemode=FullScreen, colorlinks=true,
citecolor=blue,
linkcolor=red,
backref=page,
bookmarksdepth=3,
hypertexnames=true]{hyperref}       
\usepackage[hyperpageref]{backref} 
\makeatletter
\def\neutralisetitre{\def\chapter{\@ifstar\@gobble\@gobble}}
\def\deneutralisetitre{\def\chapter{\@startsection {section}{1}{\z@}%
		{-3.5ex \@plus -1ex \@minus -.2ex}%
		{2.3ex \@plus.2ex}%
		{\normalfont\huge\bfseries}}}
\makeatother

\renewcommand*{\backref}[1]{}
\renewcommand*{\backrefalt}[4]{%
	\ifcase #1%
	{}
	\or
	{\footnotesize(cit. p.~#2)}
	\else
	{\footnotesize(cit. pp.~#2)}
	\fi
}

\usepackage{bookmark}   
\usepackage[english]{cleveref}		

\usepackage{minitoc}
\mtcsetfeature{minitoc}{pagestyle}%
{\thispagestyle{empty}\pagestyle{empty}}

\newcounter{myeqncounter}
\newenvironment{myequation}{%
	\addtocounter{equation}{-1}
	\refstepcounter{myeqncounter}

	\equation%
}{%
	\endequation%
}

\newcounter{myfigcounter}
\newenvironment{myfigure}{%
	\refstepcounter{myfigcounter}

	\begin{figure}}
	{\end{figure}}

\newcounter{myFigcounter}
\newenvironment{myFigure}{%
	\refstepcounter{myFigcounter}

	\begin{figure}}
	{\end{figure}}

\newtheorem*{Theorem*}{Theorem}
\newtheorem*{Definition*}{Definition}
\newtheorem*{Proposition*}{Proposition}
\newtheorem*{Corollary*}{Corollary}
\newtheorem*{Lemma*}{Lemma}
\newtheorem*{Remark*}{Remark}
\newtheorem*{Assumption*}{Assumption}
\newtheorem*{Claim*}{Claim}
\newtheorem*{Conjecture*}{Conjecture}

\newcommand{\eps}{\varepsilon}
\newcommand{\sep}{\epsilon}

\newcommand{\eqdef}{\stackrel{\rm def}{=}}

\newcommand{\ie}{{\em i.e.}~}
\newcommand{\eg}{{\em e.g.}~}

\newcommand{\id}[1]{\left\vert_{_{#1}}\right.}

\newcommand{\SWSW}{{\tfrac{\cramped[\scriptscriptstyle]{\rm SW}}{\cramped[\scriptscriptstyle]{\rm SW}}}}

\DeclareMathOperator{\dd}{d\!}

\DeclareMathOperator{\rot}{{rot}}

\DeclareMathOperator{\sgn}{{sgn}}
\DeclareMathOperator{\sech}{{sech}}

\DeclarePairedDelimiterX\innerb[2]{\langle}{\rangle}{#1,#2}
\DeclarePairedDelimiterX\innerp[2]{(}{)}{#1,#2}

\newcommand{\RR}{\mathbb{R}}

\newcommand{\NN}{\mathbb{N}}

\newcommand{\mfp}{\mathfrak{p}}

\newcommand{\cG}{\mathcal{G}}

\newcommand{\cO}{\mathcal{O}}

\newcommand{\bx}{{\bm{x}}}

\newcommand{\bz}{{\bm{0}}}

\newcommand{\fL}{{\sf L}}

\DeclareMathAlphabet{\mathpzc}{OT1}{cmbrm}{m}{it}

\DeclareSymbolFont{dletters}{OML}{cmbrm}{m}{it}
\DeclareSymbolFont{bletters}{OML}{cmbrm}{b}{it}
\DeclareSymbolFont{oletters}{OT1}{pzc}{m}{n}

\DeclareMathSymbol{\dbc}{\mathalpha}{bletters}{`c}
\DeclareMathSymbol{\dbe}{\mathalpha}{bletters}{`e}
\DeclareMathSymbol{\dbr}{\mathalpha}{bletters}{`r}
\DeclareMathSymbol{\dbu}{\mathalpha}{bletters}{`u}
\DeclareMathSymbol{\dbv}{\mathalpha}{bletters}{`v}
\DeclareMathSymbol{\dbx}{\mathalpha}{bletters}{`x}
\DeclareMathSymbol{\dby}{\mathalpha}{bletters}{`y}
\DeclareMathSymbol{\dbR}{\mathalpha}{bletters}{`R}
\DeclareMathSymbol{\dbK}{\mathalpha}{bletters}{`K}
\DeclareMathSymbol{\dbU}{\mathalpha}{bletters}{`U}
\DeclareMathSymbol{\dmA}{\mathalpha}{oletters}{`A}
\DeclareMathSymbol{\dmC}{\mathalpha}{oletters}{`C}
\DeclareMathSymbol{\dmG}{\mathalpha}{oletters}{`G}
\DeclareMathSymbol{\dmg}{\mathalpha}{oletters}{`g}
\DeclareMathSymbol{\dmH}{\mathalpha}{oletters}{`H}
\DeclareMathSymbol{\dmI}{\mathalpha}{oletters}{`I}
\DeclareMathSymbol{\dmL}{\mathalpha}{oletters}{`L}
\DeclareMathSymbol{\dmM}{\mathalpha}{oletters}{`M}
\DeclareMathSymbol{\dmZ}{\mathalpha}{oletters}{`Z}

\DeclareMathSymbol{\dA}{\mathalpha}{dletters}{`A}
\DeclareMathSymbol{\dC}{\mathalpha}{dletters}{`C}
\DeclareMathSymbol{\dH}{\mathalpha}{dletters}{`d}
\DeclareMathSymbol{\ddH}{\mathalpha}{dletters}{`H}
\DeclareMathSymbol{\dI}{\mathalpha}{dletters}{`I}
\DeclareMathSymbol{\ddL}{\mathalpha}{dletters}{`L}
\DeclareMathSymbol{\dL}{\mathalpha}{dletters}{21}
\DeclareMathSymbol{\dM}{\mathalpha}{dletters}{`M}
\DeclareMathSymbol{\dP}{\mathalpha}{dletters}{`P}
\DeclareMathSymbol{\dS}{\mathalpha}{dletters}{`S}
\DeclareMathSymbol{\dT}{\mathalpha}{dletters}{`T}
\DeclareMathSymbol{\dU}{\mathalpha}{dletters}{`U}

\DeclareMathSymbol{\da}{\mathalpha}{dletters}{`a}
\DeclareMathSymbol{\db}{\mathalpha}{dletters}{`b}
\DeclareMathSymbol{\dc}{\mathalpha}{dletters}{`c}
\DeclareMathSymbol{\ddd}{\mathalpha}{dletters}{`d}
\DeclareMathSymbol{\df}{\mathalpha}{dletters}{`f}
\DeclareMathSymbol{\dg}{\mathalpha}{dletters}{`g}
\DeclareMathSymbol{\ddh}{\mathalpha}{dletters}{`h}
\DeclareMathSymbol{\dm}{\mathalpha}{dletters}{`m}
\DeclareMathSymbol{\ddp}{\mathalpha}{dletters}{`p}
\DeclareMathSymbol{\dq}{\mathalpha}{dletters}{`q}
\DeclareMathSymbol{\dr}{\mathalpha}{dletters}{`r}
\DeclareMathSymbol{\ds}{\mathalpha}{dletters}{`s}
\DeclareMathSymbol{\dt}{\mathalpha}{dletters}{`t}
\DeclareMathSymbol{\du}{\mathalpha}{dletters}{`u}
\DeclareMathSymbol{\dv}{\mathalpha}{dletters}{`v}
\DeclareMathSymbol{\dw}{\mathalpha}{dletters}{`w}
\DeclareMathSymbol{\dx}{\mathalpha}{dletters}{`x}
\DeclareMathSymbol{\dy}{\mathalpha}{dletters}{`y}
\DeclareMathSymbol{\dz}{\mathalpha}{dletters}{`z}
\DeclareMathSymbol{\dR}{\mathalpha}{dletters}{`R}
\DeclareMathSymbol{\dX}{\mathalpha}{dletters}{`X}

\DeclareMathSymbol{\dGamma}{\mathalpha}{dletters}{0}
\DeclareMathSymbol{\dDelta}{\mathalpha}{dletters}{1}
\DeclareMathSymbol{\dTheta}{\mathalpha}{dletters}{2}
\DeclareMathSymbol{\dLambda}{\mathalpha}{dletters}{3}
\DeclareMathSymbol{\dXi}{\mathalpha}{dletters}{4}
\DeclareMathSymbol{\dPi}{\mathalpha}{dletters}{5}
\DeclareMathSymbol{\dSigma}{\mathalpha}{dletters}{6}
\DeclareMathSymbol{\dUpsilon}{\mathalpha}{dletters}{7}
\DeclareMathSymbol{\dPhi}{\mathalpha}{dletters}{8}
\DeclareMathSymbol{\dPsi}{\mathalpha}{dletters}{9}
\DeclareMathSymbol{\dOmega}{\mathalpha}{dletters}{10}
\DeclareMathSymbol{\dalpha}{\mathalpha}{dletters}{11}
\DeclareMathSymbol{\dbeta}{\mathalpha}{dletters}{12}
\DeclareMathSymbol{\dgamma}{\mathalpha}{dletters}{13}
\DeclareMathSymbol{\ddelta}{\mathalpha}{dletters}{14}
\DeclareMathSymbol{\depsilon}{\mathalpha}{dletters}{15}
\DeclareMathSymbol{\dzeta}{\mathalpha}{dletters}{16}
\DeclareMathSymbol{\deta}{\mathalpha}{dletters}{17}
\DeclareMathSymbol{\dtheta}{\mathalpha}{dletters}{18}
\DeclareMathSymbol{\diota}{\mathalpha}{dletters}{19}
\DeclareMathSymbol{\dkappa}{\mathalpha}{dletters}{20}
\DeclareMathSymbol{\dlambda}{\mathalpha}{dletters}{21}
\DeclareMathSymbol{\dmu}{\mathalpha}{dletters}{22}
\DeclareMathSymbol{\dnu}{\mathalpha}{dletters}{23}
\DeclareMathSymbol{\dxi}{\mathalpha}{dletters}{24}
\DeclareMathSymbol{\dpi}{\mathalpha}{dletters}{25}
\DeclareMathSymbol{\drho}{\mathalpha}{dletters}{26}
\DeclareMathSymbol{\dsigma}{\mathalpha}{dletters}{27}
\DeclareMathSymbol{\dtau}{\mathalpha}{dletters}{28}
\DeclareMathSymbol{\dupsilon}{\mathalpha}{dletters}{29}
\DeclareMathSymbol{\dphi}{\mathalpha}{dletters}{30}
\DeclareMathSymbol{\dchi}{\mathalpha}{dletters}{31}
\DeclareMathSymbol{\dpsi}{\mathalpha}{dletters}{32}
\DeclareMathSymbol{\domega}{\mathalpha}{dletters}{33}
\DeclareMathSymbol{\dvarepsilon}{\mathalpha}{dletters}{34}
\DeclareMathSymbol{\dvartheta}{\mathalpha}{dletters}{35}
\DeclareMathSymbol{\dvarpi}{\mathalpha}{dletters}{36}
\DeclareMathSymbol{\dvarrho}{\mathalpha}{dletters}{37}
\DeclareMathSymbol{\dvarsigma}{\mathalpha}{dletters}{38}
\DeclareMathSymbol{\dvarphi}{\mathalpha}{dletters}{39}

\DeclareMathSymbol{\dbxi}{\mathalpha}{bletters}{24}

\begin{document}

\begin{titlepage}

\centering

\vspace{2.5cm}

{\Huge \sf \raisebox{\depth}{\rotatebox{180}{W}}any \raisebox{\depth}{\rotatebox{180}{W}}odels  for Water Waves \par}
\vspace{1cm}
{\Large	\sf A unified theoretical approach \par}
\vspace{2cm}
{\Large Vincent \textsc{Duchêne} \par}
\vspace{.5cm}
{\normalsize Centre national de la recherche scientifique (CNRS)\\
 Institut de Recherche Mathématique de Rennes (IRMAR), Université de Rennes 1 \par}
\vspace{.5cm}
\href{mailto:vincent.duchene@univ-rennes1.fr}{vincent.duchene@univ-rennes1.fr}
\vfill

\vfill

\raggedleft

\setcounter{page}{0} 
\end{titlepage}

\pagenumbering{alph} 
\pagestyle{empty}


\addtocontents{toc}{\vspace{-50pt}} 
\addtocontents{lof}{ \vspace{-50pt}
Figures listed in {\bf bold} expose interplays between models presented in this work.
 \vspace{20pt}
} 

\dominitoc
\dominilof


\paragraph*{Caveat lector}

This document is an announcement and preview\footnote{The present document contains the prologue, slightly edited forewords of all chapters, and the bibliography.} of a memoir whose full version is available on the \href{https://www.ams.org/open-math-notes/omn-view-listing?listingId=111309}{Open Math Notes} repository of the American Mathematical Society. In this memoir, I try to provide a fairly comprehensive picture of (mostly shallow water) asymptotic models for water waves. The work and presentation is heavily inspired by the book of D. Lannes \cite{Lannes}, yet extends the discussion into several directions, notably high order and fully dispersive models, and internal/interfacial waves. 

 I plan to update this memoir from time to time when novel material fitting in the picture will arise. Please do not hesitate to contact me when you notice typos or mistakes, or if you have any question, comment or query, using the email address provided on the front page.
 
\bigskip

In the memoir one derives, discusses, and justifies as much as possible a large class of models describing in an approximate manner the propagation of waves at the surface of water, at the interface between two homogeneous fluids, or in the bulk of a continuously density-stratified fluid.
In the considered idealized frameworks, these waves propagate from an initial perturbation of the rest state under the influence of gravity forces.
Let me unveil a little bit of the material contained in the memoir in order to warn the potentially disappointed reader.
\begin{itemize}
	\item The motivation is theoretical, in the sense that practical direct use of the results is not the main objective. The problem of the propagation of water waves is one example of partial differential equations which may be written under a compact formulation but forecasts a fascinating variety of phenomena, while enjoying a rich mathematical structure. It is hence a formidable toy on which one can apply advanced tools of modelization. Yet it is impossible not to have in mind that practical applications are just a few steps away, and many choices in the modelization procedure are grounded on applicative views, for instance robustness of the models or easy numerical implementation.
	\item The \og master\fg{} equations, that is the system of equations from which all subsequent simplified models are derived, already incorporates many idealizations. To name a few, earth curvature, the Coriolis force, wind forcing, any dissipative effect and---most of the time---surface tension are neglected. In the \og water waves\fg{} case one considers homogeneous fluids and potential flows. Moreover, the analysis is restricted to laminar (\ie regular) rather than turbulent flows.
	\item While the equations at stake are of dispersive nature, little or none of the advanced tools on dispersive equations is used, and we barely report on the latest mathematical developments involving paradifferential calculus, normal forms, KAM theory, {\em etc.} The main mathematical tools that are put to use are rather old but robust: on one hand the elliptic theory to derive models from approximate solutions to a Laplace problem; and on the other hand the energy method to justify rigorously the resulting evolution equations (being of quasilinear hyperbolic nature). The heart of the matter consists in using these tools in a sufficiently refined manner so as to offer error estimates uniform with respect to the relevant parameters at stake. 
	\item Given their number and diversity, it is impossible to present all relevant models based on the water waves system, even restricting to a specific asymptotic regime (the shallow water regime in our case). The memoir focuses on models which preserve as much as possible the structure (and in particular the Hamiltonian formulation) of the master equations, as well as mathematical properties (typically the well-posedness of the initial-value problem). That such models are often historical and among the most studied is not, to my opinion, a coincidence. Hence most of this work is dedicated to fairly standard models in oceanography. The aim of this document is to present such models together with more recent ones in a unified framework, and to address the state of their rigorous justification.
\end{itemize}

\tableofcontents
\listoffigures



\frontmatter
\pagestyle{prologue} 
\renewcommand{\thesection}{\roman{section}}
\renewcommand{\headrulewidth}{.5pt}

\chapter*{Prologue}\label{C.prologue}
\addstarredchapter{Prologue}
\phantomsection 

\epigraph{Y'a tant de vagues, et tant d'idées qu'on n'arrive plus à décider le faux du vrai}{--- \textsc{Michel Berger}, \textit{Le paradis blanc}}
\minitoc
\section*{Foreword}

In this monograph we aim at describing the evolution in time of a body of fluid---typically water.
Of course the features of the dynamics depend greatly on the framework, and in particular on the scales involved.
As a rule of thumb, we will be motivated by the description of the motion of the surface of water as seen by a human eye. These are often referred to as surface gravity waves, or simply {\em water waves}. As any wanderer knows, despite the restrictive framework, water waves are still remarkably diverse, and this is what makes them a fascinating subject of study.\footnote{To quote Feynman during his well-known Lectures on Physics (Vol. I, Ch. 51: Waves): ``Now, the next waves of interest, that are easily seen by everyone and which are usually used as an example of waves in elementary courses, are water waves. As we shall soon see, they are the worst possible example, because they are in no respects like sound and light; they have all the complications that waves can have.'' }
 In order to get a grasp at the behavior of water waves in a given situation, one typically uses simplified models. Below we give examples of a few such models\footnote{We do not attempt at exhaustiveness. The relentless reader will find more in the present document and much more in the literature, using for instance \cite{MadsenFuhrman10,Lannes,Saut13,BridgesGrovesNicholls16,Lannes20} as starting points. }  which appeared in the early literature,\footnote{The interested reader will find in \cite{Darrigol03} a detailed historical account on the early studies on water waves.} with the aim at emphasizing the diversity of possible waves and the hope of giving an insight at the possible mechanisms involved in the full picture. The models described further on in this work are refinements of such models.

 There are many ways to formally derive the models presented below. Considering the Saint-Venant system for instance, a typical way consists in integrating the horizontal velocity over the fluid layer and invoking a closure formula, based on physical principles such as energy conservation. One can also use some {\em ad hoc} hypotheses, such as columnar motion and hydrostatic approximation. Or a loose assumption that derivatives of a function are smaller than the function itself. Our strategy, called {\em asymptotic modeling}, is akin to the latter one, and provides a justification of the former ones, with quantitative estimates of the inaccuracies. We start with the so-called {\em full Euler system} (or more precisely, for models in this Prologue, the {\em water waves system}) whose solutions are regarded as ``exact'' (although, admittedly, the derivation of the equations relies on many oversimplifications). Using the typical scales of the flow, we can extract dimensionless parameters describing the strength of the main mechanisms involved. The asymptotic models are obtained through a description of the operators involved in the water waves system using assumptions on the size of these parameters, which will be called the {\em asymptotic regime}.

 The complete rigorous justification of models in a given asymptotic regime typically proceeds in two steps. First we prove that sufficiently regular solutions to the water waves system satisfy the equations of the model---or the other way around---up to a small remainder term, measured by the size of the dimensionless parameters and data in a prescribed metric space; this is called {\em consistency}.  Anticipating with future notations and results, we find that the water waves system is consistent with the acoustic wave equation~\eqref{eq.intro..Acoustic} with precision $\cO(\mu+\eps)$, with the linearized (Airy) equations~\eqref{eq.intro..WW0} with precision $\cO(\eps)$, with the Saint-Venant system~\eqref{eq.intro..SV} with precision $\cO(\mu)$, with all the Boussinesq systems~\eqref{eq.intro..Boussinesq} with precision $\cO(\mu^2+\mu\eps)$, {\em etc.} This is however not sufficient, and there remains to prove that for a large class of sufficiently regular initial data (typically a neighborhood of the rest state in the aforementioned metric space), there exist unique solutions to both the water waves system and the asymptotic model, and that the two remain close on the relevant timescale. Following Lannes~\cite{Lannes}, we call the former property (uniform) {\em well-posedness}, and the latter {\em convergence}. 
 
 An important portion of this monograph is dedicated to the rigorous justification in the above sense---together with the study of a few basic properties---of standard and less-standard models for the propagation of surface, interfacial and internal gravity waves.

\newpage

\section{The linear acoustic wave equation}\label{S.intro..acoustic}

Arguably the simplest (partial differential) equation describing the motion of water waves, already put forward by Lagrange~\cite{Lagrange81}, is the following:
\begin{myequation}\label{eq.intro..Acoustic}
	\partial_\dt^2\dzeta=\dg\dH\, \Delta_\bx \dzeta\,.
\end{myequation}
Here $\dzeta$ represents the deformation of the free surface, in the sense that the surface of the body of water at time $\dt$ is parameterized as
\[\dGamma_{\rm top}= \{(\dbx,\dz)\in\RR^{d+1}\ : \  \dz=\dzeta(\dt,\dbx)\}.\]
Hence the function $\dzeta$ depends on time, $\dt$, and horizontal space variable, $\dbx$. For simplicity we assume that the horizontal variable lies in the full space $\RR^d$. The constant $\dg$ denotes the gravity acceleration and $\dH$ is the depth of the layer. \Cref{eq.intro..Acoustic} is called the {\bf{\em linear acoustic wave equation}} as it governs the propagation of infinitesimally small acoustic waves through a material medium. It is only a coincidence that it also describes---very roughly, remember Feynman's quote---water waves. In fact the above equation describes infinitely small and infinitely long water waves.

In the special case of horizontal dimension $d=1$,\footnote{The one-dimensional framework $d=1$ is relevant for instance for waves propagating along a narrow channel. } the solution of the initial-value problem is easily found as
\[\dzeta(\dt,\dx)=\frac12\Big(\dzeta(0,\dx+ \dc_0\dt)+\dzeta(0,\dx- \dc_0\dt)\Big)+\frac1{2\dc_0}\int_{\dx- \dc_0\dt}^{\dx+ \dc_0\dt}\partial_t\zeta(0,\dy)\dd\dy.\]
with $\dc_0=\sqrt{\dg\dH}$. Hence the wave decomposes into the superposition of a right-going and a left-going components, both translating with velocity $\dc_0$.
This is shown in \Cref{F.intro..1DWW0} where the evolution of the surface deformation when taken initially as Gaussians (with zero initial velocities) according to \cref{eq.intro..Acoustic} and \cref{eq.intro..WW0} in dimension $d=1$ is represented.
\begin{myfigure}[htb]
	\begin{subfigure}{.5\textwidth}
		\movie{1DAiry}{60}
		\caption{Linear acoustic wave equation, \cref{eq.intro..Acoustic}}
	\end{subfigure}%
	\begin{subfigure}{.5\textwidth}
		\movie{1DWW0}{60}
		\caption{Linearized water waves system, \cref{eq.intro..WW0}}
	\end{subfigure}
	\caption[Disintegration of heap of water, according to the linear acoustic wave equation and linearized water waves system, in dimension $d=1$.]{Disintegration of Gaussian initial data, $\dzeta(\dt=0,\dx)=0.01\exp(-(0.1\,\dx)^2)$ (left) and $\dzeta(\dt=0,\dx)=0.01\exp(-\dx^2)$ (right), with zero initial velocities. $\dg=9.81\, {\rm m}.{\rm s}^{-2}$, $\dH=1\,{\rm m}$. }
	\label{F.intro..1DWW0}
\end{myfigure}

In dimension $d=2$, the solution is less explicit, but a formula can still be written---at least for sufficiently regular initial data---with the use of Green's function (we could also use Fourier representation as in the next section):
\[\dzeta(\dt,\dbx)=\frac1{2\pi\dc_0}\int_{|\dby-\dbx|\leq \dc_0 \dt} \frac{\partial_\dt\dzeta(0,\dby)}{\sqrt{(\dc_0\dt)^2-|\dby-\dbx|^2}}\dd\dby
+\frac1{2\pi\dc_0}\partial_\dt\int_{|\dby-\dbx|\leq \dc_0 \dt} \frac{\dzeta(0,\dby)}{\sqrt{(\dc_0\dt)^2-|\dby-\dbx|^2}}\dd\dby.\]
We can observe that the solution satisfies causality (but not Huygens' principle): waves must be  given enough time to propagate between two specified points. Again, $\dc_0$ is a good measure of the (scalar) velocity of waves according to \cref{eq.intro..Acoustic}. Less obvious is the fact that for sufficiently smooth and decaying initial data, the amplitude of the solution decays for large time as $(\dc_0\dt)^{-1/2}$.
\Cref{F.intro..2DWW0} represents the evolution of the surface deformation when taken initially as Gaussians (with zero initial velocities) according to \cref{eq.intro..Acoustic} and \cref{eq.intro..WW0}, in dimension $d=2$.

\begin{myfigure}[!ht]
	\begin{subfigure}{.5\textwidth}
		\movie{2DAiry}{60}
		\caption{Linear acoustic wave equation, \cref{eq.intro..Acoustic}}
	\end{subfigure}%
	\begin{subfigure}{.5\textwidth}
		\movie{2DWW0}{60}
		\caption{Linearized water waves system, \cref{eq.intro..WW0}}
	\end{subfigure}
	\caption[Disintegration of a heap of water, according to the linear acoustic wave equation and linearized water waves system, in dimension $d=2$]{Disintegration of Gaussian initial data, with zero initial velocities.\\ $\dg=9.81\, {\rm m}.{\rm s}^{-2}$, $\dH=1\,{\rm m}$. The bottom plot represents the  solution on $\{(\dx,\dy) \ : \ \dy=0\}$.}
	\label{F.intro..2DWW0}
\end{myfigure}
\section{The linearized (Airy) water waves equations}\label{S.intro..Airy}

The following equations describe the propagation of infinitesimally small waves without the long wave assumption of the previous section: it is the linearized system about the rest-state solution to the water waves equations, whose solutions shall be considered as ``exact'', and which is introduced in \Cref{C.master}. Consider the {\bf{\em  linearized water waves equations}} as
\begin{myequation}\label{eq.intro..WW0}
	\left\{\begin{array}{l}
		\partial_\dt\dzeta-\dmG_0\dpsi=0,\\[1ex]
		\partial_\dt\dpsi +\dg\dzeta=0
	\end{array}\right.
\end{myequation}
where $\dmG_0=|D|\tanh(\dH|D|)$ is the Fourier multiplier operator defined on sufficiently regular solutions by
\[\forall\dbxi\in\RR^d, \qquad \widehat{\dmG_0\dpsi}(\dbxi)= |\dbxi|\tanh(\dH|\dbxi|)\widehat\dpsi(\dbxi).\]
Here, $\dg$, $\dH$ and $\dzeta$ are as above and $\dpsi$ represents the trace of the velocity potential at the surface.
\Cref{eq.intro..WW0} is a system of linear constant-coefficient equations of the form
\[\partial_\dt \begin{pmatrix}\dzeta\\ \dpsi
\end{pmatrix} =\fL(D)\begin{pmatrix}\dzeta\\ \dpsi
\end{pmatrix}\]
where $\fL(D)$ is a matrix with Fourier multiplier coefficients.

Formally taking the limit $\dH\to0$, we may replace $\tanh(\dH|\dbxi|)$ with $\dH|\dbxi|$ in $\dmG_0$, and then we recover the acoustic wave equation, \cref{eq.intro..Acoustic}. In fact using \cref{eq.intro..WW0} instead of \cref{eq.intro..Acoustic} in the left side of \Cref{F.intro..1DWW0} and \Cref{F.intro..2DWW0} yields a very similar outcome; such is not the case for the narrower initial data used for right sides.

\paragraph*{Modal analysis}
Plane waves of the form
$(\dzeta,\dpsi)=(\dzeta_0 e^{i(\dbxi\cdot\dbx-\domega \dt)},\dpsi_0e^{i(\dbxi\cdot\dbx-\domega \dt)})$
are solutions to \cref{eq.intro..WW0} provided that $i\domega\dpsi_0=\dg\dzeta_0$ and the {\em dispersion relation} holds~\cite{Kelland40,Airy45}:
\[\domega(\dbxi)^2=\dg|\dbxi|\tanh(\dH|\dbxi|).\]
In other words, we can explicitly solve the equation in the Fourier space:
\[\begin{pmatrix}\widehat\dzeta(\dt,\dbxi) \\ \widehat\dpsi(\dt,\dbxi)
\end{pmatrix} =\exp(\fL(\dbxi) \dt) \begin{pmatrix}\widehat\dzeta(0,\dbxi) \\ \widehat\dpsi(0,\dbxi)
\end{pmatrix}  = \begin{pmatrix}
\cos(|\domega(\dbxi)| \dt) &  \frac{|\domega(\dbxi)|}{\dg}\sin(|\domega(\dbxi)| \dt)\\
-\frac{\dg}{|\domega(\dbxi)|}\sin(|\domega(\dbxi)| \dt) &  \cos(|\domega(\dbxi)| \dt)
\end{pmatrix}\begin{pmatrix}\widehat\dzeta(0,\dbxi) \\ \widehat\dpsi(0,\dbxi)
\end{pmatrix}  .\]
For such plane wave solutions, $\domega$ is called the (angular) {\em frequency}, $\dbxi$ the (angular) {\em wave vector} (wavenumber if $d=1$), and $|\dbxi|$ the (angular) {\em wavenumber}.
{\em Phase velocities} describe the velocity in a given direction of a plane wave with wave vector $\dbxi$, and satisfy
\[ {\dbc}_p\cdot \dbxi = \domega(\dbxi).\]
The {\em group velocity} represents the traveling velocity of a wave packet about wave vector $\dbxi$, and is given by
\[\dbc_g = \nabla_{\dbxi} (\domega(\dbxi)).\]
Misusing these definitions, we shall also refer to
\[\dc_p =  \frac{|\domega(\dbxi)|}{|\dbxi|} = \sqrt{\dg\dH}\Big(\frac{\tanh(\dH|\dbxi|)}{\dH|\dbxi|}\Big)^{1/2}\]
as the phase velocity, and to
\[\dc_g=|\dbc_g|=\sqrt{\dg\dH}\left(\frac12\Big(\frac{\tanh(\dH|\dbxi|)}{\dH|\dbxi|}\Big)^{1/2} +\frac{\sech^2(\dH |\dbxi|)}2\Big(\frac{\dH|\dbxi|}{\tanh(\dH|\dbxi|)}\Big)^{1/2}  \right)\]
as the group velocity. They are represented in \Cref{F.intro..phase_and_group_velocities}. That the phase velocity is different (and greater) than the group velocity manifests the essential feature of the (linearized) water waves equations as being dispersive. Notice however that for small-normed wave vectors, $\dH|\dbxi|\ll 1$, both velocities converge to $\dc_0=\sqrt{\dg\dH}$, the velocity of (non-dispersive) infinitely long waves. In the opposite direction, for $\dH|\dbxi|\gg 1$, we have $\dc_g\sim\frac12\dc_p\sim \frac{\sqrt{\dg}}{2|\dbxi|^{1/2}}$.

\begin{myfigure}[!ht]
	\begin{center}
		\includegraphics[width=.6\textwidth]{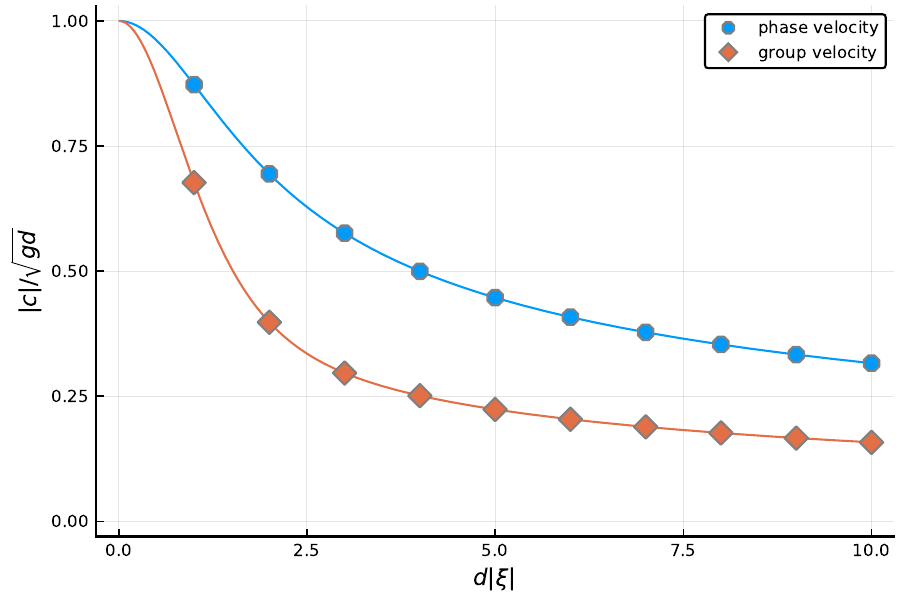}
	\end{center}
	\caption[Phase and group velocities of the linearized water waves system]{Phase and group velocities of the linearized water waves system.}
	\label{F.intro..phase_and_group_velocities}
\end{myfigure}

\paragraph*{Large-time behavior}
We can infer the large-time behavior of the solution, at least in dimension $d=1$, through the stationary phase theorem on oscillatory integrals; see \eg~\cite{Stein93}. Indeed, for any $\dc\in\RR$, and initial data such that $(\widehat\dzeta(0,\cdot) ,|\domega|(\cdot)\widehat\dpsi(0,\cdot))\in L^1(\RR)^2$, we have from the above
\[\dzeta(\dt,\dc \dt) = \frac1{4\pi} \int_{\RR} e^{i(\dc \dxi-\domega(\dxi) )\dt} \left(\widehat\dzeta(0,\dxi)+i \frac{\domega(\dxi) }{\dg}\widehat\dpsi(0,\dxi)\right)+ e^{i((\dc \dxi +\domega(\dxi) )\dt} \left(\widehat\dzeta(0,\dxi)-i \frac{\domega(\dxi) }{\dg}\widehat\dpsi(0,\dxi)\right) \dd\dxi\]
where we denote $\domega(\xi)=\sgn(\xi)(\dg|\xi|\tanh(\dH|\dxi|))^{1/2}$, and use a standard convention for the Fourier transform. We deduce that the following holds for sufficiently decaying and regular initial data.
\begin{enumerate}
	\item For any $\dc\in(-\infty,-\sqrt{\dg\dH})\cup (\sqrt{\dg\dH},+\infty)$, one has for any $n\in\NN$,
	\[ |\dzeta|(\dt,\dc \dt) =\cO(\dt^{-n}) .\]
	\item For any $\dc\in (-\sqrt{\dg\dH},\sqrt{\dg\dH})\setminus\{0\}$, one has
	\[ |\dzeta|(\dt,\dc \dt) \sim_{t\to\infty}  \frac1{4\pi}(2!)^{1/2}\Gamma(\tfrac32) |\dA(\dxi_\dc)|\big( |\domega''(\dxi_\dc)|\, \dt\big)^{-\frac12} \]
	where $\dxi_\dc$ is defined by the relation $\dc=\domega'(\dxi_\dc)$ and $\dA(\dxi_\dc)\eqdef\widehat\dzeta(0,\dxi_\dc)+\sgn(\dc)i \frac{\domega(\dxi_c) }{\dg}\widehat\dpsi(0,\dxi_\dc)$; unless $\dA(\dxi_\dc)=0$ in which case the decay is at least $\cO(\dt^{-1})$.
	\item If $\dc\in\{-\sqrt{\dg\dH},\sqrt{\dg\dH}\}$
	, one has
	\[ |\dzeta|(\dt,\dc \dt) \sim_{t\to\infty}  \frac1{4\pi}(3!)^{\frac13}\Gamma( \tfrac43) |\dA(0)|\big( \dH^2\sqrt{\dg\dH}\, \dt\big)^{-\frac13} \approx  \da\, ((\dH^2/\dL^2) \sqrt{\dg \dH}/\dL \dt )^{-\frac13}, \]
	with $ \dA(0)\eqdef \lim_{\dxi\to 0}\dA(\dxi)$ (notice we require regularity only on $\dxi\widehat\dpsi(0,\dxi)$); unless $\dA(0)=0$, in which case the decay is at least $\cO(\dt^{-\frac23})$. The last approximation is meant in a loose sense, where we set $\dA(0)\approx \da\dL$. This allows to hint at the timescale for which dispersive mechanisms have a bearing on the behavior of the flow, which is large compared with the time period of long waves, $\dT\eqdef \dL/\sqrt{\dg \dH}$, when $\dH^2/\dL^2\ll 1$.
\end{enumerate}
Above, $\Gamma$ is the Euler Gamma function: $\Gamma(s)\eqdef \int_0^{+\infty} \tau^{s-1}e^{-\tau}\dd \tau$.
A loose interpretation of the above is that for large time, the dominant part of the wave which will remain visible is the large wavelength component, traveling at velocity $|\dc|\approx \dc_0=\sqrt{\dg\dH}$.


\newpage

\section{The Saint-Venant system}\label{S.intro..Saint-Venant}

Our first nonlinear model is the so-called shallow water, or {\bf {\em Saint-Venant system}} system~\cite{Saint-Venant71}:
\begin{myequation}\label{eq.intro..SV}
	\left\{
	\begin{array}{l}
		\partial_\dt \dzeta+\nabla\cdot(\ddh\dbu)=0,\\[1ex] 
		\partial_\dt\dbu+\dg\nabla\dzeta+ (\dbu\cdot\nabla)\dbu=\bz,
	\end{array}
	\right.
\end{myequation}
where $\ddh\eqdef\dH+\dzeta$ represents the water depth, and $\dbu$ a horizontal velocity (it can be the layer-averaged horizontal velocity, velocity at a certain depth, or $\nabla\dpsi$).
Pursuing the analogy of \Cref{S.intro..acoustic}, one can notice that the Saint-Venant system is equivalent to the isentropic, {\em compressible} Euler equation for ideal gases with the pressure law $\ddp(\drho)\propto \drho^2$ (identifying $\drho$ with $\ddh$).

System~\eqref{eq.intro..SV}  is hence a prototype of quasilinear hyperbolic systems. Hyperbolicity amounts to the non-cavitation assumption, that is restricting data to  $\{(\dzeta,\dbu)\ : \ \dH+\dzeta>0\}$.\footnote{\label{f.intro:SV}Sufficiently regular solutions with initial data in the hyperbolicity domain cannot leave the hyperbolicity domain due to first equation (mass conservation) in \cref{eq.intro..SV}. Indeed, denoting  $\ddh_{\dbx_\star}(\dt)=\ddh(\dt,\dbx_{\dx_\star}(\dt))$  where $\dbx_{\dbx_\star}(\dt)$ is defined by the final condition $\dbx_{\dbx_\star}(\dt_\star)=\dbx_\star$ and the ordinary differential equation $ \dbx_{\dbx_\star}'(\dt)=\dbu(\dt,\dbx_{\dbx_\star}(\dt))$ for $\dt\in[0,\dt_\star]$, we find  $\ddh(\dt_\star,\dbx_\star)=\ddh_{\dbx_\star}(\dt_\star)=\ddh(0,\dbx_{\dbx_\star}(0))\exp(-\int_0^{\dt_\star} (\nabla\cdot\dbu)(\dt,\dbx_{\dbx_\star}(\dt))\dd\dt)>0$.} Indeed, the system in dimension $d=2$ reads
\[\partial_\dt \begin{pmatrix}
\dzeta \\ \du_\dx \\ \du_\dy
\end{pmatrix}+\begin{pmatrix}
\du_\dx & \ddh & 0 \\ \dg& \du_\dx & 0\\ 0&0&\du_\dx
\end{pmatrix}\partial_\dx\begin{pmatrix}
\dzeta \\ \du_\dx \\ \du_\dy
\end{pmatrix}+\begin{pmatrix}
\du_\dy &0& \ddh \\ 0& \du_\dy & 0\\ \dg &0&\du_\dy
\end{pmatrix}\partial_\dy\begin{pmatrix}
\dzeta \\ \du_\dx \\ \du_\dy
\end{pmatrix}=\begin{pmatrix}
0 \\ 0 \\ 0
\end{pmatrix}\]
and the eigenvalues of the associated symbol (see \eg \cite{Metivier08}) are
\[ \dbu\cdot\dbxi \quad \text{ and } \quad  \dbu\cdot\dbxi  \pm \sqrt{\dg\ddh} |\dbxi|\,.\]
Notice here again the ``sound speed'' of long surface gravity waves as being $\dc_0=\sqrt{\dg\dH}$.

In dimension $d=1$, as any quasilinear system of two scalar balance laws, \cref{eq.intro..SV} enjoys a basis of Riemann invariants. The Riemann invariants are explicit in this case: setting $\dr_\pm= \du\pm 2\sqrt{\dg\ddh}$, the system~\eqref{eq.intro..SV} is equivalent to
\begin{myequation}\label{eq.intro..SV-R}
	\left\{
	\begin{array}{l}
		\partial_\dt \dr_++\frac{3\dr_++\dr_-}4\partial_\dx \dr_+=0,\\[1ex]
		\partial_\dt \dr_-+\frac{3\dr_-+\dr_+}4\partial_\dx \dr_-=0.
	\end{array}
	\right.
\end{myequation}
Notice that $\frac{3r_++r_-}4= \du+\sqrt{\dg\ddh}$ and $\frac{3r_-+r_+}4= \du-\sqrt{\dg\ddh}$, consistently with the hyperbolicity discussion.
The diagonal formulation, \cref{eq.intro..SV-R}, allows to construct {\em simple waves}, \ie solutions of the form
\[(\dr_+,\dr_-)=\dbR(\dtheta(t,x))\]
where $\theta$ is a scalar function. For instance, any sufficiently regular solution to \cref{eq.intro..SV-R} with initial data satisfying $ \du\id{\dt=0}=2\sqrt{\dg\ddh\id{\dt=0}}-2\sqrt{\dg\dH}$, the second equation yields $\dr_-\equiv -2\sqrt{\dg\dH}$ for all times, from which we deduce $\dr_+= 2\sqrt{\dg\dH}+2 \du$, where $\du(\dt,\dx)$ satisfies the {\bf {\em inviscid Burgers}} (or Hopf) equation
\begin{myequation}\label{eq.intro..Hopf}
	\partial_\dt \du+\big(\sqrt{\dg\dH}+\tfrac{3}2\du\big)\partial_\dx \du = 0.
\end{myequation}
Conversely, any solution to the above equation provides a particular solution to \cref{eq.intro..SV-R} by setting $(\dr_-,\dr_+)=(2\sqrt{\dg\dH}+2 \du,-2\sqrt{\dg\dH})$, or equivalently a solution to \cref{eq.intro..SV} with $\dzeta=\dg^{-1}\big( \sqrt{\dg\dH}\du+\tfrac14\du^2\big)$.

\Cref{eq.intro..Hopf} may be solved by the hodograph transform, or the characteristics method, and exhibits a new phenomenon with respect to the linear equations discussed in previous sections: finite-time singularity formation. Assume $\du $ is a Lipschitz solution to \cref{eq.intro..Hopf} and define, for any $\dx_0\in\RR$, $\dv_{\dx_0}(\dt)\eqdef \du(t,\dx_{\dx_0}(\dt))$ where $\dx_{\dx_0}(\dt)$ is defined by the initial condition $\dx_{\dx_0}(t)=\dx_0$ and the ordinary differential equation $ \dx_{\dx_0}'(\dt)=\sqrt{\dg\dH}+\tfrac{3}2\du(\dt,\dx_{\dx_0}(\dt))$. Chain rule and \cref{eq.intro..Hopf} yields $ \dv_{\dx_0}'(\dt)
=0$, and hence $\dv_{\dx_0}(\dt)
=\du(0,\dx_0)$ and finally
$\dx_{\dx_0}(t)=\dx_0+(\sqrt{\dg\dH}+\tfrac{3}2\du(0,\dx_0))\dt$. In other words, the solution is constant along the characteristics defined by $\dx_{\dx_0}(\dt)$, for any $\dx_0\in\RR$, and the characteristics are straight lines. This allows to define and describe solutions as long as two characteristics do not cross, \ie as long as for any $\dt\in(0,\dT)$, there does not exists $\dx_0\neq \dx_1\in\RR$ with
\[\dx_0+\big(\sqrt{\dg\dH}+\tfrac{3}2\du(0,\dx_0)\big)\dt=\dx_1+\big(\sqrt{\dg\dH}+\tfrac{3}2\du(0,\dx_1)\big)\dt \quad \Longleftrightarrow \quad \tfrac{\du(0,\dx_1)-\du(0,\dx_0)}{\dx_1-\dx_0}=-\tfrac{2}{3 \dt}.\]
Hence we see that for any Lipschitz initial data $\du(\dt=0,\cdot)=\du_0\in W^{1,\infty}(\RR)$, the solution described above (which is unique) exists on the time domain $[0,\dT^\star)$ where $\dT^\star=-\frac2{3  }(\inf_{\RR} \du_0')^{-1}$ with the convention $\dT^\star=\infty$ if $\inf_{\RR} \du_0'\geq 0$. In the situation where $\inf_{\RR} \du_0'<0$ (in particular for any non-trivial $\du_0$ such that $\du_0\to 0$ as $|\dx|\to\infty$), there exists indeed a singularity formation as $t\to \dT^\star$: since the solution remains bounded but $\inf_{\RR}\partial_\dx u(\dt,\cdot)\to -\infty$ as $\dt\nearrow \dT^\star$, we say that a shock, or a {\em wavebreaking}, occurs. We represent this situation in \Cref{F.intro..wavebreaking}.

\begin{myfigure}[htbp]
	\begin{subfigure}{.47\textwidth}
		\movie{ShockSV}{60}
		\caption{Evolution in time}
	\end{subfigure}%
	\begin{subfigure}{.53\textwidth}
		\includegraphics[width=\textwidth]{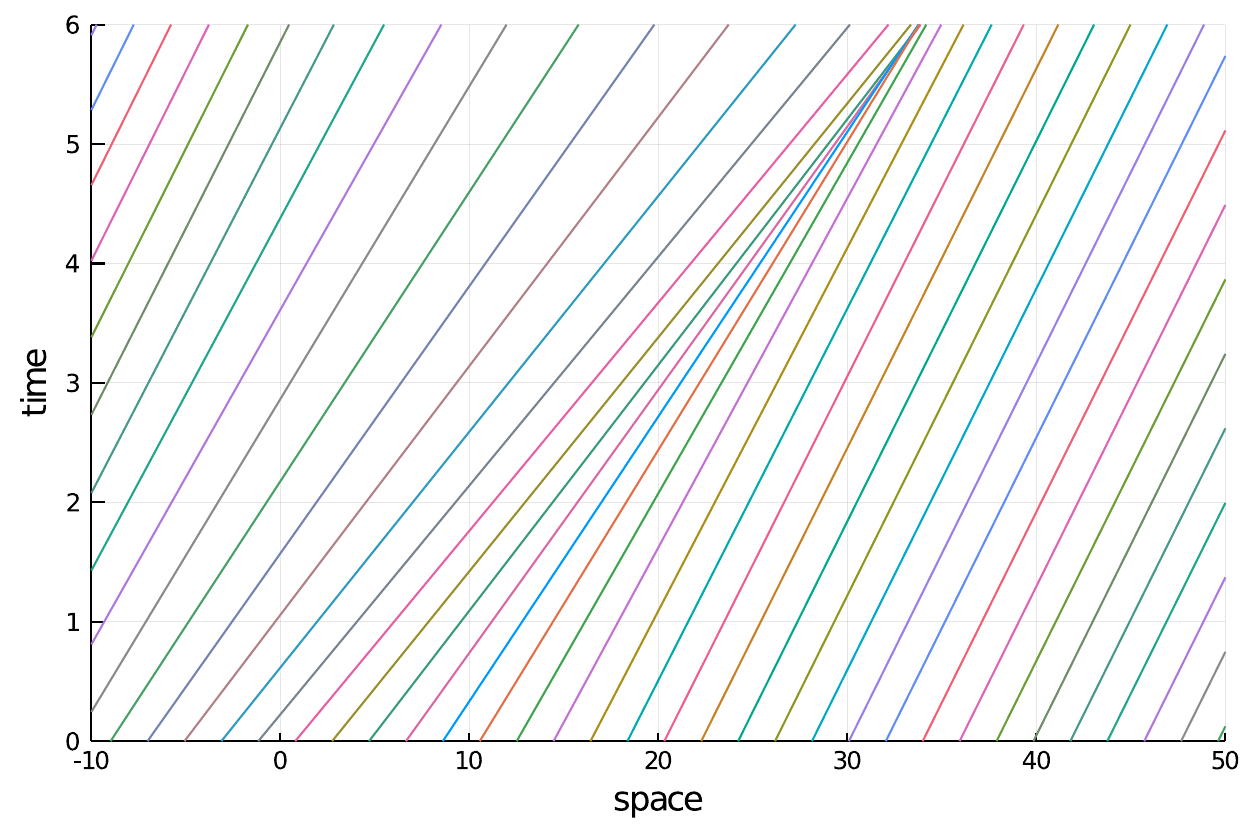}
		\caption{Characteristics}
	\end{subfigure}
	\caption[Wavebreaking according to the inviscid Burgers equation]{Wavebreaking of a simple wave according to \cref{eq.intro..Hopf}. The initial data for $\dzeta$ is the Gaussian $\dzeta(\dt=0,\dx)=0.5\exp(-(0.1\,\dx)^2)$ and corresponding velocity. $\dg=9.81\, {\rm m}.{\rm s}^{-2}$, $\dH=1\,{\rm m}$. }
	\label{F.intro..wavebreaking}
\end{myfigure}

Going back to the system case, \cref{eq.intro..SV-R}, each of the Riemann invariants, $\dr_\pm$, is constant along characteristics curves defined by
\[ \dx_{\pm,\dx_0}(0)=\dx_0, \quad \dx_{\pm,\dx_0}'(\dt)=\frac14(3\dr_\pm+\dr_\mp)(\dt,\dx_{\pm,\dx_0}(\dt)).\]
However the characteristics curves are no longer straight lines in general. Still we can infer the behavior of solutions for instance if we assume that  initial data $(\dzeta(t=0,\cdot),\du(t=0,\cdot))\eqdef(\dzeta_0,\du_0)$ have compact support, say in $(-\dL,\dL)$, and
are are sufficiently small so that there exists $\dc\in(0,\dc_0)$ with
\[ \dr_{+,0}\eqdef  \du_0+2\sqrt{\dH+\dzeta_0} \in (2\dc_0-\dc,2\dc_0+\dc) \quad \text{ and } \quad  \dr_{-,0}\eqdef  \du_0-2\sqrt{\dH+\dzeta_0} \in (-2\dc_0-\dc,-2\dc_0+\dc).\]
Because the Riemann invariants are constant along characteristics, we have, as long as the solution remains regular, $\frac{3\dr_++\dr_-}4\in (\dc_0-\dc,\dc_0+\dc)$ and $ \frac{3\dr_-+\dr_+}4\in (-\dc_0-\dc,-\dc_0+\dc)$, and as a consequence
\[ \dr_+(\dt,\dx)\equiv 2\dc_0 \text{ if $\dx\leq -\dL+(\dc_0-\dc)\dt$}  \quad \text{ and } \quad \dr_-(\dt,\dx)\equiv -2\text{ if $\dx\geq \dL-(\dc_0-\dc)\dt$}.\]
If the initial data is sufficiently small in order to ensure that no shock formation occurs before $\dT_\star=\frac{\dL}{\dc-\dc_0}$, we can afterwards decompose the flow as the superposition of two simple waves described by Hopf equations,
and in particular a shock inevitably occurs after sufficiently large time.

\section{Boussinesq systems}\label{S.intro..Boussinesq}

In his celebrated manuscript~\cite{Boussinesq72}, Boussinesq introduced the first models for the propagation of surface gravity waves taking into account both (first order) nonlinear and dispersive effects. While restricted in the original work to unidirectional waves, models with similar flavor were later on obtained for general waves. Eventually, one may obtain a full family of systems~\cite{MadsenSchaeffer99,BonaChenSaut02}, often called ($abcd$) {\bf{\em Boussinesq systems}}, of the form%
\footnote{The transport term $(\dbu\cdot\nabla)\dbu$ is often replaced with $\tfrac12\nabla(|\dbu|^2)$, trading the direct comparison with the Saint-Venant system, \cref{eq.intro..SV}, with conservative form. The change is immaterial in dimension $d=1$, or when $\rot\dbu=0$. Similar systems can be derived using momentum-type variables instead of velocity variables, thus slightly altering the nonlinear/dispersive interplay; see~\cite{FilippiniBellecColinEtAl15}. These systems, sometimes called Abbott systems~\cite{AbbottPetersenSkovgaard78,AbbottMcCowanWarren84}, have conservative form. Other {\em ad hoc} transformations can be performed, for instance to improve the mathematical properties of the system; see~\cite{BonaColinLannes05}. Finally, the models can also be written as second order scalar equations similar to \cref{eq.intro..Acoustic}, as in the original work of Boussinesq \cite[(26), p.~75]{Boussinesq72}:
	\begin{myequation}\label{eq.intro..Boussinesq-bad}
		\partial_\dt^2\dzeta=\dg\dH\, \Delta \big(\dzeta+\tfrac{3}{2\dH}\dzeta^2+ \tfrac{\dH^2}{3}\Delta\dzeta\big)\,.
	\end{myequation}
}
\begin{myequation}\label{eq.intro..Boussinesq}
	\left\{
	\begin{array}{l}
		\partial_\dt (\dzeta-b\dH^2\Delta\dzeta)+\nabla\cdot(\ddh\dbu+a\dH^3 \nabla\nabla\cdot \dbu )=0,\\[1ex] 
		\partial_\dt(\dbu-d\dH^2\nabla\nabla\cdot \dbu)+\dg\nabla(\dzeta+ c\dH^2\Delta\dzeta)+ (\dbu\cdot\nabla)\dbu=\bz,
	\end{array}
	\right.
\end{myequation}
where $\mfp=(a,b,c,d)\in\RR^4$ is such that (when neglecting surface tension) $a+b+c+d=\tfrac{1}{3}$. 
In \cref{eq.intro..Boussinesq} the precise meaning of the velocity variable depends on the choice of the parameters. The freedom in the choice of $(a,b,c,d)\in\mfp$ is at the same time a blessing---for instance one may tune parameters so as to enhance the accuracy of the dispersion relation---and a curse, since important properties of the system will typically depend on the choice of $(a,b,c,d)\in\mfp$. In particular, the initial-value problem of a subfamily of \cref{eq.intro..Boussinesq} is strongly ill-posed, as can be seen from modal instabilities of the linearized equations about the rest state: the dispersion relation being
\[\domega(\dbxi)^2=\dg\dH|\dbxi|^2\frac{(1-a|\dH\dbxi|^2)(1-c|\dH\dbxi|^2)}{(1+b|\dH\dbxi|^2)(1+d|\dH\dbxi|^2)}\,\]
with right-hand side taking arbitrarily large negative values at large wavenumbers, $|\dbxi|$, for ill-chosen $(a,b,c,d)\in\mfp$.
Incidentally, this is also the case for the original ``bad'' Boussinesq equation, \cref{eq.intro..Boussinesq-bad}.
This is a useful reminder that consistency is not the only property to look for in a model.

In the other way, it is expected that for ``good'' choices of $(a,b,c,d)\in\mfp$, dispersive properties of the Boussinesq systems prevent the wavebreaking scenario in the Saint-Venant model, \cref{eq.intro..SV}. As a matter of fact, for several families of parameters, $(a,b,c,d)\in\mfp$, global-in-time existence and uniqueness of solutions have been proved (see~\cite[Remark~1.1]{SautXu}) and---to the author's knowledge---the emergence of finite-time singularity has not been proved or numerically witnessed on any of the models, at least for solutions maintaining positive layer depth. In the situation of long waves and relatively large amplitude, the solution typically generates a zone of rapid oscillations (or modulations) often called {\em dispersive shock wave}, in place of the shock predicted by the Saint-Venant system. Properties of these dispersive shock waves will typically depend on the choice of $(a,b,c,d)\in\mfp$, and is not expected to accurately describe the real-life phenomenon.

An important property of nonlinear {\em and} dispersive equations such as \cref{eq.intro..Boussinesq} is that they allow the existence of {\em traveling waves}, that is solutions that maintain their shape while propagating at a constant velocity, including {\em solitary waves} which in addition bear finite energy. Once again the reader will find in \cite{Darrigol03} the fascinating and tumultuous story of the discovery and progressive  acceptance of these waves. Existence and properties of traveling waves again typically depend on the choice of $(a,b,c,d)\in\mfp$. We however expect that they exist at least for small supercritical velocities, $0<\dc-\dc_0\ll1$, and grow in amplitude with the velocity parameter; see \eg \cite{Dinvay}. We show examples in \Cref{F.intro..Boussinesq}.
\begin{myfigure}
	\begin{subfigure}{.5\textwidth}
		\movie{Boussinesq1}{100}
		\caption{Disintegration of Gaussian initial data}
	\end{subfigure}%
	\begin{subfigure}{.5\textwidth}
		\movie{Boussinesq2}{100}
		\caption{Traveling waves}
	\end{subfigure}
	\caption[Disintegration of a heap of water and solitary wave solutions of a Boussinesq system]{Left: Disintegration of the Gaussian $\dzeta(\dt=0, \dx)=0.25\exp(-0.1\dx^2)$, with zero initial velocity. Right: Solitary wave solutions with velocities  $\dc=1.05\dc_0$ and $\dc=1.01\dc_0$. \\ Both according to system~\eqref{eq.intro..Boussinesq} with  $-a=b=d=\tfrac{1}{3}$, $c=0$, $\dg=9.81\, {\rm m}.{\rm s}^{-2}$, $\dH=1\,{\rm m}$.}
	\label{F.intro..Boussinesq}
\end{myfigure}

It would be impossible to review all known results on Boussinesq systems and closely related (symmetric, Abbott, {\em etc.}) variants. Let me lazily refer to~\cite{DougalisMitsotakis08,Lannes20,SautXu}---in addition to previous references---and references therein,
and conclude with a last warning. The Boussinesq systems typically lose important properties of the original water waves equations and in particular its Hamiltonian structure. Hence unless the parameters $(a,b,c,d)\in\mfp$ are well-chosen,%
we do not expect energy conservation, or Galilean invariance, {\em etc.}

\section{The Korteweg--de Vries and Whitham equations}\label{S.intro..KdV}

It was mentioned in the previous section that Boussinesq's original motivation was the study of unidirectional waves, and in particular solitary waves. Using such assumption one may derive%
\footnote{We will not discuss in this document the interesting question of justifying such scalar equations from aforementioned systems of equations. Let me just mention that this justification is relatively straightforward for well-prepared initial data accounting for the assumption of unidirectional propagation, and much more involved for general initial data where we want to express that the flow can be decomposed at first order as the superposition of two counter-propagating unidirectional waves. Let me also refer---once again---to~\cite{Lannes} and references therein (see also \cite{Bambusi20} for a recent development) for all details concerning the Korteweg--de Vries equation and to~\cite{EmeraldII} for the Whitham equation.}
(as did Boussinesq) simplified scalar equations, of which the most famous is the {\bf{\em Korteweg--de Vries equation}}~\cite{Boussinesq77,KortewegDe95} for right-going waves in dimension $d=1$:
\begin{myequation}\label{eq.intro..KdV}
	\partial_\dt\dzeta+\sqrt{\dg\dH}\partial_\dx(\dzeta+\tfrac3{4\dH}\dzeta^2+\tfrac16\dH^2\partial_\dx^2\dzeta)=0\,.
\end{myequation}
One of the many reasons for the importance of the Korteweg--de Vries equation is the family of explicit solitary wave solutions\footnote{Far from being simply some entertaining special solutions, solitary waves play a very important role as they allow to describe the large-time dynamics of generic solutions; a phenomenon designated as {\em soliton resolution}. We will not discuss further on this feature as it relies on the integrability of the Korteweg--de Vries equation, a property which is not shared by other models in this document. }
\[\dzeta(\dt,\dx)=\dzeta_\dc(\dx-\dc t), \qquad \dzeta_\dc(\dx)\eqdef 2\dH(\tfrac{\dc}{\dc_0}-1)\sech^2\big(\sqrt{\tfrac3{2\dH^2}(\tfrac{\dc}{\dc_0}-1)} \,\dx \big)\]
where the velocity variable, $\dc$, may take any value $\dc>\dc_0=\sqrt{\dg\dH}$.


The existence of traveling waves with arbitrarily large amplitude and arbitrarily large velocity may found undesirable as nonphysical \cite{Stokes47}. Such is the case also for the global-in-time well-posedness properties, preventing the aforedescribed wavebreaking scenario. With this in mind, Whitham~\cite{Whitham67} introduced the following equation\footnote{He also proposed~\cite[\S13.14]{Whitham} 
	\begin{myequation}\label{eq.intro..Whitham-2}
		\partial_\dt\dzeta+\sqrt{\dg\tfrac{\tanh(\dH|D|)}{|D|}}\partial_\dx\dzeta+\big(3\sqrt{\dg(\dH+\dzeta)}-3\sqrt{\dg\dH}\big)\partial_\dx\zeta=0\,,
	\end{myequation}
	where the advection term fits the decomposition in Riemann invariants of the Saint-Venant system.}
which is now called {\bf{\em Whitham equation}}:
\begin{myequation}\label{eq.intro..Whitham}
	\partial_\dt\dzeta+\sqrt{\dg\dH}\partial_\dx\Big(\sqrt{\tfrac{\tanh(\dH|D|)}{\dH|D|}}\dzeta+\tfrac3{4\dH}\dzeta^2\Big)=0\,.
\end{myequation}
arguing that the fact that its linear dispersion relation reproduces exactly one branch of the dispersion relation of \cref{eq.intro..WW0} would authorize wavebreaking and peaked traveling waves of extreme height. This prediction turned out to be valid, as recently shown in~\cite{Hur17,EhrnstromWahlen19,TruongWahlenWheeler,SautWang}. A numerical comparison of solitary wave solutions to the Korteweg--de Vries and Whitam equations is shown in \Cref{F.intro..solitary}.
%

\begin{myfigure}[!htb]
	\begin{center}
		\begin{subfigure}{.44\textwidth}
			\includegraphics[width=\textwidth]{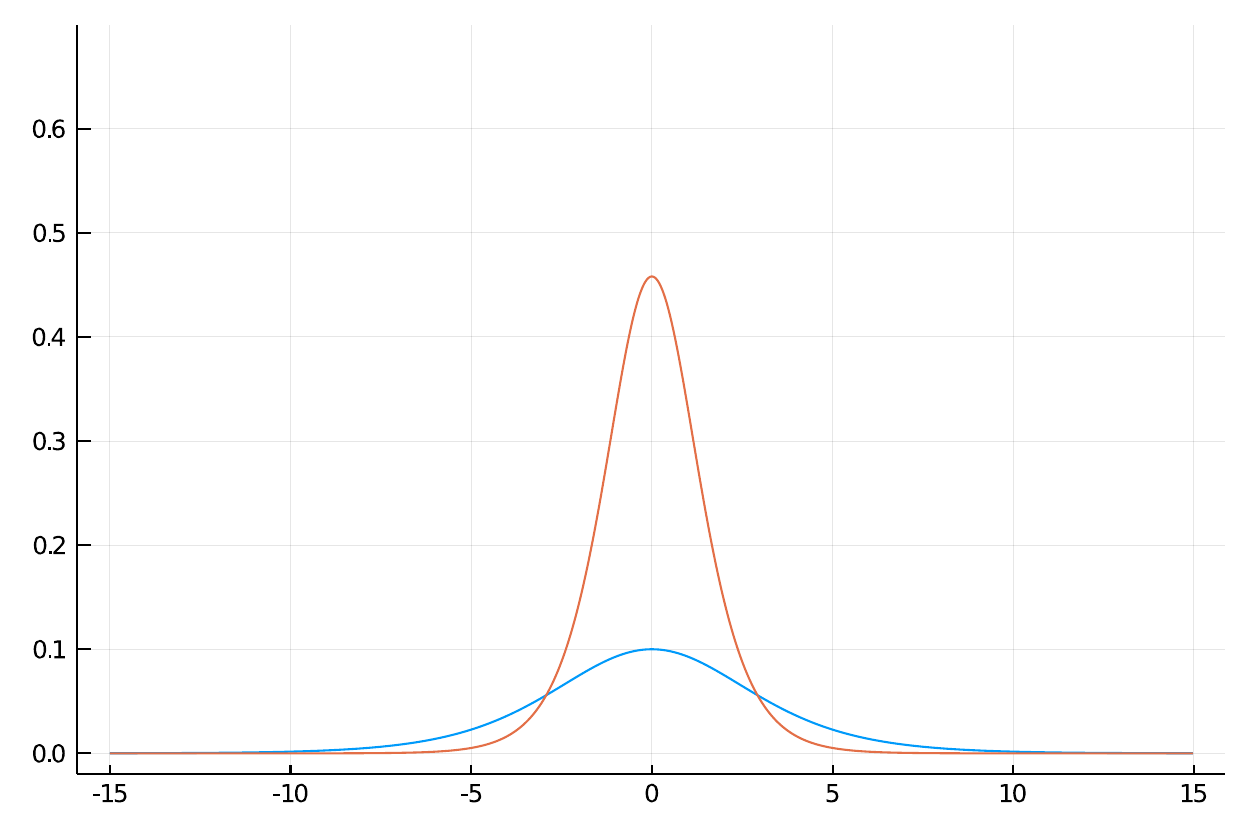}
			\caption{KdV equation, \cref{eq.intro..KdV}}
		\end{subfigure}%
		\begin{subfigure}{.44\textwidth}
			\includegraphics[width=\textwidth]{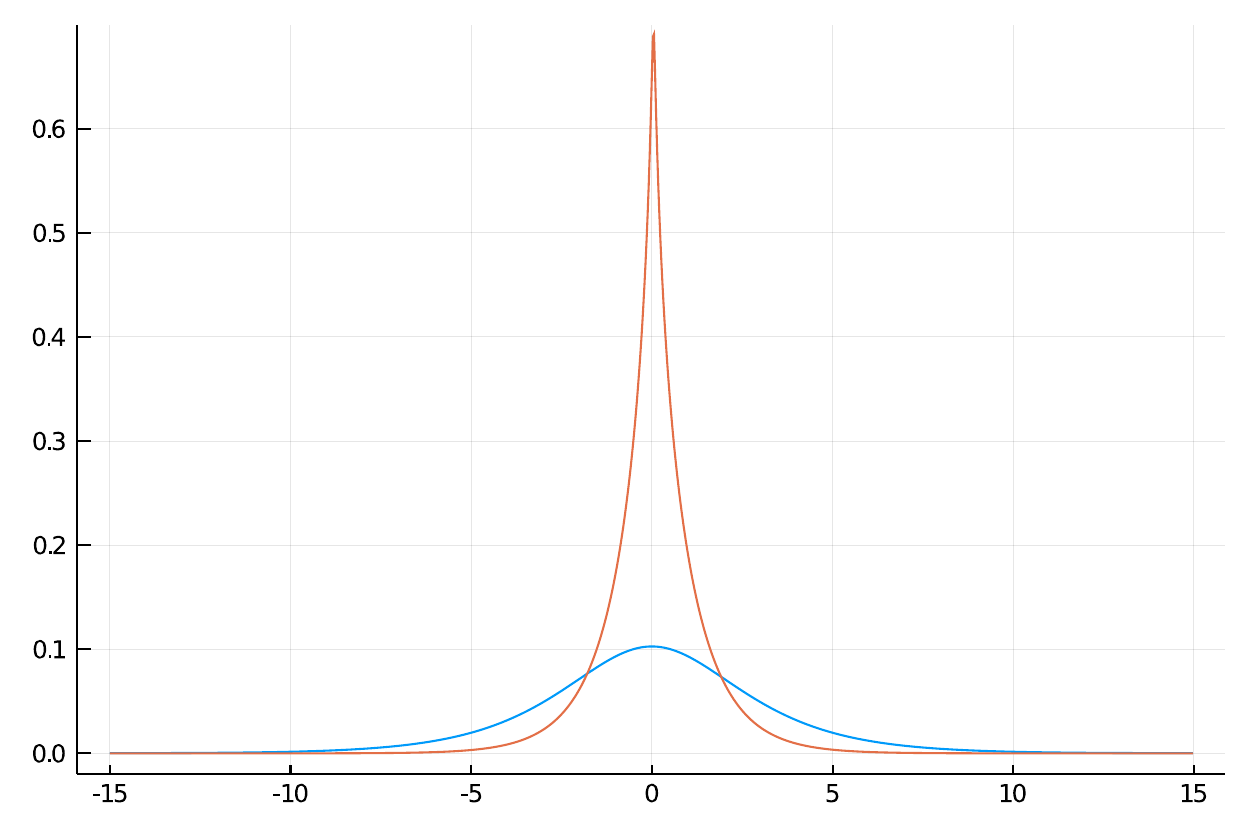}
			\caption{Whitham equation, \cref{eq.intro..Whitham}}
		\end{subfigure}
	\end{center}
	\caption[Solitary wave solutions to the Korteweg--de Vries and Whitham equations]{Solitary waves of unidirectional models with velocity $\dc=1.05 \dc_0$ (blue, smaller) and $\dc=1.2290408\dc_0$ (red, larger). $\dg=9.81\, {\rm m}.{\rm s}^{-2}$, $\dH=1\,{\rm m}$, $\dc_0=\sqrt{\dg\dH}$.}
	\label{F.intro..solitary}
\end{myfigure}




\mainmatter
\pagestyle{mystyle}
\renewcommand{\thesection}{\arabic{section}}
\setcounter{section}{0}

\chapter{The ``master'' equations}\label{C.master}
\vspace{-.2cm}
\epigraph{Le problème de l'établissement [...] des équations différentielles du mouvement, et ensuite de leur intégration approchée, aura encore sa difficulté souvent grande. Mais il ne présentera plus, envisagé ainsi, cette désespérante énigme contre laquelle des esprits distingués se sont heurtés en vain. 
}{--- \textsc{Adhémar Barré de Saint-Venant}, \textit{Comptes rendus des séances de l'Académie des sciences, séance du 18 mars 1872}}
\vspace{-.2cm}
\minitoc

\section*{Foreword}

In this chapter, we introduce and provide a preliminary study of the systems of equations from which asymptotic models are derived in subsequent chapters. The presentation, as well as most of the notations, are borrowed from Lannes' book~\cite{Lannes}. However concision has been pursued and I cannot encourage enough a thorough reading of the book for a detailed account.

We first write down the most general system of equations which is considered in this work, that is Euler equations for a layer of (non-necessarily homogeneous) incompressible ideal fluid, coupled with boundary conditions accounting for the impermeable bottom and the free surface. The only external force acting on the system will be the gravity force, assumed constant and vertical. We refer to the system we obtain as the {\bf{\em full Euler system}}. Then we focus on particular settings.

The homogeneous and irrotational framework is particularly rewarding, as it allows to rewrite the whole system as two evolution equations for unknown functions of time and horizontal space variables only. This system is referred to as the {\bf{\em water waves equations}}. 

Prominently important in the water waves equations is the {\em Dirichlet-to-Neumann operator}, which is defined after solving a Laplace problem on the fluid domain with Dirichlet and Neumann boundary conditions. Its study, and in particular the asymptotic expansions which allow to derive asymptotic models, are briefly reviewed.

Meanwhile we make a small step outside the world of homogeneous and potential flows to consider {\bf{\em interfacial waves}} between two layers of  homogeneous fluids with irrotational velocities. Additional Dirichlet/Neumann operators appearing in this framework are tackled.

\chapter{Hydrostatic models}\label{C.hydro}
\epigraph{``Begin at the beginning," the King said gravely, ``and go on till you
come to the end: then stop."}{--- \textsc{Lewis Carroll}, \textit{Alice in Wonderland}}
\minitoc

\begin{myFigure}[ht]
\begin{center}
\begin{tikzpicture}[every text node part/.style={align=center}]
\mytikzstyle

\draw[obvious] (0,-2)--(2,-2)node[midway,above]{\scriptsize specific case};
\draw[proved] (0,-3)--(2,-3)node[midway,above]{\scriptsize rigorous asymptotic};
\draw[formal] (0,-4)--(2,-4)node[midway,above]{\scriptsize formal asymptotic};

\node[model]  (ww) at (5.5,3) {water waves\\ \commentedref{eq.WW..dim}};
\node[model]  (ww2) at (11,3) {interfacial waves\\ \commentedref{eq.WW2..dim}};
\node[model]  (fE) at (0,3) {full Euler\\ \commentedref{eq.fE}};

\node[model]  (SV) at (5.5,0) {{\color{dgreen}{\em Saint-Venant}}\\ \commentedref{eq.SV..u-dim}};
\node[model]  (SV2) at (11,0) {{\color{dgreen}{\em bilayer hydrostatic}}\\ \commentedref{eq.SV2FS..u-dim}};
\node[model]  (SV2RL) at (11,-3) {{\color{dgreen}{\em bilayer hydrostatic}} \\  {\color{dgreen}{\em with rigid lid and Boussinesq}} \\ \commentedref{eq.SV2RL..u-p-dim-Bouss}};
\node[model]  (H) at (0,0) {{\color{dgreen}{\em hydrostatic equations}} \\ \commentedref{eq.hydro..dim}};


\draw[obvious] (fE)--(ww)node[midway,above]{\scriptsize homogen.}node[midway,below]{\scriptsize potential};
\draw[obvious] (ww2)--(ww)node[midway,above]{\scriptsize zero }node[midway,below]{\scriptsize upper density};

\draw[obvious] (H)--(SV)node[midway,above]{\scriptsize homogen.}node[midway,below]{\scriptsize columnar};
\draw[obvious] (SV2)--(SV)node[midway,above]{\scriptsize zero }node[midway,below]{\scriptsize upper density};

\draw[formal] (fE)--(H)node[midway]{\scriptsize shallow water, \\ $\mu\ll 1$};
\draw[proved] (ww2)--(SV2)node[midway]{\scriptsize shallow water, \\ $\mu\ll 1$};
\draw[proved] (ww)--(SV)node[midway]{\scriptsize shallow water, \\ $\mu\ll 1$};

\draw[proved] (SV2)--(SV2RL)node[midway]{\scriptsize weak density contrast, \\ $\gamma \nearrow 1$};


\end{tikzpicture}
\end{center}
\caption[{\bf Models in \Cref{C.hydro} and some filiations}]{Models in \Cref{C.hydro} (in {\color{dgreen}\em green}) and some filiations.}
\label{F.models-hydro-overview}
\end{myFigure}

~

\newpage

 \section*{Foreword}

We start our journey towards asymptotic models with ones among the oldest and simplest-looking. The so-called {\em hydrostatic models} can be formally derived from the ``master'' full Euler equations by using the {\em hydrostatic assumption}, that is approximating the pressure terms using an explicit formula stemming from neglecting the velocity advection terms in the horizontal momentum conservation equation, specifically
\[-\partial_\dz\dP=\drho\dg, \]
which we can integrate using the known pressure at the free surface. Additionally, one often adds the assumption of {\em columnar motion}, stating that the horizontal velocity (approximately) does not depend on the vertical variable. When both assumptions are made, then we quickly obtain models with the rewarding properties that the vertical space variable has disappeared from the equations and only (first order) differential operators are involved.

Yet we shall not assume {\em a priori}---but indirectly prove---the hydrostatic assumption nor the columnar motion and will rather justify models asymptotically---with quantitative error estimates---in the shallow water regime, as $\mu\ll1$; and using the irrotationality assumption in lieu of columnar motion. 

Our first model is derived from the water waves system, that is assuming that the density is homogeneous and the flow potential (as this allows to discard the vertical variable except in the Dirichlet-to-Neumann operator). We then obtain the well-known and much-studied {\bf{\em Saint-Venant system}}, already introduced in \Cref{S.intro..Saint-Venant}. Its derivation and rigorous justification, together with a very short description of some of its properties, is the subject of a first section.

Then we move to the {\bf{\em bilayer}} framework, with two layers of homogeneous potential flows. The situation is slightly messier as models differ whether we use the {\em free-surface} framework or the {\em rigid-lid} framework, and in the latter one often uses the so-called {\em Boussinesq approximation}. It turns out the rigid-lid assumption and Boussinesq approximation both follow from the same assumption of  weak density contrast.

Finally we quickly extend the analysis to $N\geq 2$ layers as above. While this {\bf{\em multilayer}} framework may appear artificial, it is expected to approximate (as $N\gg1$) the setting of continuously stratified flows, in view of withdrawing the assumptions of homogeneous density and potential flows while keeping the hydrostatic approximation in the shallow water regime.

Hydrostatic equations for {\bf{\em continuously stratified}} flows are also discussed. It turns out very little is known on these equations, despite the fact that they are at the core of the {\em primitive equations} which are widely used in studies and numerical simulations of geophysical flows. This offers stimulating mathematical challenges.

\chapter{Weakly dispersive models}\label{C.weakly}
\epigraph{parce que, [les Anciens] s'étant élevés jusqu'à un certain degré où ils nous ont portés, le moindre effort nous fait monter plus haut, et avec moins de peine et moins de gloire nous nous trouvons au-dessus d'eux. C'est de là que nous pouvons découvrir des choses qu'il leur était impossible d'apercevoir. }{--- \textsc{Blaise Pascal}, \textit{traité du vide}}
\medskip

\minitoc

\begin{myFigure}[ht]
\begin{center}\vspace{1cm}
\begin{tikzpicture}[every text node part/.style={align=center}]
\mytikzstyle

\draw[obvious] (6,6)--(8,6)node[midway,above]{\scriptsize specific case};
\draw[proved] (9,6)--(11,6)node[midway,above]{\scriptsize rigorous asymptotic}node[midway,below]{\tiny (convergence)};

\node[model]  (ww) at (2.5,6) {water waves\\ \commentedref{eq.WW..dim}};
\node[model]  (GN) at (5,4) {{\color{dgreen}\em Green--Naghdi}\\ \commentedref{eq.GN..ddot-dim}};
\node[model]  (WGN) at (0,4) {{\color{dgreen}\em Whitham--Green--Naghdi}\\ \commentedref{eq.WGN..u-dim}};
\node[model]  (FG) at (10,4) {{\color{dgreen}\em Favrie--Gavrilyuk}\\ \commentedref{eq.FG..dim}};
\node[model]  (B) at (5,2) {{\color{dgreen}\em Boussinesq}\\ \commentedref{eq.B..u-dim}};
\node[model]  (WB) at (0,2) {{\color{dgreen}\em Whitham--Boussinesq}\\ \commentedref{eq.WB..psi-dim}};
\node[model]  (SV) at (5,0) {Saint-Venant\\ \commentedref{eq.SV..u-dim}};
\node[model]  (L) at (0,0) {linearized (Airy)\\ \commentedref{eq.WW..lin-dim}};

\draw[proved] (ww.south east)--(GN)node[midway,right]{ $\cO(\mu^2)$};
\draw[proved] (ww.south west)--(WGN)node[midway,right]{$\cO(\mu^2(\eps+\beta))$};
\draw[proved] (GN)--(FG)node[midway,above]{ $\cO(\lambda^{-1}\mu)$};
\draw[proved] (GN)--(B)node[midway,right]{ $\cO(\mu(\eps+\beta))$ };
\draw[proved] (WGN)--(WB)node[midway,right]{ $\cO(\mu(\eps+\beta))$ };
\draw[obvious] (B)--(SV)node[midway,right]{ $\mu=0$ };
\draw[obvious] (WB)--(L)node[midway,right]{ $\eps =  0$ };
\draw[obvious] (WB.south east)--(SV.north west)node[midway,right]{ $\ \mu=0$};
\draw[obvious] (GN.south east) to[out = -30,in=30] node[pos=0.5,right]{$\mu=0$} (SV.north east);
\draw[proved] (WGN)--(GN)node[midway,above]{$\cO(\mu^2)$};
\draw[proved] (WB)--(B)node[midway,above]{$\cO(\mu^2)$};

\end{tikzpicture}
\end{center}
\caption[{\bf Models in \Cref{C.weakly} and some filiations}]{Models in \Cref{C.weakly} (in {\color{dgreen}\em green}) and some filiations.}
\label{F.models-weakly}
\end{myFigure}

~

\newpage

\section*{Foreword}

This chapter is devoted to the derivation and analysis of {\em weakly dispersive models}. These models refine the hydrostatic equations studied in \Cref{C.hydro} (in fact, more precisely, the Saint-Venant equations since here we restrict the analysis to the water waves framework; see \Cref{C.bilayer} for an extension to the bilayer framework) by introducing dispersive effects at first order. More refined models are presented in \Cref{C.higher}. 

The first section of this chapter 
 has been meant as a showcase for a thorough study of water waves models. Here we introduce and analyze the so-called {\bf{\em(Serre--)Green--Naghdi model}}. Firstly the model is quickly derived as an asymptotic model from the expansion of the Dirichlet-to-Neumann operator obtained in \Cref{C.master}. Yet the result which follows (namely the {\em consistency} of the model) is far from being sufficient to validate the Green--Naghdi equations as a good model for water waves. Firstly, its rigorous justification must be completed by {\em well-posedness}, {\em stability} and {\em convergence} results. They follow from careful energy estimates in suitable functional spaces. In a looser way, we also expect ``good'' models to retain important properties of the master equations (here the water waves system). Here we focus mostly on the variational structure of the equations: we observe that the Green--Naghdi equations not only preserve Zakharov's canonical Hamiltonian structure of the water waves system, but it also enjoys a deeper Lagrangian formalism which embeds the system inside a natural family of conservation laws, which can be interpreted as equations for compressible fluids with inertia effects. Hence the structure of the equations becomes richer as we simplify the equations from the water waves system to the Green--Naghdi equations (and then from the Green--Naghdi equations to the Saint-Venant system). This explains in my opinion why the Green--Naghdi equations, among many other loosely equivalent models, has attracted so much attention from diverse communities. 
We review some basic properties of the Green--Naghdi equations: preserved quantities, modal analysis and dispersion relation, solitary wave solutions. Finally, some open questions are discussed.

Of course I do not claim that the Green--Naghdi model is perfect! One of its main drawback is certainly that numerically approximating the equations turns out to be quite costly. In a second section we explore some equations which have been proposed by Favrie and Gavrilyuk~\cite{FavrieGavrilyuk17} to circumvent this issue. The equations are constructed using the aforementioned Lagrangian formalism, using a strategy akin to relaxation limits. Hence the system contains additional unknowns as well as a (large) parameter which is expected to measure the precision of solutions to the augmented equations with respect to solutions to the original Green--Naghdi equations, at least when initial data are well-prepared. The rigorous study of this {\em singular limit} is based on~\cite{Duchene19}. Again, the section is concluded by perspectives and open questions.

In a third section we introduce a {\em fully dispersive} analogue of the Green--Naghdi system, which we name {\bf{\em Whitham--Green--Naghdi}}. When linearized about trivial equilibrium solutions, fully dispersive models coincide with the corresponding (Airy) linearized water waves equations, as introduced in \Cref{S.intro..Airy}. Interest in such fully dispersive models in the context of long water waves started with the work of Whitham, which proposed \cref{eq.intro..Whitham} and \cref{eq.intro..Whitham-2} as suitable modifications of the standard Korteweg-de Vries equation, \cref{eq.intro..KdV}, in view of reproducing at least qualitatively important features of water waves such as wavebreaking and peaked traveling waves. Much more recently, the interest was renewed as Whitham's predictions were proved to be valid~\cite{Hur17,EhrnstromWahlen19,TruongWahlenWheeler,SautWang}. Yet the question of validating fully dispersive models as asymptotic models with improved accuracy with respect to their standard counterparts was mostly left aside. The precision of the Whitham--Green--Naghdi model (respectively {\bf{\em Whitham--Boussinesq}}) introduced in this section significantly improves the precision of the Green--Naghdi (respectively {\bf{\em Boussinesq}}) model for weak nonlinearities and mild bottom variations with the important price to pay that nonlocal operators (Fourier multipliers) are involved. 
These models also allow to rigorously justify the Whitham equations and observe a similar improvement with respect to the Korteweg--de Vries equation~\cite{EmeraldII}.

\chapter{Higher order models}\label{C.higher}
\vspace{-20pt}
\epigraph{Jésus a dit : « Que celui qui cherche ne cesse pas de chercher, jusqu’à ce qu’il trouve. Et quand il aura trouvé, il sera troublé ; quand il sera troublé, il sera émerveillé, et il régnera sur le Tout. » }{--- \textsc{Thomas l'apôtre}, \textit{évangile selon Thomas}}
\minitoc

\begin{myFigure}[ht]
\begin{center}\vspace{1cm}
\begin{tikzpicture}[every text node part/.style={align=center}]
\mytikzstyle

\draw[obvious] (-1,-2)--(1,-2)node[midway,above]{\scriptsize specific case};
\draw[proved] (2,-2)--(4,-2)node[midway,above]{\scriptsize rigorous asymptotic}node[midway,below]{\tiny (convergence)};
\draw[formal] (5,-2)--(7,-2)node[midway,above]{\scriptsize formal asymptotic}node[midway,below]{\tiny (consistency)};
\draw[conjecture] (8,-2)--(10,-2)node[midway,above]{\scriptsize expected asymptotic}node[midway,below]{\tiny (conjecture)};

\node[model]  (ww) at (0,3) {water waves\\ \commentedref{eq.WW..dim}};
\node[model]  (IK) at (5,6) {{\color{dgreen}\em Isobe--Kakinuma}\\ \commentedref{eq.IK..dim}};
\node[model]  (F) at (5,4) {{\color{dgreen}\em Friedrichs-type}\\ 
	 \commentedref{eq.HOF..Mat-N-adim}};
\node[model]  (aGN) at (5,2) {{\color{dgreen}\em augmented}\\ {\color{dgreen}\em Green--Naghdi}\\ \commentedref{eq.HOG..aGN-dim}};
\node[model]  (mGN) at (5,0) {{\color{dgreen}\em ``multilayer''}\\ {\color{dgreen}\em Green--Naghdi}\\ \commentedref{eq.HOG..mGN-dim}};
\node[model]  (GN) at (10,2) { Green--Naghdi\\ \commentedref{eq.GN..ddot-dim}};
\node[model]  (SV) at (10,4) {Saint-Venant\\ \commentedref{eq.SV..u-dim}};

\draw[proved] (ww)--(IK)node[pos=0.5,left]{$\cO(\mu^{2N+1})$};
\draw[formal] (ww)--(F)node[pos=0.5,above]{ $\cO(\mu^N)$};
\draw[conjecture] (ww)--(aGN)node[pos=0.5,below]{$\cO(\mu^{2N})$};
\draw[conjecture] (ww)--(mGN)node[pos=0.5,left]{$\cO(\mu^2 N^{-2})\ $};

\draw[obvious] (IK)--(SV)node[pos=0.5,above]{ $\ \ \ N=0$ };
\draw[obvious] (F)--(SV)node[pos=0.5,above]{ $N=1$ };
\draw[obvious] (F)--(GN)node[pos=0.5,above]{ $\ \ \ N=2$ };
\draw[obvious] (aGN)--(GN)node[pos=0.5,below]{ $N=1$ };
\draw[obvious] (mGN)--(GN)node[pos=0.5,below]{ $\ \ \ N=1$ };


\end{tikzpicture}
\end{center}
\caption[{\bf Models in \Cref{C.higher} and some filiations}]{Models in \Cref{C.higher} (in {\color{dgreen}\em green}) and some filiations.}
\label{F.models-higher}
\end{myFigure}

\section*{Foreword}

In this chapter we introduce and discuss {\em higher order models} for the water waves system, building upon the Saint-Venant system (\Cref{C.hydro}) and the Green--Naghdi system (\Cref{C.weakly}). These are hierarchies of models, that is families of system depending on a parameter---always denoted $N$---which we call the rank of the model, of which the Saint-Venant and/or the Green--Naghdi system are typically the first rank elements. 
The Saint-Venant (resp. Green--Naghdi) system has been rigorously justified as a shallow water model for the water waves system, in---roughly speaking---the following way: the size of the difference between solutions to the dimensionless water waves system and the corresponding solutions to the model equations grows proportionally to the size of the initial data with a prefactor bounded as $C\,\mu\, t$ (resp. $C\,\mu^2\,t$) over a relevant time interval (being of size inversely proportional to the size of the initial data), where $\mu$ is the shallow water parameter, and $C$ depends on an upper bound on the size of the admissible initial data (together with a lower bound on the minimal depth of the layer, an upper bound on admissible values for $\mu$, and the norms measuring the size of the data). In good cases we expect that a similar result holds for all elements in a hierarchy of models, with different prefactors $C_N\,\mu^{\alpha_N}\,t$. There are typically two situations:
\begin{enumerate}
\item the order as a shallow water model increases with $N$, that is $\alpha_N\to \infty$ as $N\to\infty$;
\item the accuracy of the model improves with $N$, that is $C_N\to 0$ as $N\to\infty$.
\end{enumerate}
In the latter but not in the former we can hope that the hierarchy provides a robust tool for the approximation of {\em any} (sufficiently regular) solution to the water waves system, and can be useful for instance to devise strategies for their numerical integration.

This chapter is decomposed into three sections, corresponding to three different strategies, each producing a variety of families of higher order models.\footnote{\label{f.spectral}The list is by no means complete. In particular it lacks {\em spectral methods} based on expansions with respect to the steepness parameter, $\sep=\eps\sqrt\mu$, initiated in \cite{DommermuthYue87,WestBruecknerJandaEtAl87,CraigSulem93} (see \eg \cite{LeTouze03,Schaeffer08,WilkeningVasan15,Nicholls16} for a detailed account and comparisons). Among them the strategy brought to light by Craig and Sulem in~\cite{CraigSulem93}, consisting in expanding the Dirichlet-to-Neumann operator, $\cG^\mu[\eps\zeta]$, along the variable $\eps\zeta$, is particularly elegant and effective. Contrarily to the models introduced in this chapter, the family of models involve Fourier multipliers in addition to differential operators; see \cite{Choi19a} for discussion, references and the explicit display of the models up to fifth order. This method has been extended and successfully employed in many situations (see \cite{Guyenne19} and references therein), despite the claim---based on numerical experiments and formal arguments---in \cite{AmbroseBonaNicholls14,MelinandDuchene} that the Cauchy problem associated with systems in the family are ill-posed in Sobolev spaces.  } 

In the first section, we use an expansion due to Boussinesq~\cite{Boussinesq73} and Rayleigh~\cite{Rayleigh76} (see~\cite[\S4.1]{Dingemans97} for discussion and other relevant references) of the velocity potential---as a solution to the Laplace problem---as a series involving powers of the shallow water parameter, $\mu$. We are hence typically in the framework of the first aforementioned situation, and this section emphasizes its possible shortcomings. Among the different models which can be naturally constructed by this way---which we call Friedrichs-type systems in acknowledgment to his Appendix to~\cite{Stoker48}---we introduce explicitly two families of models: the {\bf {\em high order shallow water models}} and the {\bf{\em extended Green--Naghdi models}}. These models involve differential operators of increasing order as the rank of the model grows, which yields several complications. Firstly, half of these models suffer from very serious high frequency instabilities which prevent any hope as for the well-posedness of the Cauchy problem. But even in good cases, it is expected that for fixed initial data the solutions to the systems---if they exist---do not converge towards the corresponding solution to the water waves system, as $N\to\infty$. This can be seen in particular when studying the dispersion relation of the models, which converge towards the dispersion relation of the water waves system only for wavenumbers in a finite-size neighborhood of the origin.

In the second section we set up a Galerkin dimension reduction strategy to a reformulation of the Laplace problem, to devise the approximate formula for the velocity potential---or, more precisely, the horizontal velocity. As a second step, the usual procedure consists in using this approximate formula in the Hamiltonian functional of the water waves system, and express the model as the canonical Hamiltonian equations associated with the approximate Hamiltonian. This procedure produces a different model for any (reasonable) choice of subspace of real-valued functions of the fluid domain used in the Galerkin method. Natural examples of such spaces in the shallow water framework are functions of the form
\begin{equation}\tag{$\star$} \label{eq.HO..decomposition}
\Phi(t,\bx,z) = \sum_{i=1}^N \phi_i(t,\bx)\,  \Psi_i(\bx,z)  
 \end{equation}
where $\phi_i(t,\bx)$ are variable unknowns of the resulting model, which is characterized by the choice of the vertical distribution, $ \{\Psi_i\}_{i\in\{1,\dots,N\}} $. We explore the outcome of vertical distribution defined, following the finite element method, as piecewise polynomials in the vertical variable, $z$. We particularly emphasize two families of models (respectively playing with the degrees of the polynomials and the number of elements in the vertical discretization): the {\bf {\em augmented Green--Naghdi models}} and 
the  {\bf{\em ``multilayer'' Green--Naghdi models}}. In each case, the system consists in two evolution equations coupled with a system of differential equations of order two mimicking the Laplace problem. The first family is a higher order shallow water hierarchy comparable to the models of the preceding section, yet instead of involving high order differential operators, the size of the system of differential equations grows with the rank, $N$. The second family has different  properties, akin to the second situation described above. The term ``multilayer'' stems from the fact that the models can be interpreted as resulting from the vertical discretization of the fluid layer in $N$ prescribed---typically proportional---sublayers.

In the third section we describe the strategy referred in \cite{KlopmanGroesenDingemans10,Klopman10} as  ``variational'' (see~\cite{Papoutsellis17} for an overview of related earlier and subsequent works). Of course the preceding strategy was also variational in nature, and we argue that the two strategies in fact differ only by the choice of the variational formulation of the Laplace problem. Yet in the latter, we plug directly the decomposition into Luke's Lagrangian action for the water waves system, and let Hamilton's principle do all the work in one single step. The outcome is surprising at first, as we obtain an overdetermined/underdetermined composite system  of $N$ evolution equations for the surface deformation, $\zeta$, and only one evolution equation for $(\phi_1,\ldots,\phi_N)$. Yet the systems can in fact be written---as in the above hierarchies---under a canonical Hamiltonian formulation of two evolution equations coupled with a system of differential equations of order two. Again each choice of the vertical distribution, $ \{\Psi_i\}_{i\in\{1,\dots,N\}} $, yields a different model. We only quickly mention the ``multilayer'' systems and instead focus on the shallow water system named the {\bf{\em Isobe--Kakinuma models}} in reference to~\cite{Isobe94a,Kakinuma01}. Indeed the latter benefit from a rigorous justification theory, thanks to the work of Iguchi and collaborators, culminating with~\cite{Iguchi18a}.

\chapter{Non-hydrostatic models for interfacial waves}\label{C.bilayer}

\epigraph{Toujours vouloir tout essayer, et recommencer}{--- \textsc{Michel Berger}, \textit{Le paradis blanc}}
\minitoc

\begin{myFigure}[ht]
\begin{center}
\begin{tikzpicture}[every text node part/.style={align=center}]
\mytikzstyle

\draw[obvious] (-7,1.5)--(-5,1.5)node[midway,above]{\scriptsize specific case};
\draw[proved] (-7,0.5)--(-5,0.5)node[midway,above]{\scriptsize rigorous asymptotic}node[midway,below]{\tiny (limit system is well-posed)};
\draw[formal] (-7,-0.5)--(-5,-0.5)node[midway,above]{\scriptsize formal asymptotic}node[midway,below]{\tiny (limit system is ill-posed)};

\node[model]  (ww2) at (-6,3) {interfacial waves\\ \commentedref{eq.WW2..dim}};
\node[model]  (IK) at (5.5,6) {Isobe--Kakinuma\\ \commentedref{eq.IK..dim}};
\node[model]  (K) at (0,6) {{\color{dgreen}\em Kakinuma}\\ \commentedref{eq.K..dim}};
\node[model]  (WGN) at (5.5,4) {Whitham--Green--Naghdi\\ \commentedref{eq.WGN..u-dim}};
\node[model]  (WCC) at (0,4) {{\color{dgreen}\em Whitham--Choi--Camassa}\\ \commentedref{eq.DIT..u-dim,eq.DIT..def-WCC}};
\node[model]  (SGN) at (5.5,2) {Serre--Green--Naghdi\\ \commentedref{eq.GN..ddot-dim}};
\node[model]  (MCC) at (0,2) {{\color{dgreen}\em Miyata--Choi--Camassa}\\ \commentedref{eq.MCC..ddot-dim}};
\node[model]  (SV) at (5.5,0) {Saint-Venant\\ \commentedref{eq.SV..u-dim}};
\node[model]  (SV2) at (0,0) {bilayer hydrostatic\\ \commentedref{eq.SV2RL..u-p-dim}};

\draw[proved] (ww2)--(K)node[pos=0.7,left]{$\cO(\mu^{2N+1})$};
\draw[formal] (ww2)--(WCC)node[pos=0.5,above]{$\cO(\mu^2(\eps+\beta))$};
\draw[formal] (ww2)--(MCC)node[pos=0.5,below]{ $\cO(\mu^2)$};
\draw[proved] (ww2)--(SV2)node[pos=0.5,below]{ $\cO(\mu)$ };

\draw[obvious] (K)--(IK)node[midway,above]{ $\gamma=0$ }node[midway,below]{ $\delta=1$ };
\draw[obvious] (WCC)--(WGN)node[midway,above]{ $\gamma=0$ }node[midway,below]{ $\delta=1$ };
\draw[obvious] (MCC)--(SGN)node[midway,above]{ $\gamma=0$ }node[midway,below]{ $\delta=1$ };
\draw[obvious] (SV2)--(SV)node[midway,above]{ $\gamma=0$ }node[midway,below]{ $\delta=1$ };

\end{tikzpicture}
\end{center}

\caption[{\bf Models in \Cref{C.bilayer} and some filiations}]{Models in \Cref{C.bilayer} (in {\color{dgreen}\em green}) and some filiations.  }
\label{F.models-bilayer}
\end{myFigure}

~

\newpage

\section*{Foreword}

In \Cref{C.master} the equations extending the water waves system to the framework of interfacial waves between two layers of incompressible, homogeneous, inviscid and immiscible fluids with potential flows are introduced. The physical motivation for studying such systems is the reported (and ubiquitous) existence of coherent waves traveling at the sharp interface between, say, fresh and warm water above denser salted cold water. One can refer to \eg~\cite{Jackson04,HelfrichMelville06} for a small peek at the vast literature on the subject. The main features of these waves is that they have tremendously large amplitudes---sometimes of the order of magnitude of the layer itself---, very long wave length, and travel over very long distances. Hence the assumptions of the shallow water regime, and in particular the fact that we do not impose any smallness assumption on the amplitude of the wave, is perfectly suited to the study of such waves. It is therefore very tempting to introduce asymptotic models for interfacial waves which are analogous to the asymptotic models for the water waves system. This is done in the hydrostatic framework in \Cref{C.hydro} and non-hydrostatic models are the subject of this chapter.

In addition to the physical motivation, there are interesting new features and challenges when studying interfacial waves. First and foremost, three additional dimensionless parameters come into view, namely
\[\alpha=\frac{\da_{\rm top}}{\da_{\rm int}} \quad ; \quad \delta=\frac{\dH_1}{\dH_2} \quad ; \quad \gamma = \frac{\drho_1}{\drho_2}\]
respectively the amplitude ratio of the free surface and interface, the depth ratio between the two layers, and the density ratio. Hence there are plethora of interesting limits to consider. We will focus here on the framework which is the most similar to the one-layer case%
\footnote{see \eg~\cite{ChoiCamassa99,BonaLannesSaut08} for related studies in other physically relevant asymptotic regimes.} 
and in particular we will assume that the two layers are of comparable depth, both small with respect to the typical horizontal wavelength of the flow.
The relation between the limit of small density contrast, $\gamma \nearrow 1$ and the rigid-lid hypothesis, is discussed in details in the hydrostatic framework in \Cref{C.hydro}. Somewhat inconsistently, we will restrict henceforth to the rigid-lid situation%
 \footnote{see \eg~\cite{ChoiCamassa96,Duchene10,Duchene12} for related studies in the free-surface framework.}
 without assuming the Boussinesq approximation, yet allowing $\gamma$ to approach unity. 
 To summarize, our results hold for parameters in the following set.
\begin{Definition*}[Shallow water/Shallow water asymptotic regime]
Given $\mu^\star,\delta_\star,\delta^\star>0$, we let
\[\mfp_{\SWSW}=\big\{(\mu,\eps,\beta,\delta,\gamma)\ : \ \mu\in(0,\mu^\star],\ \eps\in[0,1],\ \beta\in [0,1],\  \delta\in[\delta_\star,\delta^\star],\gamma\in[0,1)\big\}.\]
\end{Definition*}

One of the main striking difference between the water waves system and the corresponding interfacial waves system is the emergence of {\em Kelvin--Helmholtz instabilities} in the latter. Recall that the provided modal analysis shows that large wavenumber modes are unstable, and that the exponential growth rate takes arbitrarily large values as the wavenumber goes to infinity. This explains why the initial-value problem associated with the nonlinear system is strongly  ill-posed outside of the analytic framework. This appears to contradict the fact that, as we said, large interfacial waves do exist and appear remarkably stable! An answer to this paradox has been given by Lannes in~\cite{Lannes13}, by introducing interfacial tension effects: it is shown that well-posedness is restored, and more importantly the time of existence of solutions grows as $\mu\searrow 0$, consistently with the fact that the hydrostatic equations for interfacial waves are well-posed. It should be emphasized however that interfacial tension is not physical at the interface between two miscible fluids such as fresh and salted water; here it plays the role of a regularizing operator acting mostly on the high (spatial) frequency component of the flow. The real physical explanation is that mixing occurs, yet on a very thin transition layer: the {\em pycnocline}. In the absence---to my knowledge---of a simple expression revealing the effective influence of such mixing in the equations, we shall discard any effect when deriving asymptotic models.

It should be emphasized however that the models studied in this manuscript behave very differently regarding Kelvin--Helmholtz instabilities. Indeed, as the derivation focuses on the low frequency (large wavelength) component of the flow, the high frequency behavior can be very dissimilar between the different models, and hence with respect to the original interfacial waves system. A key revelation of the forthcoming study is the following.
\begin{itemize}
\item The {\bf {\em Miyata--Choi--Camassa}} model, which is analogous to the Green--Naghdi system (see \Cref{C.weakly}), overestimates Kelvin--Helmholtz instabilities.
\item This unfortunate behavior can be corrected through artificial---but harmless for the precision (in the sense of consistency) of the asymptotic model---modifications, which naturally yields fully dispersive systems named {\bf {\em Whitham--Choi--Camassa}}; or regularized systems. 
\item The {\bf{\em Kakinuma}} model, which extends the Isobe--Kakinuma model (see \Cref{C.higher}) to the bilayer framework, inherently tames Kelvin--Helmholtz instabilities.
\end{itemize}

The latter model can be expected to be useful for understanding the propagation of long interfacial waves, focusing on the large-scale dynamics of the flow, and discarding small-scale effects as irrelevant. Once again, this should not blurry the fact that mixing do occur, and may in some circumstances play an important role on the large-scale dynamics. Models with the aim of tracking these effects---at least at first order---should use the continuously stratified Euler equations as a starting point. Yet as mentioned in \Cref{C.hydro}, very little is known for this system in the shallow water regime. An important reference---in my opinion---dealing with long weakly dispersive internal (and not interfacial) waves is~\cite{DesjardinsLannesSaut19}. The Perspectives section in that reference supports and complements the present discussion.

\bookmarksetup{startatroot}
\appendix

\pagestyle{appendix}
\renewcommand{\theHsection}{Romansection.\thesection} 
\renewcommand{\thesection}{\Roman{section}}
\phantomsection 
\addtocontents{lof}{\vspace{10pt}} 
\addstarredchapter{Appendix}


\minitoc


  
  \pagestyle{biblio}             
  \phantomsection 
  \addcontentsline{toc}{chapter}{Bibliography}
  \providecommand{\arXiv}[1]{arXiv preprint:\href{https://arxiv.org/abs/#1}{#1}}
\providecommand{\noopsort}[1]{}\def\cprime{$'$}

 \thispagestyle{biblio}          



\begin{thebibliography}{100}
	
	\bibitem{AbarbanelHolmMarsdenEtAl86}
	H.~D.~I. Abarbanel, D.~D. Holm, J.~E. Marsden, and T.~S. Ratiu.
	\newblock Nonlinear stability analysis of stratified fluid
	equilibria\label{S.bibliostart}.
	\newblock {\em Philos. Trans. Roy. Soc. London Ser. A}, 318(1543):349--409,
	1986.
	
	\bibitem{AbbottPetersenSkovgaard78}
	M.~Abbott, H.~Petersen, and O.~Skovgaard.
	\newblock On the numerical modelling of short waves in shallow water.
	\newblock {\em J. Hydr. Res.}, 16(3):173--204, 1978.
	
	\bibitem{AbbottMcCowanWarren84}
	M.~B. Abbott, A.~D. McCowan, and I.~R. Warren.
	\newblock Accuracy of short-wave numerical models.
	\newblock {\em J. Hydr. Eng. ASCE}, 110(10):1287--1301, 1984.
	
	\bibitem{AbramowitzStegun64}
	M.~Abramowitz and I.~A. Stegun.
	\newblock {\em Handbook of mathematical functions with formulas, graphs, and
		mathematical tables}, volume~55 of {\em National Bureau of Standards Applied
		Mathematics Series}.
	\newblock U.S. Government Printing Office, Washington, D.C., 1964.
	
	\bibitem{Aceves-SanchezMinzoniPanayotaros13}
	P.~Aceves-S{\'a}nchez, A.~A. Minzoni, and P.~Panayotaros.
	\newblock Numerical study of a nonlocal model for water-waves with variable
	depth.
	\newblock {\em Wave Motion}, 50(1):80--93, 2013.
	
	\bibitem{AiIfrimTataru}
	A.~Ai, M.~Ifrim, and D.~Tataru.
	\newblock {Two dimensional gravity waves at low regularity II: Global
		solutions}.
	\newblock \arXiv{2009.11513}.
	
	\bibitem{Airy45}
	G.~B. Airy.
	\newblock Tides and waves.
	\newblock {\em Encycl. Metropolitana}, 5:291--396, 1845.
	
	\bibitem{AissioueneBristeauGodlewskiEtAl20}
	N.~A\"{\i}ssiouene, M.-O. Bristeau, E.~Godlewski, A.~Mangeney, C.~Par\'{e}s
	Madro\~{n}al, and J.~Sainte-Marie.
	\newblock A two-dimensional method for a family of dispersive shallow water
	models.
	\newblock {\em SMAI J. Comput. Math.}, 6:187--226, 2020.
	
	\bibitem{Alazard08}
	T.~Alazard.
	\newblock A minicourse on the low {M}ach number limit.
	\newblock {\em Discrete Contin. Dyn. Syst. Ser. S}, 1(3):365--404, 2008.
	
	\bibitem{AlazardBurqZuily13}
	T.~Alazard, N.~Burq, and C.~Zuily.
	\newblock The water-wave equations: from {Z}akharov to {E}uler.
	\newblock In {\em Studies in phase space analysis with applications to {PDE}s},
	volume~84 of {\em Progr. Nonlinear Differential Equations Appl.}, pages
	1--20. Birkh\"{a}user/Springer, New York, 2013.
	
	\bibitem{AlazardBurqZuily14}
	T.~Alazard, N.~Burq, and C.~Zuily.
	\newblock On the {C}auchy problem for gravity water waves.
	\newblock {\em Invent. Math.}, 198(1):71--163, 2014.
	
	\bibitem{AlazardBurqZuily16}
	T.~Alazard, N.~Burq, and C.~Zuily.
	\newblock Cauchy theory for the gravity water waves system with non-localized
	initial data.
	\newblock {\em Ann. Inst. H. Poincar\'{e} Anal. Non Lin\'{e}aire},
	33(2):337--395, 2016.
	
	\bibitem{AlinhacGerard}
	S.~Alinhac and P.~G{\'e}rard.
	\newblock {\em Op\'erateurs pseudo-diff\'erentiels et th\'eor\`eme de
		{N}ash-{M}oser}.
	\newblock Savoirs Actuels. InterEditions et Editions du CNRS, Paris, 1991.
	
	\bibitem{Alvarez-SamaniegoLannes08}
	B.~Alvarez-Samaniego and D.~Lannes.
	\newblock Large time existence for 3{D} water-waves and asymptotics.
	\newblock {\em Invent. Math.}, 171(3):485--541, 2008.
	
	\bibitem{Alvarez-SamaniegoLannes08a}
	B.~Alvarez-Samaniego and D.~Lannes.
	\newblock A {N}ash-{M}oser theorem for singular evolution equations.
	{A}pplication to the {S}erre and {G}reen-{N}aghdi equations.
	\newblock {\em Indiana Univ. Math. J.}, 57(1):97--131, 2008.
	
	\bibitem{Ambrose03}
	D.~M. Ambrose.
	\newblock Well-posedness of vortex sheets with surface tension.
	\newblock {\em SIAM J. Math. Anal.}, 35(1):211--244 (electronic), 2003.
	
	\bibitem{AmbroseBonaNicholls14}
	D.~M. Ambrose, J.~L. Bona, and D.~P. Nicholls.
	\newblock On ill-posedness of truncated series models for water waves.
	\newblock {\em Proc. R. Soc. Lond. Ser. A Math. Phys. Eng. Sci.},
	470(2166):20130849, 16, 2014.
	
	\bibitem{AmbroseCamassaMarzuolaEtAl}
	D.~M. Ambrose, R.~Camassa, J.~L. Marzuola, R.~M. McLaughlin, Q.~Robinson, and
	J.~Wilkening.
	\newblock {Numerical Algorithms for Water Waves with Background Flow over
		Obstacles and Topography}.
	\newblock \arXiv{2108.01786}.
	
	\bibitem{AmbroseMasmoudi07}
	D.~M. Ambrose and N.~Masmoudi.
	\newblock Well-posedness of 3{D} vortex sheets with surface tension.
	\newblock {\em Commun. Math. Sci.}, 5(2):391--430, 2007.
	
	\bibitem{Ambrosi00}
	D.~Ambrosi.
	\newblock Hamiltonian formulation for surface waves in a layered fluid.
	\newblock {\em Wave Motion}, 31(1):71--76, 2000.
	
	\bibitem{Amick84}
	C.~J. Amick.
	\newblock Regularity and uniqueness of solutions to the {B}oussinesq system of
	equations.
	\newblock {\em J. Differential Equations}, 54(2):231--247, 1984.
	
	\bibitem{AmickFraenkelToland82}
	C.~J. Amick, L.~E. Fraenkel, and J.~F. Toland.
	\newblock On the {S}tokes conjecture for the wave of extreme form.
	\newblock {\em Acta Math.}, 148:193--214, 1982.
	
	\bibitem{AmickToland81}
	C.~J. Amick and J.~F. Toland.
	\newblock On solitary water-waves of finite amplitude.
	\newblock {\em Arch. Rational Mech. Anal.}, 76(1):9--95, 1981.
	
	\bibitem{AndradeNachbin18}
	D.~Andrade and A.~Nachbin.
	\newblock A three-dimensional {D}irichlet-to-{N}eumann operator for water waves
	over topography.
	\newblock {\em J. Fluid Mech.}, 845:321--345, 2018.
	
	\bibitem{AntonopoulosDougalisMitsotakis}
	D.~C. Antonopoulos, V.~A. Dougalis, and D.~E. Mitsotakis.
	\newblock {On the well-posedness of the Galerkin semidiscretization of the
		periodic initial-value problem of the Serre equations}.
	\newblock \arXiv{2107.04403}.
	
	\bibitem{Armi86}
	L.~Armi.
	\newblock The hydraulics of two flowing layers with different densities.
	\newblock {\em J. Fluid Mech.}, 163:27--58, 1986.
	
	\bibitem{Asano87}
	K.~Asano.
	\newblock On the incompressible limit of the compressible {E}uler equation.
	\newblock {\em Japan J. Appl. Math.}, 4(3):455--488, 1987.
	
	\bibitem{AthanassoulisBelibassakis99}
	G.~A. Athanassoulis and K.~A. Belibassakis.
	\newblock A consistent coupled-mode theory for the propagation of
	small-amplitude water waves over variable bathymetry regions.
	\newblock {\em Journal of Fluid Mechanics}, 389:275--301, 1999.
	
	\bibitem{AthanassoulisPapoutsellis17}
	G.~A. Athanassoulis and C.~E. Papoutsellis.
	\newblock Exact semi-separation of variables in waveguides with non-planar
	boundaries.
	\newblock {\em Proc. A.}, 473(2201):20170017, 18, 2017.
	
	\bibitem{AudusseBristeauPelantiEtAl11}
	E.~Audusse, M.-O. Bristeau, M.~Pelanti, and J.~Sainte-Marie.
	\newblock Approximation of the hydrostatic {N}avier-{S}tokes system for density
	stratified flows by a multilayer model: kinetic interpretation and numerical
	solution.
	\newblock {\em J. Comput. Phys.}, 230(9):3453--3478, 2011.
	
	\bibitem{AudusseBristeauPerthameEtAl11}
	E.~Audusse, M.-O. Bristeau, B.~Perthame, and J.~Sainte-Marie.
	\newblock A multilayer {S}aint-{V}enant system with mass exchanges for shallow
	water flows. {D}erivation and numerical validation.
	\newblock {\em ESAIM Math. Model. Numer. Anal.}, 45(1):169--200, 2011.
	
	\bibitem{Babenko87}
	K.~I. Babenko.
	\newblock Some remarks on the theory of surface waves of finite amplitude.
	\newblock {\em Dokl. Akad. Nauk SSSR}, 294(5):1033--1037, 1987.
	
	\bibitem{BaeGranero-Belinchon}
	H.~Bae and R.~Granero-Belinchón.
	\newblock {Singularity formation for the Serre-Green-Naghdi equations and
		applications to abcd-Boussinesq systems}.
	\newblock \arXiv{2001.11937}.
	
	\bibitem{BahouriCheminDanchin11}
	H.~Bahouri, J.-Y. Chemin, and R.~Danchin.
	\newblock {\em Fourier analysis and nonlinear partial differential equations},
	volume 343.
	\newblock Springer, 2011.
	
	\bibitem{Baines88}
	P.~G. Baines.
	\newblock A general method for determining upstream effects in stratified flow
	of finite depth over long two-dimensional obstacles.
	\newblock {\em J. Fluid Mech.}, 188:1--22, 1988.
	
	\bibitem{Bambusi20}
	D.~Bambusi.
	\newblock Hamiltonian studies on counter-propagating water waves.
	\newblock {\em Water Waves}, pages 1--35, 2020.
	
	\bibitem{BardosBesse13}
	C.~Bardos and N.~Besse.
	\newblock The {C}auchy problem for the {V}lasov-{D}irac-{B}enney equation and
	related issues in fluid mechanics and semi-classical limits.
	\newblock {\em Kinet. Relat. Models}, 6(4):893--917, 2013.
	
	\bibitem{BardosBesse15}
	C.~Bardos and N.~Besse.
	\newblock Hamiltonian structure, fluid representation and stability for the
	{V}lasov-{D}irac-{B}enney equation.
	\newblock In {\em Hamiltonian partial differential equations and applications},
	volume~75 of {\em Fields Inst. Commun.}, pages 1--30. Fields Inst. Res. Math.
	Sci., Toronto, ON, 2015.
	
	\bibitem{BarrosChoi08}
	R.~Barros and W.~Choi.
	\newblock On the hyperbolicity of two-layer flows.
	\newblock In {\em Frontiers of applied and computational mathematics}, pages
	95--103. World Sci. Publ., Hackensack, NJ, 2008.
	
	\bibitem{Barthelemy04}
	E.~Barth{\'e}lemy.
	\newblock Nonlinear shallow water theories for coastal waves.
	\newblock {\em Surveys in Geophysics}, 25(3-4):315--337, 2004.
	
	\bibitem{BassiBonaventuraBustoEtAl20}
	C.~Bassi, L.~Bonaventura, S.~Busto, and M.~Dumbser.
	\newblock A hyperbolic reformulation of the {S}erre-{G}reen-{N}aghdi model for
	general bottom topographies.
	\newblock {\em Comput. \& Fluids}, 212:104716, 21, 2020.
	
	\bibitem{BazdenkovMorozovPogutse87}
	S.~Bazdenkov, N.~Morozov, and O.~Pogutse.
	\newblock Dispersive effects in two-dimensional hydrodynamics.
	\newblock In {\em Soviet Physics Doklady}, volume~32, page 262, 1987.
	\newblock In Russian.
	
	\bibitem{Beale77}
	J.~T. Beale.
	\newblock The existence of solitary water waves.
	\newblock {\em Comm. Pure Appl. Math.}, 30(4):373--389, 1977.
	
	\bibitem{BelibassakisAthanassoulis11}
	K.~A. Belibassakis and G.~A. Athanassoulis.
	\newblock A coupled-mode system with application to nonlinear water waves
	propagating in finite water depth and in variable bathymetry regions.
	\newblock {\em Coastal Engineering}, 58(4):337--350, 2011.
	
	\bibitem{BenjaminBridges97}
	T.~B. Benjamin and T.~J. Bridges.
	\newblock Reappraisal of the {K}elvin-{H}elmholtz problem. {I}. {H}amiltonian
	structure.
	\newblock {\em J. Fluid Mech.}, 333:301--325, 1997.
	
	\bibitem{BenjaminOlver82}
	T.~B. Benjamin and P.~J. Olver.
	\newblock Hamiltonian structure, symmetries and conservation laws for water
	waves.
	\newblock {\em J. Fluid Mech.}, 125:137--185, 1982.
	
	\bibitem{Benney73}
	D.~J. Benney.
	\newblock Some properties of long nonlinear waves.
	\newblock {\em Studies in Appl. Math.}, 52(1):45--50, 1973.
	
	\bibitem{Benton54}
	G.~S. Benton.
	\newblock The occurrence of critical flow and hydraulic jumps in a
	multi-layered fluid system.
	\newblock {\em Journal of Meteorology}, 11(2):139--150, 1954.
	
	\bibitem{Benzoni-GavageSerre07}
	S.~Benzoni-Gavage and D.~Serre.
	\newblock {\em Multidimensional hyperbolic partial differential equations.
		{F}irst-order systems and applications}.
	\newblock Oxford Mathematical Monographs. The Clarendon Press Oxford University
	Press, Oxford, 2007.
	
	\bibitem{Julia}
	J.~Bezanson, A.~Edelman, S.~Karpinski, and V.~B. Shah.
	\newblock Julia: a fresh approach to numerical computing.
	\newblock {\em SIAM Rev.}, 59(1):65--98, 2017.
	
	\bibitem{BonaChenSaut02}
	J.~L. Bona, M.~Chen, and J.-C. Saut.
	\newblock Boussinesq equations and other systems for small-amplitude long waves
	in nonlinear dispersive media. {I}. {D}erivation and linear theory.
	\newblock {\em J. Nonlinear Sci.}, 12(4):283--318, 2002.
	
	\bibitem{BonaChenSaut04}
	J.~L. Bona, M.~Chen, and J.-C. Saut.
	\newblock Boussinesq equations and other systems for small-amplitude long waves
	in nonlinear dispersive media. {II}. {T}he nonlinear theory.
	\newblock {\em Nonlinearity}, 17(3):925--952, 2004.
	
	\bibitem{BonaColinLannes05}
	J.~L. Bona, T.~Colin, and D.~Lannes.
	\newblock Long wave approximations for water waves.
	\newblock {\em Arch. Ration. Mech. Anal.}, 178(3):373--410, 2005.
	
	\bibitem{BonaLannesSaut08}
	J.~L. Bona, D.~Lannes, and J.-C. Saut.
	\newblock Asymptotic models for internal waves.
	\newblock {\em J. Math. Pures Appl. (9)}, 89(6):538--566, 2008.
	
	\bibitem{BonaSmith75}
	J.~L. Bona and R.~Smith.
	\newblock The initial-value problem for the {K}orteweg-de {V}ries equation.
	\newblock {\em Philos. Trans. Roy. Soc. London Ser. A}, 278(1287):555--601,
	1975.
	
	\bibitem{Bonnet-BenDhiaBristeauGodlewskiEtAl}
	A.-S. Bonnet-Ben~Dhia, M.-O. Bristeau, E.~Godlewski, S.~Imperiale, A.~Mangeney,
	and J.~Sainte-Marie.
	\newblock Pseudo-compressibility, dispersive model and acoustic waves in
	shallow water flows.
	\newblock preprint available at \url{https://hal.inria.fr/hal-02493518}.
	
	\bibitem{BoonkasameMilewski12}
	A.~Boonkasame and P.~Milewski.
	\newblock The stability of large-amplitude shallow interfacial non-{B}oussinesq
	flows.
	\newblock {\em Stud. Appl. Math.}, 128(1):40--58, 2012.
	
	\bibitem{BoonkasameMilewski14}
	A.~Boonkasame and P.~A. Milewski.
	\newblock A model for strongly nonlinear long interfacial waves with background
	shear.
	\newblock {\em Stud. Appl. Math.}, 133(2):182--213, 2014.
	
	\bibitem{Boussinesq72}
	J.~Boussinesq.
	\newblock Th{\'e}orie des ondes et des remous qui se propagent le long d'un
	canal rectangulaire horizontal, en communiquant au liquide contenu dans ce
	canal des vitesses sensiblement pareilles de la surface au fond.
	\newblock {\em J. Math. Pures Appl.}, 17(2):55--108, 1872.
	
	\bibitem{Boussinesq73}
	J.~Boussinesq.
	\newblock Addition au mémoire sur la théorie des ondes et des remous qui se
	propagent le long d'un canal rectangulaire, etc.
	\newblock {\em J. Math. Pures Appl.}, 17(2):47--52, 1873.
	
	\bibitem{Boussinesq77}
	J.~Boussinesq.
	\newblock Essai sur la théorie des eaux courantes.
	\newblock {\em Mém. présent. divers savants Acad. sci. Inst. Fr.}, 23:1--680,
	1877.
	
	\bibitem{BreschMetivier10}
	D.~Bresch and G.~M{\'e}tivier.
	\newblock Anelastic limits for {E}uler-type systems.
	\newblock {\em Appl. Math. Res. Express. AMRX}, 2010(2):119--141, 2010.
	
	\bibitem{BreschRenardy11}
	D.~Bresch and M.~Renardy.
	\newblock Well-posedness of two-layer shallow water flow between two horizontal
	rigid plates.
	\newblock {\em Nonlinearity}, 24:1081--1088, 2011.
	
	\bibitem{BridgesGrovesNicholls16}
	T.~J. Bridges, M.~D. Groves, and D.~P. Nicholls, editors.
	\newblock {\em Lectures on the theory of water waves}, volume 426 of {\em
		London Mathematical Society Lecture Note Series}.
	\newblock Cambridge University Press, Cambridge, 2016.
	\newblock Papers from the talks given at the Isaac Newton Institute for
	Mathematical Sciences, Cambridge, July--August, 2014.
	
	\bibitem{BristeauMangeneySainte-MarieEtAl15}
	M.-O. Bristeau, A.~Mangeney, J.~Sainte-Marie, and N.~Seguin.
	\newblock An energy-consistent depth-averaged {E}uler system: derivation and
	properties.
	\newblock {\em Discrete Contin. Dyn. Syst. Ser. B}, 20(4):961--988, 2015.
	
	\bibitem{BrowningKreiss82}
	G.~Browning and H.-O. Kreiss.
	\newblock Problems with different time scales for nonlinear partial
	differential equations.
	\newblock {\em SIAM J. Appl. Math.}, 42(4):704--718, 1982.
	
	\bibitem{BuffoniGrovesSunEtAl13}
	B.~Buffoni, M.~D. Groves, S.~M. Sun, and E.~Wahl\'{e}n.
	\newblock Existence and conditional energetic stability of three-dimensional
	fully localised solitary gravity-capillary water waves.
	\newblock {\em J. Differential Equations}, 254(3):1006--1096, 2013.
	
	\bibitem{BustoDumbserEscalanteEtAl21}
	S.~Busto, M.~Dumbser, C.~Escalante, N.~Favrie, and S.~Gavrilyuk.
	\newblock On {H}igh {O}rder {ADER} {D}iscontinuous {G}alerkin {S}chemes for
	{F}irst {O}rder {H}yperbolic {R}eformulations of {N}onlinear {D}ispersive
	{S}ystems.
	\newblock {\em J. Sci. Comput.}, 87(2):48, 2021.
	
	\bibitem{Byatt-Smith71}
	J.~G.~B. Byatt-Smith.
	\newblock An integral equation for unsteady surface waves and a comment on the
	{B}oussinesq equation.
	\newblock {\em J. Fluid Mech.}, 49:625--633, 1971.
	
	\bibitem{CamassaChenFalquiEtAl12}
	R.~Camassa, S.~Chen, G.~Falqui, G.~Ortenzi, and M.~Pedroni.
	\newblock An inertia `paradox' for incompressible stratified {E}uler fluids.
	\newblock {\em J. Fluid Mech.}, 695:330--340, 2012.
	
	\bibitem{CamassaChenFalquiEtAl13}
	R.~Camassa, S.~Chen, G.~Falqui, G.~Ortenzi, and M.~Pedroni.
	\newblock Effects of inertia and stratification in incompressible ideal fluids:
	pressure imbalances by rigid confinement.
	\newblock {\em J. Fluid Mech.}, 726:404--438, 2013.
	
	\bibitem{CamassaFalquiOrtenziEtAl19}
	R.~Camassa, G.~Falqui, G.~Ortenzi, M.~Pedroni, and C.~Thomson.
	\newblock Hydrodynamic models and confinement effects by horizontal boundaries.
	\newblock {\em J. Nonlinear Sci.}, 29(4):1445--1498, 2019.
	
	\bibitem{CamassaHolmLevermore96}
	R.~Camassa, D.~D. Holm, and C.~D. Levermore.
	\newblock Long-time effects of bottom topography in shallow water.
	\newblock {\em Phys. D}, 98(2-4):258--286, 1996.
	\newblock Nonlinear phenomena in ocean dynamics (Los Alamos, NM, 1995).
	
	\bibitem{CamassaHolmLevermore97}
	R.~Camassa, D.~D. Holm, and C.~D. Levermore.
	\newblock Long-time shallow-water equations with a varying bottom.
	\newblock {\em J. Fluid Mech.}, 349:173--189, 1997.
	
	\bibitem{CaoLiTiti16}
	C.~Cao, J.~Li, and E.~S. Titi.
	\newblock Global well-posedness of the three-dimensional primitive equations
	with only horizontal viscosity and diffusion.
	\newblock {\em Comm. Pure Appl. Math.}, 69(8):1492--1531, 2016.
	
	\bibitem{CaoLiTiti17}
	C.~Cao, J.~Li, and E.~S. Titi.
	\newblock Strong solutions to the 3{D} primitive equations with only horizontal
	dissipation: near {$H^1$} initial data.
	\newblock {\em J. Funct. Anal.}, 272(11):4606--4641, 2017.
	
	\bibitem{Capistrano-FilhoGallegoPazoto19}
	R.~A. Capistrano-Filho, F.~A. Gallego, and A.~F. Pazoto.
	\newblock On the well-posedness and large-time behavior of higher order
	{B}oussinesq system.
	\newblock {\em Nonlinearity}, 32(5):1852--1881, 2019.
	
	\bibitem{Carter18}
	J.~D. Carter.
	\newblock Bidirectional {W}hitham equations as models of waves on shallow
	water.
	\newblock {\em Wave Motion}, 82:51--61, 2018.
	
	\bibitem{CarterCienfuegos11}
	J.~D. Carter and R.~Cienfuegos.
	\newblock The kinematics and stability of solitary and cnoidal wave solutions
	of the {S}erre equations.
	\newblock {\em Eur. J. Mech. B Fluids}, 30(3):259--268, 2011.
	
	\bibitem{CastroCordobaFeffermanEtAl13}
	A.~Castro, D.~C\'{o}rdoba, C.~Fefferman, F.~Gancedo, and J.~G\'{o}mez-Serrano.
	\newblock Finite time singularities for the free boundary incompressible
	{E}uler equations.
	\newblock {\em Ann. of Math. (2)}, 178(3):1061--1134, 2013.
	
	\bibitem{CastroCordobaFeffermanEtAl12a}
	A.~Castro, D.~C\'{o}rdoba, C.~Fefferman, F.~Gancedo, and
	M.~L\'{o}pez-Fern\'{a}ndez.
	\newblock Rayleigh-{T}aylor breakdown for the {M}uskat problem with
	applications to water waves.
	\newblock {\em Ann. of Math. (2)}, 175(2):909--948, 2012.
	
	\bibitem{CastroCordobaFeffermanEtAl12}
	A.~Castro, D.~C\'{o}rdoba, C.~L. Fefferman, F.~Gancedo, and
	J.~G\'{o}mez-Serrano.
	\newblock Splash singularity for water waves.
	\newblock {\em Proc. Natl. Acad. Sci. USA}, 109(3):733--738, 2012.
	
	\bibitem{CastroLannes14}
	A.~Castro and D.~Lannes.
	\newblock Fully nonlinear long-waves models in presence of vorticity.
	\newblock {\em J. Fluid Mech.}, 759:642--675, 2014.
	
	\bibitem{CastroLannes15}
	A.~Castro and D.~Lannes.
	\newblock Well-posedness and shallow-water stability for a new {H}amiltonian
	formulation of the water waves equations with vorticity.
	\newblock {\em Indiana Univ. Math. J.}, 64(4):1169--1270, 2015.
	
	\bibitem{Cathala13}
	M.~Cathala.
	\newblock {\em Probl{\'e}matiques d’analyse num{\'e}rique et de
		mod{\'e}lisation pour {\'e}coulements de fluides environnementaux}.
	\newblock PhD thesis, Universit{\'e} Montpellier II, 2013.
	
	\bibitem{Cathala16}
	M.~Cathala.
	\newblock Asymptotic shallow water models with non smooth topographies.
	\newblock {\em Monatsh. Math.}, 179(3):325--353, 2016.
	
	\bibitem{Chandrasekhar61}
	S.~Chandrasekhar.
	\newblock {\em Hydrodynamic and hydromagnetic stability}.
	\newblock The International Series of Monographs on Physics. Clarendon Press,
	Oxford, 1961.
	
	\bibitem{ChenJin}
	R.~M. Chen and J.~Jin.
	\newblock {Global bifurcation of solitary waves to the Boussinesq $abcd$
		system}.
	\newblock \arXiv{2103.10812}.
	
	\bibitem{ChenWalsh}
	R.~M. Chen and S.~Walsh.
	\newblock Orbital stability of internal waves.
	\newblock \arXiv{2102.13590}.
	
	\bibitem{ChenWalsh16}
	R.~M. Chen and S.~Walsh.
	\newblock Continuous dependence on the density for stratified steady water
	waves.
	\newblock {\em Arch. Ration. Mech. Anal.}, 219(2):741--792, 2016.
	
	\bibitem{ChengJuSchochet18}
	B.~Cheng, Q.~Ju, and S.~Schochet.
	\newblock Three-scale singular limits of evolutionary {PDE}s.
	\newblock {\em Arch. Ration. Mech. Anal.}, 229(2):601--625, 2018.
	
	\bibitem{ChengCoutandShkoller08}
	C.-H.~A. Cheng, D.~Coutand, and S.~Shkoller.
	\newblock On the motion of vortex sheets with surface tension in
	three-dimensional {E}uler equations with vorticity.
	\newblock {\em Comm. Pure Appl. Math.}, 61(12):1715--1752, 2008.
	
	\bibitem{ChesnokovElGavrilyukEtAl17}
	A.~A. Chesnokov, G.~A. El, S.~L. Gavrilyuk, and M.~V. Pavlov.
	\newblock Stability of shear shallow water flows with free surface.
	\newblock {\em SIAM J. Appl. Math.}, 77(3):1068--1087, 2017.
	
	\bibitem{Choi00}
	W.~Choi.
	\newblock Modeling of strongly nonlinear internal gravity waves.
	\newblock In {\em Proceedings of 4thInternational Conference on Hydrodynamics,
		Yokohama, Japan}, pages 453--458, 2000.
	
	\bibitem{Choi19a}
	W.~Choi.
	\newblock Fifth-order nonlinear spectral model for surface gravity waves: From
	pseudo-spectral to spectral formulations (workshop on nonlinear water waves).
	\newblock {\em RIMS Kokyuroku}, 2109:47--60, 2019.
	
	\bibitem{Choi19}
	W.~Choi.
	\newblock On {R}ayleigh expansion for nonlinear long water waves.
	\newblock {\em J. Hydrodyn.}, 31(6):1115--1126, 2019.
	
	\bibitem{ChoiBarrosJo09}
	W.~Choi, R.~Barros, and T.-C. Jo.
	\newblock A regularized model for strongly nonlinear internal solitary waves.
	\newblock {\em J. Fluid Mech.}, 629:73--85, 2009.
	
	\bibitem{ChoiCamassa96}
	W.~Choi and R.~Camassa.
	\newblock Weakly nonlinear internal waves in a two-fluid system.
	\newblock {\em J. Fluid Mech.}, 313:83--103, 1996.
	
	\bibitem{ChoiCamassa99a}
	W.~Choi and R.~Camassa.
	\newblock Exact evolution equations for surface waves.
	\newblock {\em J. Eng. Mech.}, 125(7):756--760, 1999.
	
	\bibitem{ChoiCamassa99}
	W.~Choi and R.~Camassa.
	\newblock Fully nonlinear internal waves in a two-fluid system.
	\newblock {\em J. Fluid Mech.}, 396:1--36, 1999.
	
	\bibitem{ChumakovaMenzaqueMilewskiEtAl09b}
	L.~Chumakova, F.~E. Menzaque, P.~A. Milewski, R.~R. Rosales, E.~G. Tabak, and
	C.~V. Turner.
	\newblock Shear instability for stratified hydrostatic flows.
	\newblock {\em Comm. Pure Appl. Math.}, 62(2):183--197, 2009.
	
	\bibitem{ChumakovaMenzaqueMilewskiEtAl09}
	L.~Chumakova, F.~E. Menzaque, P.~A. Milewski, R.~R. Rosales, E.~G. Tabak, and
	C.~V. Turner.
	\newblock Stability properties and nonlinear mappings of two and three-layer
	stratified flows.
	\newblock {\em Stud. Appl. Math.}, 122(2):123--137, 2009.
	
	\bibitem{ChumakovaTabak10}
	L.~Chumakova and E.~G. Tabak.
	\newblock Simple waves do not avoid eigenvalue crossings.
	\newblock {\em Comm. Pure Appl. Math.}, 63(1):119--132, 2010.
	
	\bibitem{CienfuegosBarthelemyBonneton06}
	R.~Cienfuegos, E.~Barth\'{e}lemy, and P.~Bonneton.
	\newblock A fourth-order compact finite volume scheme for fully nonlinear and
	weakly dispersive {B}oussinesq-type equations. {I}. {M}odel development and
	analysis.
	\newblock {\em Internat. J. Numer. Methods Fluids}, 51(11):1217--1253, 2006.
	
	\bibitem{ClaassenJohnson18}
	K.~M. Claassen and M.~A. Johnson.
	\newblock Numerical bifurcation and spectral stability of wavetrains in
	bidirectional {W}hitham models.
	\newblock {\em Stud. Appl. Math.}, 141(2):205--246, 2018.
	
	\bibitem{ClamondDutykh12}
	D.~Clamond and D.~Dutykh.
	\newblock Practical use of variational principles for modeling water waves.
	\newblock {\em Phys. D}, 241(1):25--36, 2012.
	
	\bibitem{ClamondDutykh18}
	D.~Clamond and D.~Dutykh.
	\newblock Accurate fast computation of steady two-dimensional surface gravity
	waves in arbitrary depth.
	\newblock {\em J. Fluid Mech.}, 844:491--518, 2018.
	
	\bibitem{ClamondDutykhMitsotakis17}
	D.~Clamond, D.~Dutykh, and D.~Mitsotakis.
	\newblock Conservative modified {S}erre-{G}reen-{N}aghdi equations with
	improved dispersion characteristics.
	\newblock {\em Commun. Nonlinear Sci. Numer. Simul.}, 45:245--257, 2017.
	
	\bibitem{Cokelet77}
	E.~D. Cokelet.
	\newblock Steep gravity waves in water of arbitrary uniform depth.
	\newblock {\em Philos. Trans. Roy. Soc. London Ser. A}, 286(1335):183--230,
	1977.
	
	\bibitem{ColinIguchi20}
	M.~Colin and T.~Iguchi.
	\newblock Solitary wave solutions to the {I}sobe-{K}akinuma model for water
	waves.
	\newblock {\em Stud. Appl. Math.}, 145(1):52--80, 2020.
	
	\bibitem{Constantin11}
	A.~Constantin.
	\newblock {\em Nonlinear water waves with applications to wave-current
		interactions and tsunamis}, volume~81 of {\em CBMS-NSF Regional Conference
		Series in Applied Mathematics}.
	\newblock Society for Industrial and Applied Mathematics (SIAM), Philadelphia,
	PA, 2011.
	
	\bibitem{ConstantinStrauss04}
	A.~Constantin and W.~Strauss.
	\newblock Exact steady periodic water waves with vorticity.
	\newblock {\em Comm. Pure Appl. Math.}, 57(4):481--527, 2004.
	
	\bibitem{ConstantinVarvaruca11}
	A.~Constantin and E.~Varvaruca.
	\newblock Steady periodic water waves with constant vorticity: regularity and
	local bifurcation.
	\newblock {\em Arch. Ration. Mech. Anal.}, 199(1):33--67, 2011.
	
	\bibitem{CordobaFefferman18}
	D.~C\'{o}rdoba and C.~Fefferman.
	\newblock Water waves with or without surface tension.
	\newblock In {\em Handbook of mathematical analysis in mechanics of viscous
		fluids}, pages 1329--1349. Springer, Cham, 2018.
	
	\bibitem{CotterHolmPercival10}
	C.~J. Cotter, D.~D. Holm, and J.~R. Percival.
	\newblock The square root depth wave equations.
	\newblock {\em Proc. R. Soc. Lond. Ser. A Math. Phys. Eng. Sci.},
	466(2124):3621--3633, 2010.
	
	\bibitem{Craig85}
	W.~Craig.
	\newblock An existence theory for water waves and the {B}oussinesq and
	{K}orteweg-de {V}ries scaling limits.
	\newblock {\em Comm. Partial Differential Equations}, 10(8):787--1003, 1985.
	
	\bibitem{CraigGazeauLacaveEtAl18}
	W.~Craig, M.~Gazeau, C.~Lacave, and C.~Sulem.
	\newblock Bloch theory and spectral gaps for linearized water waves.
	\newblock {\em SIAM J. Math. Anal.}, 50(5):5477--5501, 2018.
	
	\bibitem{CraigGroves94}
	W.~Craig and M.~D. Groves.
	\newblock Hamiltonian long-wave approximations to the water-wave problem.
	\newblock {\em Wave Motion}, 19(4):367--389, 1994.
	
	\bibitem{CraigGroves00}
	W.~Craig and M.~D. Groves.
	\newblock Normal forms for wave motion in fluid interfaces.
	\newblock {\em Wave Motion}, 31(1):21--41, 2000.
	
	\bibitem{CraigGuyenneKalisch05}
	W.~Craig, P.~Guyenne, and H.~Kalisch.
	\newblock Hamiltonian long-wave expansions for free surfaces and interfaces.
	\newblock {\em Comm. Pure Appl. Math.}, 58(12):1587--1641, 2005.
	
	\bibitem{CraigGuyenneNichollsEtAl05}
	W.~Craig, P.~Guyenne, D.~P. Nicholls, and C.~Sulem.
	\newblock Hamiltonian long-wave expansions for water waves over a rough bottom.
	\newblock {\em Proc. R. Soc. Lond. Ser. A Math. Phys. Eng. Sci.},
	461(2055):839--873, 2005.
	
	\bibitem{CraigGuyenneSulem09}
	W.~Craig, P.~Guyenne, and C.~Sulem.
	\newblock Water waves over a random bottom.
	\newblock {\em J. Fluid Mech.}, 640:79--107, 2009.
	
	\bibitem{CraigGuyenneSulem10}
	W.~Craig, P.~Guyenne, and C.~Sulem.
	\newblock {C}oupling between internal and surface waves.
	\newblock {\em Natural Hazards}, 57(3):617--642, 2010.
	
	\bibitem{CraigGuyenneSulem12}
	W.~Craig, P.~Guyenne, and C.~Sulem.
	\newblock The surface signature of internal waves.
	\newblock {\em J. Fluid Mech.}, 710:277--303, 2012.
	
	\bibitem{CraigGuyenneSulem15}
	W.~Craig, P.~Guyenne, and C.~Sulem.
	\newblock Internal waves coupled to surface gravity waves in three dimensions.
	\newblock {\em Commun. Math. Sci}, 13:893--910, 2015.
	
	\bibitem{CraigLannesSulem12}
	W.~Craig, D.~Lannes, and C.~Sulem.
	\newblock Water waves over a rough bottom in the shallow water regime.
	\newblock {\em Ann. Inst. H. Poincar\'e Anal. Non Lin\'eaire}, 29(2):233--259,
	2012.
	
	\bibitem{CraigNicholls00}
	W.~Craig and D.~P. Nicholls.
	\newblock Travelling two and three dimensional capillary gravity water waves.
	\newblock {\em SIAM J. Math. Anal.}, 32(2):323--359, 2000.
	
	\bibitem{CraigSulem93}
	W.~Craig and C.~Sulem.
	\newblock Numerical simulation of gravity waves.
	\newblock {\em J. Comput. Phys.}, 108(1):73--83, 1993.
	
	\bibitem{CraigSulemSulem92}
	W.~Craig, C.~Sulem, and P.-L. Sulem.
	\newblock Nonlinear modulation of gravity waves: a rigorous approach.
	\newblock {\em Nonlinearity}, 5(2):497--522, 1992.
	
	\bibitem{CreedonDeconinckTrichtchenko21}
	R.~Creedon, B.~Deconinck, and O.~Trichtchenko.
	\newblock {High-Frequency Instabilities of a Boussinesq--Whitham System: A
		Perturbative Approach}.
	\newblock {\em Fluids}, 6(4):136, 2021.
	
	\bibitem{Dafermos10}
	C.~M. Dafermos.
	\newblock {\em Hyperbolic conservation laws in continuum physics}, volume 325
	of {\em Grundlehren der Mathematischen Wissenschaften [Fundamental Principles
		of Mathematical Sciences]}.
	\newblock Springer-Verlag, Berlin, third edition, 2010.
	
	\bibitem{Darrigol03}
	O.~Darrigol.
	\newblock The spirited horse, the engineer, and the mathematician: water waves
	in nineteenth-century hydrodynamics.
	\newblock {\em Arch. Hist. Exact Sci.}, 58(1):21--95, 2003.
	
	\bibitem{DeconinckTrichtchenko17}
	B.~Deconinck and O.~Trichtchenko.
	\newblock High-frequency instabilities of small-amplitude solutions of
	{H}amiltonian {PDE}s.
	\newblock {\em Discrete Contin. Dyn. Syst.}, 37(3):1323--1358, 2017.
	
	\bibitem{DegondTang11}
	P.~Degond and M.~Tang.
	\newblock All speed scheme for the low {M}ach number limit of the isentropic
	{E}uler equations.
	\newblock {\em Commun. Comput. Phys.}, 10(1):1--31, 2011.
	
	\bibitem{Delort18}
	J.-M. Delort.
	\newblock Long time existence results for solutions of water waves equations.
	\newblock In {\em Proceedings of the {I}nternational {C}ongress of
		{M}athematicians---{R}io de {J}aneiro 2018. {V}ol. {III}. {I}nvited
		lectures}, pages 2241--2260. World Sci. Publ., Hackensack, NJ, 2018.
	
	\bibitem{DenekeTesfahunTemesgen}
	T.~Deneke, A.~Tesfahun, and T.~Temesgen.
	\newblock {Dispersive estimates for linearized water wave type equations in
		$\mathbb{R}^d$}.
	\newblock \arXiv{2106.02717}.
	
	\bibitem{DesjardinsLannesSaut19}
	B.~Desjardins, D.~Lannes, and J.-C. Saut.
	\newblock Normal mode decomposition and dispersive and nonlinear mixing in
	stratified fluids.
	\newblock {\em Water Waves}, pages 1--40, 2020.
	
	\bibitem{DiasMilewski10}
	F.~Dias and P.~Milewski.
	\newblock On the fully-nonlinear shallow-water generalized {S}erre equations.
	\newblock {\em Phys. Lett., A}, 374(8):1049--1053, 2010.
	
	\bibitem{Dingemans97}
	M.~W. Dingemans.
	\newblock {\em Water Wave Propagation Over Uneven Bottoms: Non-linear wave
		propagation}, volume~13 of {\em Advanced Series on Ocean Engineering}.
	\newblock World Scientific, Cornell Univ., Hollister Hall, 1997.
	
	\bibitem{Dinvay}
	E.~Dinvay.
	\newblock {Travelling waves in the Boussinesq type systems}.
	\newblock \arXiv{2011.09543}.
	
	\bibitem{Dinvay19}
	E.~Dinvay.
	\newblock On well-posedness of a dispersive system of the
	{W}hitham-{B}oussinesq type.
	\newblock {\em Appl. Math. Lett.}, 88:13--20, 2019.
	
	\bibitem{DinvayDutykhKalisch19}
	E.~Dinvay, D.~Dutykh, and H.~Kalisch.
	\newblock A comparative study of bi-directional {W}hitham systems.
	\newblock {\em Appl. Numer. Math.}, 141:248--262, 2019.
	
	\bibitem{DinvayKuznetsov19}
	E.~Dinvay and N.~Kuznetsov.
	\newblock Modified {B}abenko's equation for periodic gravity waves on water of
	finite depth.
	\newblock {\em Quart. J. Mech. Appl. Math.}, 72(4):415--428, 2019.
	
	\bibitem{DinvayNilsson21}
	E.~Dinvay and D.~Nilsson.
	\newblock Solitary wave solutions of a {W}hitham-{B}oussinesq system.
	\newblock {\em Nonlinear Anal. Real World Appl.}, 60:103280, 2021.
	
	\bibitem{DinvaySelbergTesfahun20}
	E.~Dinvay, S.~Selberg, and A.~Tesfahun.
	\newblock Well-{P}osedness for a {D}ispersive {S}ystem of the
	{W}hitham--{B}oussinesq {T}ype.
	\newblock {\em SIAM J. Math. Anal.}, 52(3):2353--2382, 2020.
	
	\bibitem{DommermuthYue87}
	D.~G. Dommermuth and D.~K. Yue.
	\newblock A high-order spectral method for the study of nonlinear gravity
	waves.
	\newblock {\em J. Fluid Mech.}, 184:267--288, 1987.
	
	\bibitem{DorodnitsynKaptsovMeleshko}
	V.~A. Dorodnitsyn, E.~I. Kaptsov, and S.~V. Meleshko.
	\newblock {Symmetries, conservation laws, invariant solutions and difference
		schemes of the one-dimensional Green-Naghdi equations}.
	\newblock \arXiv{2008.12852}.
	
	\bibitem{DougalisDuranSaridaki}
	V.~A. Dougalis, A.~Duran, and L.~Saridaki.
	\newblock {Notes on numerical analysis and solitary wave solutions of
		Boussinesq/Boussinesq systems for internal waves}.
	\newblock \arXiv{2012.07992}.
	
	\bibitem{DougalisMitsotakis08}
	V.~A. Dougalis and D.~E. Mitsotakis.
	\newblock Theory and numerical analysis of {B}oussinesq systems: a review.
	\newblock In {\em Effective computational methods for wave propagation},
	volume~5 of {\em Numer. Insights}, pages 63--110. Chapman \& Hall/CRC, Boca
	Raton, FL, 2008.
	
	\bibitem{Duchene10}
	V.~Duchêne.
	\newblock Asymptotic shallow water models for internal waves in a two-fluid
	system with a free surface.
	\newblock {\em SIAM J. Math. Anal.}, 42(5):2229--2260, 2010.
	
	\bibitem{Duchene12}
	V.~Duchêne.
	\newblock Boussinesq/{B}oussinesq systems for internal waves with a free
	surface, and the {K}d{V} approximation.
	\newblock {\em ESAIM Math. Model. Numer. Anal.}, 46(1):145--185, 2012.
	
	\bibitem{Duchene13}
	V.~Duchêne.
	\newblock A note on the well-posedness of the one-dimensional multilayer
	shallow water model.
	\newblock hal preprint:00922045, 2013.
	
	\bibitem{Duchene14}
	V.~Duchêne.
	\newblock Decoupled and unidirectional asymptotic models for the propagation of
	internal waves.
	\newblock {\em Math. Models Methods Appl. Sci.}, 24(1):1--65, 2014.
	
	\bibitem{Duchene14a}
	V.~Duchêne.
	\newblock On the rigid-lid approximation for two shallow layers of immiscible
	fluids with small density contrast.
	\newblock {\em J. Nonlinear Sci.}, 24(4):579--632, 2014.
	
	\bibitem{Duchene16}
	V.~Duchêne.
	\newblock The multilayer shallow water system in the limit of small density
	contrast.
	\newblock {\em Asymptot. Anal.}, 98(3):189--235, 2016.
	
	\bibitem{Duchene19}
	V.~Duchêne.
	\newblock Rigorous justification of the {F}avrie-{G}avrilyuk approximation to
	the {S}erre-{G}reen-{N}aghdi model.
	\newblock {\em Nonlinearity}, 32(10):3772--3797, 2019.
	
	\bibitem{DucheneIguchi}
	V.~Duchêne and T.~Iguchi.
	\newblock {A mathematical analysis of the Kakinuma model for interfacial
		gravity waves. Part I: Structures and well-posedness}.
	\newblock \arXiv{2103.12392}.
	
	\bibitem{DucheneIguchia}
	V.~Duchêne and T.~Iguchi.
	\newblock {A mathematical analysis of the Kakinuma model for interfacial
		gravity waves. Part II: Justification as a shallow-water approximation}.
	\newblock In preparation.
	
	\bibitem{DucheneIguchi20}
	V.~Duchêne and T.~Iguchi.
	\newblock A hamiltonian structure of the isobe--kakinuma model for water waves.
	\newblock {\em Water Waves}, 3:1--19, 2020.
	
	\bibitem{DucheneIsrawi18}
	V.~Duchêne and S.~Israwi.
	\newblock Well-posedness of the {G}reen-{N}aghdi and {B}oussinesq-{P}eregrine
	systems.
	\newblock {\em Ann. Math. Blaise Pascal}, 25(1):21--74, 2018.
	
	\bibitem{DucheneIsrawiTalhouk14}
	V.~Duchêne, S.~Israwi, and R.~Talhouk.
	\newblock Shallow water asymptotic models for the propagation of internal
	waves.
	\newblock {\em Discrete Contin. Dyn. Syst. Ser. S}, 7(2):239--269, 2014.
	
	\bibitem{DucheneIsrawiTalhouk16}
	V.~Duchêne, S.~Israwi, and R.~Talhouk.
	\newblock A new class of two-layer {G}reen-{N}aghdi systems with improved
	frequency dispersion.
	\newblock {\em Stud. Appl. Math.}, 137(3):356--415, 2016.
	
	\bibitem{DucheneKlein}
	V.~Duchêne and C.~Klein.
	\newblock {Numerical study of the Serre-Green-Naghdi equations and a fully
		dispersive counterpart}.
	\newblock \arXiv{2005.13234}.
	
	\bibitem{DucheneNilssonWahlen18}
	V.~Duchêne, D.~Nilsson, and E.~Wahlén.
	\newblock Solitary {W}ave {S}olutions to a {C}lass of {M}odified
	{G}reen--{N}aghdi {S}ystems.
	\newblock {\em J. Math. Fluid Mech.}, 20(3):1059--1091, 2018.
	
	\bibitem{DuranMarche15}
	A.~Duran and F.~Marche.
	\newblock Discontinuous-{G}alerkin discretization of a new class of
	{G}reen-{N}aghdi equations.
	\newblock {\em Commun. Comput. Phys.}, 17(3):721--760, 2015.
	
	\bibitem{DutykhClamond14}
	D.~Dutykh and D.~Clamond.
	\newblock Efficient computation of steady solitary gravity waves.
	\newblock {\em Wave Motion}, 51(1):86--99, 2014.
	
	\bibitem{DutykhClamondMilewskietal13}
	D.~Dutykh, D.~Clamond, P.~Milewski, and D.~Mitsotakis.
	\newblock Finite volume and pseudo-spectral schemes for the fully nonlinear
	1{D} {S}erre equations.
	\newblock {\em European J. Appl. Math.}, 24(5):761--787, 2013.
	
	\bibitem{DyachenkoKuznetsovSpectorEtAl96}
	A.~I. Dyachenko, E.~A. Kuznetsov, M.~Spector, and V.~E. Zakharov.
	\newblock Analytical description of the free surface dynamics of an ideal fluid
	(canonical formalism and conformal mapping).
	\newblock {\em Phys. Lett. A}, 221(1-2):73--79, 1996.
	
	\bibitem{Ebin88}
	D.~G. Ebin.
	\newblock Ill-posedness of the {R}ayleigh-{T}aylor and {H}elmholtz problems for
	incompressible fluids.
	\newblock {\em Comm. Partial Differential Equations}, 13(10):1265--1295, 1988.
	
	\bibitem{EhrnstromGrovesWahlen12}
	M.~Ehrnstr{\"o}m, M.~D. Groves, and E.~Wahl{\'e}n.
	\newblock On the existence and stability of solitary-wave solutions to a class
	of evolution equations of {W}hitham type.
	\newblock {\em Nonlinearity}, 25(10):2903--2936, 2012.
	
	\bibitem{EhrnstroemJohnsonClaassen19}
	M.~Ehrnstr\"{o}m, M.~A. Johnson, and K.~M. Claassen.
	\newblock Existence of a highest wave in a fully dispersive two-way shallow
	water model.
	\newblock {\em Arch. Ration. Mech. Anal.}, 231(3):1635--1673, 2019.
	
	\bibitem{EhrnstromWahlen19}
	M.~Ehrnstr\"{o}m and E.~Wahl\'{e}n.
	\newblock On {W}hitham's conjecture of a highest cusped wave for a nonlocal
	dispersive equation.
	\newblock {\em Ann. Inst. H. Poincar\'{e} Anal. Non Lin\'{e}aire},
	36(6):1603--1637, 2019.
	
	\bibitem{ElGrimshawSmyth06}
	G.~A. El, R.~H.~J. Grimshaw, and N.~F. Smyth.
	\newblock Unsteady undular bores in fully nonlinear shallow-water theory.
	\newblock {\em Phys. Fluids}, 18(2):027104, 17, 2006.
	
	\bibitem{ElHoeferShearer16}
	G.~A. El, M.~A. Hoefer, and M.~Shearer.
	\newblock Expansion shock waves in regularized shallow-water theory.
	\newblock {\em Proc. A.}, 472(2189):20160141, 10, 2016.
	
	\bibitem{EmeraldIII}
	L.~Emerald.
	\newblock {Local well-posedness result for a class of non-local quasi-linear
		systems and its application to the Whitham-Boussinesq systems}.
	\newblock In preparation.
	
	\bibitem{EmeraldII}
	L.~Emerald.
	\newblock {Rigorous derivation of the Whitham equations from the water waves
		equations in the shallow water regime}.
	\newblock \arXiv{2101.02940}.
	
	\bibitem{EmeraldI}
	L.~Emerald.
	\newblock Rigorous derivation from the water waves equations of some full
	dispersion shallow water models.
	\newblock {\em SIAM J. Math. Anal.}, 53(4):3772--3800, 2021.
	
	\bibitem{EscalanteDumbserCastro19}
	C.~Escalante, M.~Dumbser, and M.~J. Castro.
	\newblock An efficient hyperbolic relaxation system for dispersive
	non-hydrostatic water waves and its solution with high order discontinuous
	{G}alerkin schemes.
	\newblock {\em J. Comput. Phys.}, 394:385--416, 2019.
	
	\bibitem{EscalanteMoralesdeLuna20}
	C.~Escalante and T.~Morales~de Luna.
	\newblock A general non-hydrostatic hyperbolic formulation for {B}oussinesq
	dispersive shallow flows and its numerical approximation.
	\newblock {\em J. Sci. Comput.}, 83(3):Paper No. 62, 37, 2020.
	
	\bibitem{FavrieGavrilyuk17}
	N.~Favrie and S.~Gavrilyuk.
	\newblock A rapid numerical method for solving {S}erre-{G}reen-{N}aghdi
	equations describing long free surface gravity waves.
	\newblock {\em Nonlinearity}, 30(7):2718--2736, 2017.
	
	\bibitem{FedotovaKhakimzyanovDutykh14}
	Z.~I. Fedotova, G.~S. Khakimzyanov, and D.~Dutykh.
	\newblock Energy equation for certain approximate models of long-wave
	hydrodynamics.
	\newblock {\em Russian J. Numer. Anal. Math. Modelling}, 29(3):167--178, 2014.
	
	\bibitem{FeolaGiuliani}
	R.~Feola and F.~Giuliani.
	\newblock Quasi-periodic traveling waves on an infinitely deep fluid under
	gravity.
	\newblock \arXiv{2005.08280}.
	
	\bibitem{Fernandez-NietoParisotPenelEtAl18}
	E.~D. Fern\'{a}ndez-Nieto, M.~Parisot, Y.~Penel, and J.~Sainte-Marie.
	\newblock A hierarchy of dispersive layer-averaged approximations of {E}uler
	equations for free surface flows.
	\newblock {\em Commun. Math. Sci.}, 16(5):1169--1202, 2018.
	
	\bibitem{FilippiniBellecColinEtAl15}
	A.~G. Filippini, S.~Bellec, M.~Colin, and M.~Ricchiuto.
	\newblock On the nonlinear behaviour of {B}oussinesq type models:
	{A}mplitude-velocity vs amplitude-flux forms.
	\newblock {\em Coast. Eng.}, 99:109--123, 2015.
	
	\bibitem{FouqueGarnierNachbin04}
	J.-P. Fouque, J.~Garnier, and A.~Nachbin.
	\newblock Shock structure due to stochastic forcing and the time reversal of
	nonlinear waves.
	\newblock {\em Phys. D}, 195(3-4):324--346, 2004.
	
	\bibitem{FriedrichsHyers54}
	K.~O. Friedrichs and D.~H. Hyers.
	\newblock The existence of solitary waves.
	\newblock {\em Comm. Pure Appl. Math.}, 7:517--550, 1954.
	
	\bibitem{Frings12}
	J.~T. Frings.
	\newblock {\em An adaptive multilayer model for density-layered shallow water
		flows}.
	\newblock PhD thesis, Aachen University, 2012.
	
	\bibitem{FujiwaraIguchi15}
	H.~Fujiwara and T.~Iguchi.
	\newblock A shallow water approximation for water waves over a moving bottom.
	\newblock In {\em Nonlinear dynamics in partial differential equations},
	volume~64 of {\em Adv. Stud. Pure Math.}, pages 77--88. Math. Soc. Japan,
	Tokyo, 2015.
	
	\bibitem{Gallagher05}
	I.~Gallagher.
	\newblock R\'esultats r\'ecents sur la limite incompressible.
	\newblock {\em Ast\'erisque}, 299:Exp. No. 926, vii, 29--57, 2005.
	\newblock S{\'e}minaire Bourbaki. Vol. 2003/2004.
	
	\bibitem{Gallay}
	T.~Gallay.
	\newblock {Stability of Vortices in Ideal Fluids : the Legacy of Kelvin and
		Rayleigh}.
	\newblock \arXiv{1901.02815}.
	
	\bibitem{GallaySmets19}
	T.~Gallay and D.~Smets.
	\newblock On the linear stability of vortex columns in the energy space.
	\newblock {\em J. Math. Fluid Mech.}, 21(4):Paper No. 48, 27, 2019.
	
	\bibitem{Gavrilyuk11}
	S.~Gavrilyuk.
	\newblock Multiphase flow modeling via {H}amilton's principle.
	\newblock In {\em Variational models and methods in solid and fluid mechanics},
	volume 535 of {\em CISM Courses and Lect.}, pages 163--210.
	SpringerWienNewYork, Vienna, 2011.
	
	\bibitem{GavrilyukGouin99}
	S.~Gavrilyuk and H.~Gouin.
	\newblock A new form of governing equations of fluids arising from {H}amilton's
	principle.
	\newblock {\em Internat. J. Engrg. Sci.}, 37(12):1495--1520, 1999.
	
	\bibitem{GavrilyukKalischKhorsand15}
	S.~Gavrilyuk, H.~Kalisch, and Z.~Khorsand.
	\newblock A kinematic conservation law in free surface flow.
	\newblock {\em Nonlinearity}, 28(6):1805--1821, 2015.
	
	\bibitem{GavrilyukNkongaShyueEtAl20}
	S.~Gavrilyuk, B.~Nkonga, K.-M. Shyue, and L.~Truskinovsky.
	\newblock Stationary shock-like transition fronts in dispersive systems.
	\newblock {\em Nonlinearity}, 33(10):5477--5509, 2020.
	
	\bibitem{GavrilyukTeshukov01}
	S.~L. Gavrilyuk and V.~M. Teshukov.
	\newblock Generalized vorticity for bubbly liquid and dispersive shallow water
	equations.
	\newblock {\em Contin. Mech. Thermodyn.}, 13(6):365--382, 2001.
	
	\bibitem{GentMcwilliams90}
	P.~R. Gent and J.~C. Mcwilliams.
	\newblock Isopycnal mixing in ocean circulation models.
	\newblock {\em J. Phys. Oceanogr.}, 20(1):150--155, 1990.
	
	\bibitem{GentWillebrandMcDougallEtAl95}
	P.~R. Gent, J.~Willebrand, T.~J. McDougall, and J.~C. McWilliams.
	\newblock Parameterizing eddy-induced tracer transports in ocean circulation
	models.
	\newblock {\em Journal of Physical Oceanography}, 25(4):463--474, 1995.
	
	\bibitem{GeurtsHolm06}
	B.~J. Geurts and D.~D. Holm.
	\newblock Leray and {LANS}-{$\alpha$} modelling of turbulent mixing.
	\newblock {\em J. Turbul.}, 7:Paper 10, 33, 2006.
	
	\bibitem{Gill82}
	A.~E. Gill.
	\newblock {\em Atmosphere-ocean dynamics}, volume~30 of {\em International
		geophysics series}.
	\newblock Academic Press, 1982.
	
	\bibitem{GreenNaghdi76}
	A.~E. Green and P.~M. Naghdi.
	\newblock A derivation of equations for wave propagation in water of variable
	depth.
	\newblock {\em J. Fluid Mech.}, 78(02):237--246, 1976.
	
	\bibitem{GrovesSun08}
	M.~D. Groves and S.-M. Sun.
	\newblock Fully localised solitary-wave solutions of the three-dimensional
	gravity-capillary water-wave problem.
	\newblock {\em Arch. Ration. Mech. Anal.}, 188(1):1--91, 2008.
	
	\bibitem{GuermondPopovTovarEtAl}
	J.-L. Guermond, B.~Popov, E.~Tovar, and C.~Kees.
	\newblock Hyperbolic relaxation technique for solving the dispersive {S}erre
	equations with topography.
	\newblock \arXiv{2103.01286}.
	
	\bibitem{GuermondPopovTovarEtAl19}
	J.-L. Guermond, B.~Popov, E.~Tovar, and C.~Kees.
	\newblock Robust explicit relaxation technique for solving the {G}reen-{N}aghdi
	equations.
	\newblock {\em J. Comput. Phys.}, 399:108917, 17, 2019.
	
	\bibitem{Guyenne19}
	P.~Guyenne.
	\newblock {HOS} simulations of nonlinear water waves in complex media.
	\newblock In D.~Henry, K.~Kalimeris, E.~Părău, J.~M. Vanden-Broeck, and
	E.~Wahlén, editors, {\em {N}onlinear {W}ater {W}aves}, Tutorials, Schools,
	and Workshops in the Mathematical Sciences, pages 53--69. Birkhäuser, Cham,
	2019.
	
	\bibitem{GuyenneLannesSaut10}
	P.~Guyenne, D.~Lannes, and J.-C. Saut.
	\newblock Well-posedness of the {C}auchy problem for models of large amplitude
	internal waves.
	\newblock {\em Nonlinearity}, 23(2):237--275, 2010.
	
	\bibitem{Hamilton77}
	J.~Hamilton.
	\newblock Differential equations for long-period gravity waves on fluid of
	rapidly varying depth.
	\newblock {\em J. Fluid Mech.}, 83(2):289--310, 1977.
	
	\bibitem{HaziotHurStraussEtAl}
	S.~V. Haziot, V.~M. Hur, W.~Strauss, J.~F. Toland, E.~Wahlén, S.~Walsh, and
	M.~H. Wheeler.
	\newblock Traveling water waves -- the ebb and flow of two centuries.
	\newblock \arXiv{2109.09208}.
	
	\bibitem{HechtHolmPetersenEtAl08}
	M.~W. Hecht, D.~D. Holm, M.~R. Petersen, and B.~A. Wingate.
	\newblock Implementation of the {LANS}-{$\alpha$} turbulence model in a
	primitive equation ocean model.
	\newblock {\em J. Comput. Phys.}, 227(11):5691--5716, 2008.
	
	\bibitem{HelfrichMelville06}
	K.~R. Helfrich and W.~K. Melville.
	\newblock Long nonlinear internal waves.
	\newblock In {\em Annual review of fluid mechanics. {V}ol. 38}, pages 395--425.
	2006.
	
	\bibitem{Holm88}
	D.~D. Holm.
	\newblock Hamiltonian structure for two-dimensional hydrodynamics with
	nonlinear dispersion.
	\newblock {\em Phys. Fluids}, 31(8):2371--2373, 1988.
	
	\bibitem{HolmLong89}
	D.~D. Holm and B.~Long.
	\newblock Lyapunov stability of ideal stratified fluid equilibria in
	hydrostatic balance.
	\newblock {\em Nonlinearity}, 2(1):23--35, 1989.
	
	\bibitem{Howard61}
	L.~N. Howard.
	\newblock Note on a paper of {J}ohn {W}. {M}iles.
	\newblock {\em J. Fluid Mech.}, 10:509--512, 1961.
	
	\bibitem{HughesKatoMarsden76}
	T.~J.~R. Hughes, T.~Kato, and J.~E. Marsden.
	\newblock Well-posed quasi-linear second-order hyperbolic systems with
	applications to nonlinear elastodynamics and general relativity.
	\newblock {\em Arch. Rational Mech. Anal.}, 63(3):273--294, 1976.
	
	\bibitem{HunterIfrimTataru16}
	J.~K. Hunter, M.~Ifrim, and D.~Tataru.
	\newblock Two dimensional water waves in holomorphic coordinates.
	\newblock {\em Comm. Math. Phys.}, 346(2):483--552, 2016.
	
	\bibitem{Hur17}
	V.~M. Hur.
	\newblock Wave breaking in the {W}hitham equation.
	\newblock {\em Adv. Math.}, 317:410--437, 2017.
	
	\bibitem{HurPandey19}
	V.~M. Hur and A.~K. Pandey.
	\newblock Modulational instability in a full-dispersion shallow water model.
	\newblock {\em Stud. Appl. Math.}, 142(1):3--47, 2019.
	
	\bibitem{HurTao18}
	V.~M. Hur and L.~Tao.
	\newblock Wave breaking in a shallow water model.
	\newblock {\em SIAM J. Math. Anal.}, 50(1):354--380, 2018.
	
	\bibitem{IfrimTataru16}
	M.~Ifrim and D.~Tataru.
	\newblock Two dimensional water waves in holomorphic coordinates {II}: {G}lobal
	solutions.
	\newblock {\em Bull. Soc. Math. France}, 144(2):369--394, 2016.
	
	\bibitem{IfrimTataru17}
	M.~Ifrim and D.~Tataru.
	\newblock The lifespan of small data solutions in two dimensional capillary
	water waves.
	\newblock {\em Arch. Ration. Mech. Anal.}, 225(3):1279--1346, 2017.
	
	\bibitem{Iguchi09}
	T.~Iguchi.
	\newblock A shallow water approximation for water waves.
	\newblock {\em J. Math. Kyoto Univ.}, 49(1):13--55, 2009.
	
	\bibitem{Iguchi18}
	T.~Iguchi.
	\newblock Isobe-{K}akinuma model for water waves as a higher order shallow
	water approximation.
	\newblock {\em J. Differential Equations}, 265(3):935--962, 2018.
	
	\bibitem{Iguchi18a}
	T.~Iguchi.
	\newblock A mathematical justification of the {I}sobe-{K}akinuma model for
	water waves with and without bottom topography.
	\newblock {\em J. Math. Fluid Mech.}, 20(4):1985--2018, 2018.
	
	\bibitem{IguchiTanakaTani97}
	T.~Iguchi, N.~Tanaka, and A.~Tani.
	\newblock On the two-phase free boundary problem for two-dimensional water
	waves.
	\newblock {\em Math. Ann.}, 309(2):199--223, 1997.
	
	\bibitem{Ionescu-Kruse12}
	D.~Ionescu-Kruse.
	\newblock Variational derivation of the {G}reen-{N}aghdi shallow-water
	equations.
	\newblock {\em J. Nonlinear Math. Phys.}, 19(suppl. 1):1240001, 12, 2012.
	
	\bibitem{Iselin39}
	C.~O. Iselin.
	\newblock The influence of vertical and lateral turbulence on the
	characteristics of the waters at mid-depths.
	\newblock {\em Eos, Trans. AGU}, 20(3):414--417, 1939.
	
	\bibitem{Isobe94}
	M.~Isobe.
	\newblock A proposal on a nonlinear gentle slope wave equation.
	\newblock {\em Proc. Coast. Eng. Jpn. Soc. Civ. Eng.}, 41:1--5, 1994.
	\newblock [Japanese].
	
	\bibitem{Isobe94a}
	M.~Isobe.
	\newblock Time-dependent mild-slope equations for random waves.
	\newblock In {\em Proceedings of 24th International Conference on Coastal
		Engineering}, pages 285--299. ASCE, 1994.
	
	\bibitem{Israwi11}
	S.~Israwi.
	\newblock Large time existence for 1{D} {G}reen-{N}aghdi equations.
	\newblock {\em Nonlinear Analysis: Theory, Methods \& Applications},
	74(1):81--93, 2011.
	
	\bibitem{Jackson04}
	C.~R. Jackson.
	\newblock An atlas of internal solitary-like waves and their properties, 2004.
	\newblock \url{http://www.internalwaveatlas.com/Atlas2_index.html}.
	
	\bibitem{Jin12}
	S.~Jin.
	\newblock Asymptotic preserving ({AP}) schemes for multiscale kinetic and
	hyperbolic equations: a review.
	\newblock {\em Riv. Math. Univ. Parma (N.S.)}, 3(2):177--216, 2012.
	
	\bibitem{JoChoi02}
	T.-C. Jo and W.~Choi.
	\newblock Dynamics of strongly nonlinear internal solitary waves in shallow
	water.
	\newblock {\em Stud. Appl. Math.}, 109(3):205--227, 2002.
	
	\bibitem{Kakinuma00}
	T.~Kakinuma.
	\newblock [title in japanese].
	\newblock {\em Proc. Coast. Eng. Jpn. Soc. Civ. Eng.}, 47:1--5, 2000.
	\newblock [Japanese].
	
	\bibitem{Kakinuma01}
	T.~Kakinuma.
	\newblock A set of fully nonlinear equations for surface and internal gravity
	waves.
	\newblock In {\em Coastal Engineering V: Computer Modelling of Seas and Coastal
		Regions}, pages 225--234. WIT Press, 2001.
	
	\bibitem{Kakinuma03}
	T.~Kakinuma.
	\newblock A nonlinear numerical model for surface and internal waves shoaling
	on a permeable beach.
	\newblock In {\em Coastal engineering VI: Computer Modelling and Experimental
		Measurements of Seas and Coastal Regions}, pages 227--236. WIT Press, 2003.
	
	\bibitem{KalischPilod19}
	H.~Kalisch and D.~Pilod.
	\newblock On the local well-posedness for a full-dispersion {B}oussinesq system
	with surface tension.
	\newblock {\em Proc. Amer. Math. Soc.}, 147(6):2545--2559, 2019.
	
	\bibitem{KamotskiLebeau05}
	V.~Kamotski and G.~Lebeau.
	\newblock On 2{D} {R}ayleigh-{T}aylor instabilities.
	\newblock {\em Asymptot. Anal.}, 42(1-2):1--27, 2005.
	
	\bibitem{Kano86}
	T.~Kano.
	\newblock Une th\'{e}orie trois-dimensionnelle des ondes de surface de l'eau et
	le d\'{e}veloppement de {F}riedrichs. {II}.
	\newblock {\em J. Math. Kyoto Univ.}, 26(2):157--175, 1986.
	
	\bibitem{KanoNishida79}
	T.~Kano and T.~Nishida.
	\newblock Sur les ondes de surface de l'eau avec une justification
	math\'{e}matique des \'{e}quations des ondes en eau peu profonde.
	\newblock {\em J. Math. Kyoto Univ.}, 19(2):335--370, 1979.
	
	\bibitem{KanoNishida84}
	T.~Kano and T.~Nishida.
	\newblock Water waves and {F}riedrichs expansion.
	\newblock In {\em Recent topics in nonlinear {PDE} ({H}iroshima, 1983)},
	volume~98 of {\em North-Holland Math. Stud.}, pages 39--57. North-Holland,
	Amsterdam, 1984.
	
	\bibitem{KaptsovMeleshkoSamatova20}
	E.~I. Kaptsov, S.~V. Meleshko, and N.~F. Samatova.
	\newblock The one-dimensional green–naghdi equations with a time dependent
	bottom topography and their conservation laws.
	\newblock {\em Phys. Fluids}, 32(12):123607, 2020.
	
	\bibitem{Kato72a}
	T.~Kato.
	\newblock Nonstationary flows of viscous and ideal fluids in {${\bf R}^{3}$}.
	\newblock {\em J. Functional Analysis}, 9:296--305, 1972.
	
	\bibitem{Kato95}
	T.~Kato.
	\newblock {\em Perturbation theory for linear operators}.
	\newblock Classics in Mathematics. Springer-Verlag, Berlin, 1995.
	\newblock Reprint of the 1980 edition.
	
	\bibitem{Kazakova18}
	M.~Kazakova.
	\newblock {\em Dispersive models of ocean waves propagation: Numerical issues
		and modelling}.
	\newblock PhD thesis, Universit{\'e} Toulouse 3 Paul Sabatier, 2018.
	
	\bibitem{Kelland40}
	P.~Kelland.
	\newblock On the theory of waves.
	\newblock {\em Trans. R. Soc. Edinburgh}, 14:497--545, 1840.
	
	\bibitem{KhorbatlyIsrawi}
	B.~Khorbatly and S.~Israwi.
	\newblock {Full justification for the extended Green-Naghdi system for an
		uneven bottom with surface tension}.
	\newblock Preprint hal-02994586.
	
	\bibitem{KhorbatlyLteifIsrawiEtAl}
	B.~Khorbatly, R.~Lteif, S.~Israwi, and S.~Gerbi.
	\newblock {Mathematical modeling and numerical analysis for the higher order
		Boussinesq system}.
	\newblock \arXiv{2102.08045}.
	
	\bibitem{KhorbatlyZaiterIsrwai18}
	B.~Khorbatly, I.~Zaiter, and S.~Isrwai.
	\newblock Derivation and well-posedness of the extended {G}reen-{N}aghdi
	equations for flat bottoms with surface tension.
	\newblock {\em J. Math. Phys.}, 59(7):071501, 20, 2018.
	
	\bibitem{Killworth97}
	P.~D. Killworth.
	\newblock On the parameterization of eddy transfer part i. theory.
	\newblock {\em J. Mar. Res.,}, 55(6):1171--1197, 1997.
	
	\bibitem{KimBaiErtekinEtAl01}
	J.~W. Kim, K.~J. Bai, R.~C. Ertekin, and W.~C. Webster.
	\newblock A derivation of the {G}reen-{N}aghdi equations for irrotational
	flows.
	\newblock {\em J. Engrg. Math.}, 40(1):17--42, 2001.
	
	\bibitem{Kirby04}
	J.~Kirby.
	\newblock Nonlinear ocean surface waves, 2004.
	
	\bibitem{KlainermanMajda81}
	S.~Klainerman and A.~Majda.
	\newblock Singular limits of quasilinear hyperbolic systems with large
	parameters and the incompressible limit of compressible fluids.
	\newblock {\em Comm. Pure Appl. Math.}, 34(4):481--524, 1981.
	
	\bibitem{KleinLinaresPilodEtAl18}
	C.~Klein, F.~Linares, D.~Pilod, and J.-C. Saut.
	\newblock On {W}hitham and related equations.
	\newblock {\em Stud. Appl. Math.}, 140(2):133--177, 2018.
	
	\bibitem{Klein95}
	R.~Klein.
	\newblock Semi-implicit extension of a {G}odunov-type scheme based on low
	{M}ach number asymptotics. {I}. {O}ne-dimensional flow.
	\newblock {\em J. Comput. Phys.}, 121(2):213--237, 1995.
	
	\bibitem{Klopman10}
	G.~Klopman.
	\newblock {\em Variational Boussinesq modelling of surface gravity waves over
		bathymetry}.
	\newblock PhD thesis, Univ. of Twente, 2010.
	
	\bibitem{KlopmanGroesenDingemans10}
	G.~Klopman, B.~van Groesen, and M.~W. Dingemans.
	\newblock A variational approach to {B}oussinesq modelling of fully nonlinear
	water waves.
	\newblock {\em J. Fluid Mech.}, 657:36--63, 2010.
	
	\bibitem{KortewegDe95}
	D.~J. Korteweg and G.~De~Vries.
	\newblock On the change of form of long waves advancing in a rectangular canal,
	and on a new type of long stationary waves.
	\newblock {\em Philos. Mag.}, 5(39):422--443, 1895.
	
	\bibitem{Kreisel49}
	G.~Kreisel.
	\newblock Surface waves.
	\newblock {\em Quart. Appl. Math.}, 7(1):21--44, 1949.
	
	\bibitem{KukavicaTemamVicolEtAl11}
	I.~Kukavica, R.~Temam, V.~C. Vicol, and M.~Ziane.
	\newblock Local existence and uniqueness for the hydrostatic {E}uler equations
	on a bounded domain.
	\newblock {\em J. Differential Equations}, 250(3):1719--1746, 2011.
	
	\bibitem{Lagrange81}
	J.~L. {\noopsort{Lagrange}}{de Lagrange}.
	\newblock Mémoire sur la théorie du mouvement des fluides.
	\newblock In {\em Œuvres complètes, tome 4}, pages 695--748. Nouveaux
	mémoires de l'Académie royale des sciences et belles-lettres de Berlin,
	1781.
	
	\bibitem{Lamb93}
	H.~Lamb.
	\newblock {\em Hydrodynamics}.
	\newblock Cambridge Mathematical Library. Cambridge University Press,
	Cambridge, reprint of the 1932 sixth edition. edition, 1993.
	
	\bibitem{Lannes13}
	D.~Lannes.
	\newblock A {S}tability {C}riterion for {T}wo-{F}luid {I}nterfaces and
	{A}pplications.
	\newblock {\em Arch. Ration. Mech. Anal.}, 208(2):481--567, 2013.
	
	\bibitem{Lannes}
	D.~Lannes.
	\newblock {\em The water waves problem}, volume 188 of {\em Mathematical
		Surveys and Monographs}.
	\newblock American Mathematical Society, Providence, RI, 2013.
	\newblock Mathematical analysis and asymptotics.
	
	\bibitem{Lannes20}
	D.~Lannes.
	\newblock Modeling shallow water waves.
	\newblock {\em Nonlinearity}, 33(5):R1--R57, 2020.
	
	\bibitem{LannesBonneton09}
	D.~Lannes and P.~Bonneton.
	\newblock Derivation of asymptotic two-dimensional time-dependent equations for
	surface water wave propagation.
	\newblock {\em Phys. Fluids}, 21(1):016601, 2009.
	
	\bibitem{LannesMarche15}
	D.~Lannes and F.~Marche.
	\newblock A new class of fully nonlinear and weakly dispersive {G}reen-{N}aghdi
	models for efficient 2{D} simulations.
	\newblock {\em J. Comput. Phys.}, 282:238--268, 2015.
	
	\bibitem{LannesMetivier18}
	D.~Lannes and G.~M\'{e}tivier.
	\newblock The shoreline problem for the one-dimensional shallow water and
	{G}reen-{N}aghdi equations.
	\newblock {\em J. \'{E}c. polytech. Math.}, 5:455--518, 2018.
	
	\bibitem{LannesMing15}
	D.~Lannes and M.~Ming.
	\newblock The {K}elvin-{H}elmholtz instabilities in two-fluids shallow water
	models.
	\newblock In {\em Hamiltonian partial differential equations and applications},
	volume~75 of {\em Fields Inst. Commun.}, pages 185--234. Fields Inst. Res.
	Math. Sci., Toronto, ON, 2015.
	
	\bibitem{Lavrentprimeev54}
	M.~A. Lavrent{\cprime}ev.
	\newblock I. {O}n the theory of long waves. {II}. {A} contribution to the
	theory of long waves.
	\newblock {\em Amer. Math. Soc. Translation}, 1954(102):53, 1954.
	
	\bibitem{LeGavrilyukHank10}
	O.~Le~M{\'e}tayer, S.~Gavrilyuk, and S.~Hank.
	\newblock A numerical scheme for the {G}reen-{N}aghdi model.
	\newblock {\em J. Comput. Phys.}, 229(6):2034--2045, 2010.
	
	\bibitem{LeTouze03}
	D.~Le~Touz{\'e}.
	\newblock {\em M{\'e}thodes spectrales pour la mod{\'e}lisation
		non-lin{\'e}aire d'{\'e}coulements a surface libre instationnaires}.
	\newblock PhD thesis, École Centrale de Nantess, 2003.
	
	\bibitem{Lebeau02}
	G.~Lebeau.
	\newblock R\'egularit\'e du probl\`eme de {K}elvin-{H}elmholtz pour
	l'\'equation d'{E}uler 2d.
	\newblock {\em ESAIM Control Optim. Calc. Var.}, 8:801--825 (electronic), 2002.
	\newblock A tribute to J. L. Lions.
	
	\bibitem{LevermoreOliverTiti96}
	C.~D. Levermore, M.~Oliver, and E.~S. Titi.
	\newblock Global well-posedness for models of shallow water in a basin with a
	varying bottom.
	\newblock {\em Indiana Univ. Math. J.}, 45(2):479--510, 1996.
	
	\bibitem{Levi-Civita25}
	T.~Levi-Civita.
	\newblock D\'{e}termination rigoureuse des ondes permanentes d'ampleur finie.
	\newblock {\em Math. Ann.}, 93(1):264--314, 1925.
	
	\bibitem{LiTiti18}
	J.~Li and E.~S. Titi.
	\newblock Recent advances concerning certain class of geophysical flows.
	\newblock In {\em Handbook of mathematical analysis in mechanics of viscous
		fluids}, pages 933--971. Springer, Cham, 2018.
	
	\bibitem{LiGuyenneLietal14}
	M.~Li, P.~Guyenne, F.~Li, and L.~Xu.
	\newblock High order well-balanced {CDG}-{FE} methods for shallow water waves
	by a {G}reen-{N}aghdi model.
	\newblock {\em J. Comput. Phys.}, 257(part A):169--192, 2014.
	
	\bibitem{Li01}
	Y.~A. Li.
	\newblock Linear stability of solitary waves of the {G}reen-{N}aghdi equations.
	\newblock {\em Comm. Pure Appl. Math.}, 54(5):501--536, 2001.
	
	\bibitem{Li02}
	Y.~A. Li.
	\newblock Hamiltonian structure and linear stability of solitary waves of the
	{G}reen-{N}aghdi equations.
	\newblock {\em J. Nonlinear Math. Phys.}, 9(suppl. 1):99--105, 2002.
	\newblock Recent advances in integrable systems (Kowloon, 2000).
	
	\bibitem{Li06}
	Y.~A. Li.
	\newblock A shallow-water approximation to the full water wave problem.
	\newblock {\em Comm. Pure Appl. Math.}, 59(9):1225--1285, 2006.
	
	\bibitem{LiHymanChoi04}
	Y.~A. Li, J.~M. Hyman, and W.~Choi.
	\newblock A numerical study of the exact evolution equations for surface waves
	in water of finite depth.
	\newblock {\em Stud. Appl. Math.}, 113(3):303--324, 2004.
	
	\bibitem{Lions84}
	P.-L. Lions.
	\newblock The concentration-compactness principle in the calculus of
	variations. {T}he locally compact case. {I}.
	\newblock {\em Ann. Inst. H. Poincar\'e Anal. Non Lin\'eaire}, 1(2):109--145,
	1984.
	
	\bibitem{LiskaWendroff97}
	R.~Liska and B.~Wendroff.
	\newblock Analysis and computation with stratified fluid models.
	\newblock {\em J. Comput. Phys.}, 137(1):212--244, 1997.
	
	\bibitem{LiuWang12}
	P.~L.-F. Liu and X.~Wang.
	\newblock A multi-layer model for nonlinear internal wave propagation in
	shallow water.
	\newblock {\em J. Fluid Mech.}, 695:341--365, 2012.
	
	\bibitem{Luke67}
	J.~C. Luke.
	\newblock A variational principle for a fluid with a free surface.
	\newblock {\em J. Fluid Mech.}, 27:395--397, 1967.
	
	\bibitem{LynettLiu04a}
	P.~Lynett and P.~L.-F. Liu.
	\newblock Linear analysis of the multi-layer model.
	\newblock {\em Coastal Engineering}, 51(5-6):439--454, 2004.
	
	\bibitem{LynettLiu04}
	P.~Lynett and P.~L.-F. Liu.
	\newblock A two-layer approach to wave modelling.
	\newblock {\em Proc. R. Soc. Lond. Ser. A Math. Phys. Eng. Sci.},
	460(2049):2637--2669, 2004.
	
	\bibitem{MadsenBinghamLiu02}
	P.~A. Madsen, H.~B. Bingham, and H.~Liu.
	\newblock A new {B}oussinesq method for fully nonlinear waves from shallow to
	deep water.
	\newblock {\em J. Fluid Mech.}, 462:1--30, 2002.
	
	\bibitem{MadsenBinghamSchaeffer03}
	P.~A. Madsen, H.~B. Bingham, and H.~A. Sch\"{a}ffer.
	\newblock Boussinesq-type formulations for fully nonlinear and extremely
	dispersive water waves: derivation and analysis.
	\newblock {\em R. Soc. Lond. Proc. Ser. A Math. Phys. Eng. Sci.},
	459(2033):1075--1104, 2003.
	
	\bibitem{MadsenFuhrman10}
	P.~A. Madsen and D.~R. Fuhrman.
	\newblock {High-order Boussinesq-type modelling of nonlinear wave phenomena in
		deep and shallow water}.
	\newblock In Q.~Ma, editor, {\em Advances in Numerical Simulation of Nonlinear
		Water Waves}, chapter~7, pages 245--285. World Scientific, 2010.
	
	\bibitem{MadsenSchaeffer99}
	P.~A. Madsen and H.~A. Schäffer.
	\newblock A review of boussinesq-type equations for gravity waves.
	\newblock {\em Advances in Coastal and Ocean Engineering}, 5:1--94, 1999.
	
	\bibitem{Makarenko86}
	N.~Makarenko.
	\newblock A second long-wave approximation in the cauchy-poisson problem.
	\newblock {\em Dyn. Contin. Media}, 77:56--72, 1986.
	\newblock In Russian.
	
	\bibitem{Malprimetseva89}
	Z.~L. Mal{\cprime}tseva.
	\newblock Unsteady long waves in a two-layer fluid.
	\newblock {\em Dinamika Sploshn. Sredy}, (93-94):96--110, 1989.
	
	\bibitem{MarsdenShkoller01}
	J.~E. Marsden and S.~Shkoller.
	\newblock Global well-posedness for the {L}agrangian averaged {N}avier-{S}tokes
	({LANS}-{$\alpha$}) equations on bounded domains.
	\newblock In {\em Topological methods in the physical sciences}, volume 359,
	pages 1449--1468. R. Soc. Lond. Philos. Trans. Ser. A Math. Phys. Eng. Sci.,
	2001.
	
	\bibitem{Masmoudi07}
	N.~Masmoudi.
	\newblock Examples of singular limits in hydrodynamics.
	\newblock In {\em Handbook of differential equations: evolutionary equations.
		{V}ol. {III}}, Handb. Differ. Equ., pages 195--275. Elsevier/North-Holland,
	Amsterdam, 2007.
	
	\bibitem{MasmoudiWong12}
	N.~Masmoudi and T.~K. Wong.
	\newblock On the {$H^s$} theory of hydrostatic {E}uler equations.
	\newblock {\em Arch. Ration. Mech. Anal.}, 204(1):231--271, 2012.
	
	\bibitem{Matsuno15}
	Y.~Matsuno.
	\newblock Hamiltonian formulation of the extended {G}reen-{N}aghdi equations.
	\newblock {\em Phys. D}, 301/302:1--7, 2015.
	
	\bibitem{Matsuno16}
	Y.~Matsuno.
	\newblock Hamiltonian structure for two-dimensional extended green–naghdi
	equations.
	\newblock {\em Proc. R. Soc. Lond. Ser. A Math. Phys. Eng. Sci.},
	472(2190):20160127--, 2016.
	
	\bibitem{MeiLe66}
	C.~C. Mei and B.~Le~M\'{e}haut\'{e}.
	\newblock Note on the equations of long waves over an uneven bottom.
	\newblock {\em J. Geophys. Res.}, 71:393--400, 1966.
	
	\bibitem{Meister00}
	A.~Meister.
	\newblock Asymptotic single and multiple scale expansions in the low {M}ach
	number limit.
	\newblock {\em SIAM J. Appl. Math.}, 60(1):256--271, 2000.
	
	\bibitem{Melinand}
	B.~Melinand.
	\newblock A mathematical study of meteo and landslide tsunamis: the {P}roudman
	resonance.
	\newblock {\em Nonlinearity}, 28(11):4037--4080, 2015.
	
		
	\bibitem{MelinandDuchene}
	B.~Melinand and V.~Duch\^ene.
	\newblock  Rectification of a deep water model for surface gravity waves.
	\newblock \arXiv{2203.03277}.
	
	\bibitem{Mesognon-Gireaua}
	B.~M\'esognon-Gireau.
	\newblock {The singular limit of the Water-Waves equations in the rigid lid
		regime}.
	\newblock \arXiv{1512.02424}.
	
	\bibitem{Mesognon-Gireau17a}
	B.~M\'esognon-Gireau.
	\newblock The {C}auchy problem on large time for a {B}oussinesq-{P}eregrine
	equation with large topography variations.
	\newblock {\em Adv. Differential Equations}, 22(7-8):457--504, 2017.
	
	\bibitem{Mesognon-Gireau17}
	B.~M\'esognon-Gireau.
	\newblock The {C}auchy problem on large time for the water waves equations with
	large topography variations.
	\newblock {\em Ann. Inst. H. Poincar\'e Anal. Non Lin\'eaire}, 34(1):89--118,
	2017.
	
	\bibitem{Mesognon-Gireau}
	B.~M\'{e}sognon-Gireau.
	\newblock A dispersive estimate for the linearized water-waves equations in
	finite depth.
	\newblock {\em J. Math. Fluid Mech.}, 19(3):469--500, 2017.
	
	\bibitem{Metivier08}
	G.~M{\'e}tivier.
	\newblock {\em Para-differential calculus and applications to the {C}auchy
		problem for nonlinear systems}, volume~5 of {\em Centro di Ricerca Matematica
		Ennio De Giorgi (CRM) Series}.
	\newblock Edizioni della Normale, Pisa, 2008.
	
	\bibitem{Metivier09}
	G.~M{\'e}tivier.
	\newblock The mathematics of nonlinear optics.
	\newblock In {\em Handbook of differential equations: evolutionary equations.
		{V}ol. {V}}, Handb. Differ. Equ., pages 169--313. Elsevier/North-Holland,
	Amsterdam, 2009.
	
	\bibitem{MetivierSchochet01}
	G.~M{\'e}tivier and S.~Schochet.
	\newblock The incompressible limit of the non-isentropic {E}uler equations.
	\newblock {\em Arch. Ration. Mech. Anal.}, 158(1):61--90, 2001.
	
	\bibitem{Mielke88}
	A.~Mielke.
	\newblock Reduction of quasilinear elliptic equations in cylindrical domains
	with applications.
	\newblock {\em Math. Methods Appl. Sci.}, 10(1):51--66, 1988.
	
	\bibitem{MilesSalmon85}
	J.~Miles and R.~Salmon.
	\newblock Weakly dispersive nonlinear gravity waves.
	\newblock {\em J. Fluid Mech.}, 157:519--531, 1985.
	
	\bibitem{Miles61}
	J.~W. Miles.
	\newblock On the stability of heterogeneous shear flows.
	\newblock {\em J. Fluid Mech.}, 10:496--508, 1961.
	
	\bibitem{Miles77}
	J.~W. Miles.
	\newblock On {H}amilton's principle for surface waves.
	\newblock {\em J. Fluid Mech.}, 83(1):153--158, 1977.
	
	\bibitem{MilewskiTabakTurnerEtAl04}
	P.~Milewski, E.~Tabak, C.~Turner, R.~Rosales, and F.~Menzaque.
	\newblock Nonlinear stability of two-layer flows.
	\newblock {\em Commun. Math. Sci.}, 2(3):427--442, 2004.
	
	\bibitem{MilewskiVanden-BroeckWang10}
	P.~A. Milewski, J.-M. Vanden-Broeck, and Z.~Wang.
	\newblock Dynamics of steep two-dimensional gravity--capillary solitary waves.
	\newblock {\em J. Fluid Mech.}, 664:466--477, 2010.
	
	\bibitem{MingWang20}
	M.~Ming and C.~Wang.
	\newblock Water waves problem with surface tension in a corner domain {I}: {A}
	priori estimates with constrained contact angle.
	\newblock {\em SIAM J. Math. Anal.}, 52(5):4861--4899, 2020.
	
	\bibitem{MingWang21}
	M.~Ming and C.~Wang.
	\newblock Water-waves problem with surface tension in a corner domain {II}: the
	local well-posedness.
	\newblock {\em Comm. Pure Appl. Math.}, 74(2):225--285, 2021.
	
	\bibitem{MitsotakisDutykhCarter17}
	D.~Mitsotakis, D.~Dutykh, and J.~Carter.
	\newblock On the nonlinear dynamics of the traveling-wave solutions of the
	{S}erre system.
	\newblock {\em Wave Motion}, 70:166--182, 2017.
	
	\bibitem{MitsotakisDutykhCarter17a}
	D.~Mitsotakis, D.~Dutykh, and J.~Carter.
	\newblock On the nonlinear dynamics of the traveling-wave solutions of the
	{S}erre system.
	\newblock {\em Wave Motion}, 70:166--182, 2017.
	
	\bibitem{MitsotakisSynolakisMcGuinness17}
	D.~Mitsotakis, C.~Synolakis, and M.~McGuinness.
	\newblock A modified {G}alerkin/finite element method for the numerical
	solution of the {S}erre-{G}reen-{N}aghdi system.
	\newblock {\em Internat. J. Numer. Methods Fluids}, 83(10):755--778, 2017.
	
	\bibitem{Miyata85}
	M.~Miyata.
	\newblock An internal solitary wave of large amplitude.
	\newblock {\em La mer}, 23(2):43--48, 1985.
	
	\bibitem{Miyata87}
	M.~Miyata.
	\newblock Long internal waves of large amplitude.
	\newblock In {\em Nonlinear Water Waves: IUTAM Symposium}, pages 399--405,
	Tokyo, Aug. 1987. Springer.
	
	\bibitem{MoldabayevKalischDutykh15}
	D.~Moldabayev, H.~Kalisch, and D.~Dutykh.
	\newblock The {W}hitham equation as a model for surface water waves.
	\newblock {\em Phys. D}, 309:99--107, 2015.
	
	\bibitem{MolinetTalhoukZaiter20}
	L.~Molinet, R.~Talhouk, and I.~Zaiter.
	\newblock The classical {B}oussinesq system revisited.
	\newblock \arXiv{2001.11870}.
	
	\bibitem{Monjarreta}
	R.~Monjarret.
	\newblock Local well-posedness of the multi-layer shallow-water model with free
	surface.
	\newblock \arXiv{1411.2342}.
	
	\bibitem{Monjarret15}
	R.~Monjarret.
	\newblock Local well-posedness of the two-layer shallow water model with free
	surface.
	\newblock {\em SIAM J. Appl. Math.}, 75(5):2311--2332, 2015.
	
	\bibitem{MontanelliBootland20}
	H.~Montanelli and N.~Bootland.
	\newblock Solving periodic semilinear stiff {PDE}s in 1{D}, 2{D} and 3{D} with
	exponential integrators.
	\newblock {\em Math. Comput. Simulation}, 178:307--327, 2020.
	
	\bibitem{Montgomery40}
	R.~Montgomery.
	\newblock The present evidence on the importance of lateral mixing processes in
	the ocean.
	\newblock {\em Bull. Amer. Meteor. Soc.}, 21(3):87--94, 1940.
	
	\bibitem{Munz02}
	C.-D. Munz.
	\newblock Computational fluid dynamics and aeroacoustics for low {M}ach number
	flow.
	\newblock In {\em Hyperbolic partial differential equations ({H}amburg, 2001)},
	pages 269--320. Friedr. Vieweg, Braunschweig, 2002.
	
	\bibitem{MunzRollerKleinEtAl03}
	C.-D. Munz, S.~Roller, R.~Klein, and K.~J. Geratz.
	\newblock The extension of incompressible flow solvers to the weakly
	compressible regime.
	\newblock {\em Comput. \& Fluids}, 32(2):173--196, 2003.
	
	\bibitem{MurakamiIguchi15}
	Y.~Murakami and T.~Iguchi.
	\newblock Solvability of the initial value problem to a model system for water
	waves.
	\newblock {\em Kodai Math. J.}, 38(2):470--491, 2015.
	
	\bibitem{Nachbin03}
	A.~Nachbin.
	\newblock A terrain-following {B}oussinesq system.
	\newblock {\em SIAM J. Appl. Math.}, 63(3):905--922, 2003.
	
	\bibitem{NachbinSolna03}
	A.~Nachbin and K.~S\o~lna.
	\newblock Apparent diffusion due to topographic microstructure in shallow
	waters.
	\newblock {\em Phys. Fluids}, 15(1):66--77, 2003.
	
	\bibitem{NakayamaKakinuma10}
	K.~Nakayama and T.~Kakinuma.
	\newblock Internal waves in a two-layer system using fully nonlinear
	internal-wave equations.
	\newblock {\em Internat. J. Numer. Methods Fluids}, 62(5):574--590, 2010.
	
	\bibitem{Nalimov74}
	V.~I. Nalimov.
	\newblock The {C}auchy-{P}oisson problem.
	\newblock {\em Dinamika Splo\v sn. Sredy}, (Vyp. 18 Dinamika Zidkost. so
	Svobod. Granicami):104--210, 254, 1974.
	
	\bibitem{Nekrasov51}
	A.~I. Nekrasov.
	\newblock {\em To\v{c}naya teoriya voln ustanoviv\v{s}egocya vida na
		poverhnosti tya\v{z}elo\u{\i} \v{z}idkosti}.
	\newblock Izdat. Akad. Nauk SSSR, Moscow, 1951.
	
	\bibitem{NemotoIguchi18}
	R.~Nemoto and T.~Iguchi.
	\newblock Solvability of the initial value problem to the {I}sobe-{K}akinuma
	model for water waves.
	\newblock {\em J. Math. Fluid Mech.}, 20(2):631--653, 2018.
	
	\bibitem{NguyenDias08}
	H.~Y. Nguyen and F.~Dias.
	\newblock A {B}oussinesq system for two-way propagation of interfacial waves.
	\newblock {\em Phys. D}, 237(18):2365--2389, 2008.
	
	\bibitem{Nicholls16}
	D.~P. Nicholls.
	\newblock High-order perturbation of surfaces short course: boundary value
	problems.
	\newblock In {\em Lectures on the theory of water waves}, volume 426 of {\em
		London Math. Soc. Lecture Note Ser.}, pages 1--18. Cambridge Univ. Press,
	Cambridge, 2016.
	
	\bibitem{NilssonWang19}
	D.~Nilsson and Y.~Wang.
	\newblock Solitary wave solutions to a class of {W}hitham-{B}oussinesq systems.
	\newblock {\em Z. Angew. Math. Phys.}, 70(3):Paper No. 70, 13, 2019.
	
	\bibitem{Nwogu93}
	O.~Nwogu.
	\newblock Alternative form of boussinesq equations for nearshore wave
	propagation.
	\newblock {\em Journal of waterway, port, coastal, and ocean engineering},
	119(6):618--638, 1993.
	
	\bibitem{Oliver97}
	M.~Oliver.
	\newblock Classical solutions for a generalized {E}uler equation in two
	dimensions.
	\newblock {\em J. Math. Anal. Appl.}, 215(2):471--484, 1997.
	
	\bibitem{Olver84}
	P.~J. Olver.
	\newblock Hamiltonian and non-{H}amiltonian models for water waves.
	\newblock In {\em Trends and applications of pure mathematics to mechanics
		({P}alaiseau, 1983)}, volume 195 of {\em Lecture Notes in Phys.}, pages
	273--290. Springer, Berlin, 1984.
	
	\bibitem{Orszag71}
	S.~A. Orszag.
	\newblock On the elimination of aliasing in finite-difference schemes by
	filtering high-wavenumber components.
	\newblock {\em J. Atmos. Sci.}, 28:1074, 1971.
	
	\bibitem{Orszag72}
	S.~A. Orszag.
	\newblock Comparison of pseudospectral and spectral approximation.
	\newblock {\em Stud. Appl. Math.}, 51:253--259, 1972.
	
	\bibitem{Ovsjannikov74}
	L.~V. Ovsjannikov.
	\newblock To the shallow water theory foundation.
	\newblock {\em Arch. Mech. (Arch. Mech. Stos.)}, 26:407--422, 1974.
	
	\bibitem{Ovsjannikov76}
	L.~V. Ovsjannikov.
	\newblock Cauchy problem in a scale of {B}anach spaces and its application to
	the shallow water theory justification.
	\newblock In {\em Applications of methods of functional analysis to problems in
		mechanics ({J}oint {S}ympos., {IUTAM}/{IMU}, {M}arseille, 1975)}, pages
	426--437. Lecture Notes in Math., 503. 1976.
	
	\bibitem{Ovsjannikov79}
	L.~V. Ovsjannikov.
	\newblock Models of two-layered ``shallow water''.
	\newblock {\em Zh. Prikl. Mekh. i Tekhn. Fiz.}, (2):3--14, 180, 1979.
	
	\bibitem{Pandey19}
	A.~K. Pandey.
	\newblock The effects of surface tension on modulational instability in
	full-dispersion water-wave models.
	\newblock {\em Eur. J. Mech. B Fluids}, 77:177--182, 2019.
	
	\bibitem{Papoutsellis17}
	C.~Papoutsellis.
	\newblock {\em Nonlinear water waves over varying bathymetry: theoretical and
		numerical study using variational methods}.
	\newblock PhD thesis, National Technical University of Athens, 2017.
	
	\bibitem{PapoutsellisAthanassoulis17}
	C.~E. Papoutsellis and G.~A. Athanassoulis.
	\newblock A new efficient hamiltonian approach to the nonlinear water-wave
	problem over arbitrary bathymetry.
	\newblock \arXiv{1704.03276}.
	
	\bibitem{PeiWang19}
	L.~Pei and Y.~Wang.
	\newblock A note on well-posedness of bidirectional {W}hitham equation.
	\newblock {\em Appl. Math. Lett.}, 98:215--223, 2019.
	
	\bibitem{PetcuTemamZiane09}
	M.~Petcu, R.~M. Temam, and M.~Ziane.
	\newblock Some mathematical problems in geophysical fluid dynamics.
	\newblock In {\em Handbook of numerical analysis. {V}ol. {XIV}. {S}pecial
		volume: computational methods for the atmosphere and the oceans}, volume~14
	of {\em Handb. Numer. Anal.}, pages 577--750. Elsevier/North-Holland,
	Amsterdam, 2009.
	
	\bibitem{PittZoppouRoberts18}
	J.~P.~A. Pitt, C.~Zoppou, and S.~G. Roberts.
	\newblock Behaviour of the {S}erre equations in the presence of steep gradients
	revisited.
	\newblock {\em Wave Motion}, 76:61--77, 2018.
	
	\bibitem{Plotnikov02}
	P.~I. Plotnikov.
	\newblock Proof of the {S}tokes conjecture in the theory of surface waves.
	\newblock {\em Stud. Appl. Math.}, 108(2):217--244, 2002.
	\newblock Translated from Dinamika Sploshn. Sredy No. 57 (1982), 41--76 [
	MR0752600 (85f:76036)].
	
	\bibitem{PlotnikovToland01}
	P.~I. Plotnikov and J.~F. Toland.
	\newblock Nash-{M}oser theory for standing water waves.
	\newblock {\em Arch. Ration. Mech. Anal.}, 159(1):1--83, 2001.
	
	\bibitem{dePoyferre19}
	T.~{\noopsort{Poyferré}}{de Poyferré}.
	\newblock A priori estimates for water waves with emerging bottom.
	\newblock {\em Arch. Ration. Mech. Anal.}, 232(2):763--812, 2019.
	
	\bibitem{Rayleigh76}
	J.~W.~S. Rayleigh.
	\newblock On waves.
	\newblock {\em Philos. Mag.}, 1(5):251--271, 1876.
	
	\bibitem{Renardy09}
	M.~Renardy.
	\newblock Ill-posedness of the hydrostatic {E}uler and {N}avier-{S}tokes
	equations.
	\newblock {\em Arch. Ration. Mech. Anal.}, 194(3):877--886, 2009.
	
	\bibitem{Ripa91}
	P.~Ripa.
	\newblock General stability conditions for a multi-layer model.
	\newblock {\em J. Fluid Mech.}, 222:119--137, 1991.
	
	\bibitem{RosalesPapanicolaou83}
	R.~R. Rosales and G.~C. Papanicolaou.
	\newblock Gravity waves in a channel with a rough bottom.
	\newblock {\em Stud. Appl. Math.}, 68(2):89--102, 1983.
	
	\bibitem{RuizdeZarateVigoNachbinEtAl09}
	A.~Ruiz~de Z\'{a}rate, D.~G.~A. Vigo, A.~Nachbin, and W.~Choi.
	\newblock A higher-order internal wave model accounting for large bathymetric
	variations.
	\newblock {\em Stud. Appl. Math.}, 122(3):275--294, 2009.
	
	\bibitem{SaadSchultz86}
	Y.~Saad and M.~H. Schultz.
	\newblock G{MRES}: a generalized minimal residual algorithm for solving
	nonsymmetric linear systems.
	\newblock {\em SIAM J. Sci. Statist. Comput.}, 7(3):856--869, 1986.
	
	\bibitem{Said}
	A.~R. Said.
	\newblock A geometric proof of the quasi-linearity of the water-waves system.
	\newblock \arXiv{2002.02940}.
	
	\bibitem{Saint-Venant71}
	B.~{\noopsort{Saint-Venant}}{de Saint-Venant}.
	\newblock Th{\'e}orie du mouvement non-permanent des eaux, avec application aux
	crues des rivi{\`e}res et {\`a} l'introduction des mar{\'e}es dans leur lit.
	\newblock {\em C.R. Acad. Sci. Paris}, 73:147--154, 1871.
	
	\bibitem{Salmon88}
	R.~Salmon.
	\newblock Hamiltonian fluid mechanics.
	\newblock {\em Annual Review of Fluid Mechanics}, 20(1):225--256, 1988.
	
	\bibitem{SanchezFernandez-NietoLunaEtAl21}
	C.~E. S\'{a}nchez, E.~D. Fern\'{a}ndez-Nieto, T.~M. de~Luna, Y.~Penel, and
	J.~Sainte-Marie.
	\newblock Numerical {S}imulations of a {D}ispersive {M}odel {A}pproximating
	{F}ree-{S}urface {E}uler {E}quations.
	\newblock {\em J. Sci. Comput.}, 89(3):Paper No. 55, 2021.
	
	\bibitem{Saut13}
	J.-C. Saut.
	\newblock {\em Asymptotic models for surface and internal waves}.
	\newblock Publica\c{c}\~{o}es Matem\'{a}ticas do IMPA. [IMPA Mathematical
	Publications]. Instituto Nacional de Matem\'{a}tica Pura e Aplicada (IMPA),
	Rio de Janeiro, 2013.
	\newblock 29${\sp{{}}{\rm{o}}}$ Col\'{o}quio Brasileiro de Matem\'{a}tica.
	[29th Brazilian Mathematics Colloquium].
	
	\bibitem{SautWang}
	J.-C. Saut and Y.~Wang.
	\newblock {The wave breaking for Whitham-type equations revisited}.
	\newblock \arXiv{2006.03803}.
	
	\bibitem{SautXu}
	J.-C. Saut and L.~Xu.
	\newblock Long time existence for a two-dimensional strongly dispersive
	{B}oussinesq system.
	\newblock {\em Comm. Partial Differential Equations}, 2021.
	
	\bibitem{Schaeffer08}
	H.~A. Sch{\"a}ffer.
	\newblock Comparison of dirichlet--neumann operator expansions for nonlinear
	surface gravity waves.
	\newblock {\em Coastal Engineering}, 55(4):288--294, 2008.
	
	\bibitem{SchneiderBottaGeratzEtAl99}
	T.~Schneider, N.~Botta, K.~J. Geratz, and R.~Klein.
	\newblock Extension of finite volume compressible flow solvers to
	multi-dimensional, variable density zero {M}ach number flows.
	\newblock {\em J. Comput. Phys.}, 155(2):248--286, 1999.
	
	\bibitem{Schochet86a}
	S.~Schochet.
	\newblock Symmetric hyperbolic systems with a large parameter.
	\newblock {\em Comm. Partial Differential Equations}, 11(15):1627--1651, 1986.
	
	\bibitem{Schochet05}
	S.~Schochet.
	\newblock The mathematical theory of low {M}ach number flows.
	\newblock {\em M2AN Math. Model. Numer. Anal.}, 39(3):441--458, 2005.
	
	\bibitem{Schonbek81}
	M.~E. Schonbek.
	\newblock Existence of solutions for the {B}oussinesq system of equations.
	\newblock {\em J. Differential Equations}, 42(3):325--352, 1981.
	
	\bibitem{Schwartz74}
	L.~W. Schwartz.
	\newblock Computer extension and analytic continuation of {S}tokes' expansion
	for gravity waves.
	\newblock {\em J. Fluid Mech.}, 62(3):553--578, 1974.
	
	\bibitem{Seabra-SantosRenouardTemperville87}
	F.~J. Seabra-Santos, D.~P. Renouard, and A.~M. Temperville.
	\newblock Numerical and experimental study of the transformation of a solitary
	wave over a shelf or isolated obstacle.
	\newblock {\em J. Fluid Mech.}, 176:117--134, 3 1987.
	
	\bibitem{Serre53}
	F.~Serre.
	\newblock Contribution {\`a} l'{\'e}tude des {\'e}coulements permanents et
	variables dans les canaux.
	\newblock {\em La Houille Blanche}, (6):830--872, 1953.
	
	\bibitem{ShatahZeng08}
	J.~Shatah and C.~Zeng.
	\newblock A priori estimates for fluid interface problems.
	\newblock {\em Comm. Pure Appl. Math.}, 61(6):848--876, 2008.
	
	\bibitem{ShatahZeng11}
	J.~Shatah and C.~Zeng.
	\newblock Local well-posedness for fluid interface problems.
	\newblock {\em Arch. Ration. Mech. Anal.}, 199(2):653--705, 2011.
	
	\bibitem{Shepherd90}
	T.~G. Shepherd.
	\newblock {S}ymmetries, conservation laws, and {H}amiltonian structure in
	geophysical fluid dynamics.
	\newblock {\em Advances in Geophysics}, 32:287--338, 1990.
	
	\bibitem{ShkollerSideris19}
	S.~Shkoller and T.~C. Sideris.
	\newblock Global existence of near-affine solutions to the compressible {E}uler
	equations.
	\newblock {\em Arch. Ration. Mech. Anal.}, 234(1):115--180, 2019.
	
	\bibitem{Simon05}
	B.~Simon.
	\newblock Sturm oscillation and comparison theorems.
	\newblock In {\em Sturm-{L}iouville theory}, pages 29--43. Birkh\"{a}user,
	Basel, 2005.
	
	\bibitem{StefanovWright20}
	A.~Stefanov and J.~D. Wright.
	\newblock Small amplitude traveling waves in the full-dispersion {W}hitham
	equation.
	\newblock {\em J. Dynam. Differential Equations}, 32(1):85--99, 2020.
	
	\bibitem{Stein93}
	E.~M. Stein.
	\newblock {\em Harmonic analysis: real-variable methods, orthogonality, and
		oscillatory integrals}, volume~43 of {\em Princeton Mathematical Series}.
	\newblock Princeton University Press, Princeton, NJ, 1993.
	\newblock With the assistance of Timothy S. Murphy, Monographs in Harmonic
	Analysis, III.
	
	\bibitem{StewartDellar13}
	A.~L. Stewart and P.~J. Dellar.
	\newblock Multilayer shallow water equations with complete coriolis force. part
	3. hyperbolicity and stability under shear.
	\newblock {\em J. Fluid Mech.}, 723:289--317, 5 2013.
	
	\bibitem{Stoker48}
	J.~J. Stoker.
	\newblock The formation of breakers and bores. {T}he theory of nonlinear wave
	propagation in shallow water and open channels.
	\newblock {\em Communications on Appl. Math.}, 1:1--87, 1948.
	
	\bibitem{Stokes47}
	G.~G. Stokes.
	\newblock On the theory of oscillatory waves.
	\newblock {\em Trans. Cambridge Philos. Soc.}, 8:441--455, 1847.
	
	\bibitem{Struik26}
	D.~J. Struik.
	\newblock D\'{e}termination rigoureuse des ondes irrotationelles
	p\'{e}riodiques dans un canal \`a profondeur finie.
	\newblock {\em Math. Ann.}, 95(1):595--634, 1926.
	
	\bibitem{Sturm34}
	C.~Sturm.
	\newblock Mémoire sur les équations différentielles linéaires du second
	ordre.
	\newblock {\em J. Math. Pures Appl. (1)}, 1:106--186, 1834.
	
	\bibitem{SuGardner69}
	C.~H. Su and C.~S. Gardner.
	\newblock Korteweg-de {V}ries equation and generalizations. {III}. {D}erivation
	of the {K}orteweg-de {V}ries equation and {B}urgers equation.
	\newblock {\em J. Mathematical Phys.}, 10:536--539, 1969.
	
	\bibitem{Tao06}
	T.~Tao.
	\newblock {\em Nonlinear dispersive equations}, volume 106 of {\em CBMS
		Regional Conference Series in Mathematics}.
	\newblock Published for the Conference Board of the Mathematical Sciences,
	Washington, DC; by the American Mathematical Society, Providence, RI, 2006.
	\newblock Local and global analysis.
	
	\bibitem{TaoBlog}
	T.~Tao.
	\newblock Incompressible fluid equations, 2018.
	\newblock Blog posts available at
	\url{https://terrytao.wordpress.com/category/teaching/254a-incompressible-fluid-equations/}.
	
	\bibitem{TaylorIII}
	M.~E. Taylor.
	\newblock {\em Partial differential equations. {III} {N}onlinear equations},
	volume 117 of {\em Applied Mathematical Sciences}.
	\newblock Springer-Verlag, New York, 1997.
	
	\bibitem{Toland78}
	J.~F. Toland.
	\newblock On the existence of a wave of greatest height and {S}tokes's
	conjecture.
	\newblock {\em Proc. Roy. Soc. London Ser. A}, 363(1715):469--485, 1978.
	
	\bibitem{Trefethen00}
	L.~N. Trefethen.
	\newblock {\em Spectral methods in {MATLAB}}, volume~10 of {\em Software,
		Environments, and Tools}.
	\newblock Society for Industrial and Applied Mathematics (SIAM), Philadelphia,
	PA, 2000.
	
	\bibitem{TruongWahlenWheeler}
	T.~Truong, E.~Wahl\'{e}n, and M.~H. Wheeler.
	\newblock Global bifurcation of solitary waves for the {W}hitham equation.
	\newblock {\em Math. Ann.}, 2021.
	
	\bibitem{Ukai86}
	S.~Ukai.
	\newblock The incompressible limit and the initial layer of the compressible
	{E}uler equation.
	\newblock {\em J. Math. Kyoto Univ.}, 26(2):323--331, 1986.
	
	\bibitem{Vargas-MaganaPanayotaros16}
	R.~M. Vargas-Maga\~{n}a and P.~Panayotaros.
	\newblock A {W}hitham-{B}oussinesq long-wave model for variable topography.
	\newblock {\em Wave Motion}, 65:156--174, 2016.
	
	\bibitem{Vargas-MaganaPanayotarosMinzoni19}
	R.~M. Vargas-Maga\~{n}a, P.~Panayotaros, and A.~A. Minzoni.
	\newblock Linear modes for channels of constant cross-section and approximate
	{D}irichlet--{N}eumann operators.
	\newblock {\em Water Waves}, 1:343--370, 2019.
	
	\bibitem{ViottiDutykhDias14}
	C.~Viotti, D.~Dutykh, and F.~Dias.
	\newblock The conformal-mapping method for surface gravity waves in the
	presence of variable bathymetry and mean current.
	\newblock {\em Procedia IUTAM}, 11:110--118, 2014.
	
	\bibitem{VirissimoMilewski19}
	F.~d.~M. Vir\'{\i}ssimo and P.~A. Milewski.
	\newblock Three-layer flows in the shallow water limit.
	\newblock {\em Stud. Appl. Math.}, 142(4):487--512, 2019.
	
	\bibitem{VirissimoMilewski20}
	F.~d.~M. Vir\'{\i}ssimo and P.~A. Milewski.
	\newblock Nonlinear stability of two-layer shallow water flows with a free
	surface.
	\newblock {\em Proc. A.}, 476(2236):20190594, 20, 2020.
	
	\bibitem{Wang20}
	Y.~Wang.
	\newblock Well-posedness to the cauchy problem of a fully dispersive boussinesq
	system.
	\newblock {\em J. Dynam. Differential Equations}, pages 1--12, 2020.
	
	\bibitem{WeiKirbyGrilliEtAl95}
	G.~Wei, J.~T. Kirby, S.~T. Grilli, and R.~Subramanya.
	\newblock A fully nonlinear {B}oussinesq model for surface waves. {I}. {H}ighly
	nonlinear unsteady waves.
	\newblock {\em J. Fluid Mech.}, 294:71--92, 1995.
	
	\bibitem{WestBruecknerJandaEtAl87}
	B.~J. West, K.~A. Brueckner, R.~S. Janda, D.~M. Milder, and R.~L. Milton.
	\newblock A new numerical method for surface hydrodynamics.
	\newblock {\em J. Geophys. Res.}, 92:11803--11824, 1987.
	
	\bibitem{Whitham67}
	G.~B. Whitham.
	\newblock Variational methods and applications to water waves.
	\newblock {\em Proc. R. Soc. Lond. Ser. A Math. Phys. Eng. Sci.}, 299, 06 1967.
	
	\bibitem{Whitham}
	G.~B. Whitham.
	\newblock {\em Linear and nonlinear waves}.
	\newblock Pure and Applied Mathematics (New York). John Wiley \& Sons Inc., New
	York, 1999.
	\newblock Reprint of the 1974 original, A Wiley-Interscience Publication.
	
	\bibitem{WilkeningVasan15}
	J.~Wilkening and V.~Vasan.
	\newblock Comparison of five methods of computing the {D}irichlet-{N}eumann
	operator for the water wave problem.
	\newblock In {\em Nonlinear wave equations: analytic and computational
		techniques}, volume 635 of {\em Contemp. Math.}, pages 175--210. Amer. Math.
	Soc., Providence, RI, 2015.
	
	\bibitem{Wu}
	S.~Wu.
	\newblock The quartic integrability and long time existence of steep water
	waves in 2{D}.
	\newblock \arXiv{2010.09117}.
	
	\bibitem{Wu97}
	S.~Wu.
	\newblock Well-posedness in {S}obolev spaces of the full water wave problem in
	{$2$}-{D}.
	\newblock {\em Invent. Math.}, 130(1):39--72, 1997.
	
	\bibitem{Wu06}
	S.~Wu.
	\newblock Mathematical analysis of vortex sheets.
	\newblock {\em Comm. Pure Appl. Math.}, 59(8):1065--1206, 2006.
	
	\bibitem{Wu01}
	T.~Y. Wu.
	\newblock A unified theory for modeling water waves.
	\newblock {\em Adv. in Appl. Mech.}, 37:1--88, 2001.
	
	\bibitem{Xu10}
	Z.~Xu.
	\newblock {\em Asymptotic analysis and numerical analysis of the
		{B}enjamin-{O}no equation}.
	\newblock PhD thesis, University of Michigan, 2010.
	
	\bibitem{YangFangTangEtAl10}
	Y.~J. Yang, Y.~C. Fang, T.~Y. Tang, and S.~R. Ramp.
	\newblock Convex and concave types of second baroclinic mode internal solitary
	waves.
	\newblock {\em Nonlin. Processes Geophys.}, 17(6):605--614, 2010.
	
	\bibitem{YatesBenoit15}
	M.~L. Yates and M.~Benoit.
	\newblock Accuracy and efficiency of two numerical methods of solving the
	potential flow problem for highly nonlinear and dispersive water waves.
	\newblock {\em Internat. J. Numer. Methods Fluids}, 77(10):616--640, 2015.
	
	\bibitem{Yosihara82}
	H.~Yosihara.
	\newblock Gravity waves on the free surface of an incompressible perfect fluid
	of finite depth.
	\newblock {\em Publ. Res. Inst. Math. Sci.}, 18(1):49--96, 1982.
	
	\bibitem{YuHoward12}
	J.~Yu and L.~N. Howard.
	\newblock Exact {F}loquet theory for waves over arbitrary periodic
	topographies.
	\newblock {\em J. Fluid Mech.}, 712:451--470, 2012.
	
	\bibitem{Zakharov68}
	V.~E. Zakharov.
	\newblock Stability of periodic waves of finite amplitude on the surface of a
	deep fluid.
	\newblock {\em J. Appl. Mech. Tech. Phys.}, 9:190--194, 1968.
	
	\bibitem{Zakharov80}
	V.~E. Zakharov.
	\newblock Benney equations and quasiclassical approximation in the inverse
	problem method.
	\newblock {\em Funktsional. Anal. i Prilozhen.}, 14(2):15--24, 1980.
	
	\bibitem{ZakharovDyachenkoVasilyev02}
	V.~E. Zakharov, A.~I. Dyachenko, and O.~A. Vasilyev.
	\newblock New method for numerical simulation of a nonstationary potential flow
	of incompressible fluid with a free surface.
	\newblock {\em Eur. J. Mech. B Fluids}, 21(3):283--291, 2002.
	
	\bibitem{ZouFangLiu10}
	Z.~L. Zou, K.~Z. Fang, and Z.~B. Liu.
	\newblock {Inter-comparisons of different forms of higher-order Boussinesq
		equations}.
	\newblock In Q.~Ma, editor, {\em Advances in Numerical Simulation of Nonlinear
		Water Waves}, chapter~8, pages 287--323. World Scientific, 2010.
	
\end{thebibliography}
\end{document}